\documentclass[11pt]{amsart}
\usepackage{amsmath,amsfonts}
\usepackage[all]{xy}

\newcommand{\Alb}{\operatorname{Alb}}
\newcommand{\bul}{\bullet}

\newcommand{\Hilb}{\operatorname{Hilb}}
\newcommand{\Hifu}{\underline{\operatorname{Hilb}}}
\newcommand{\Hom}{\operatorname{Hom}}

\newcommand{\LPN}{L_{\underline{t}}(\underline{n}, \underline{\zeta})%
                 (\np^1 \times E)}
\renewcommand{\mod}{\operatorname{mod}}

\newcommand{\nund}{{\underline{n}}}
\newcommand{\pibar}{{\overline{\pi}}}
\newcommand{\Pic}{\operatorname{Pic}}

\newcommand{\rk}{\operatorname{rk}}
\newcommand{\shat}{\hat{s}} 

\newcommand{\rhotil}{\tilde{\rho}}

\newcommand{\Spec}{\operatorname{Spec}}

\newcommand{\tund}{{\underline{t}}}

\newcommand{\ext}{\operatorname{Ext}}
\newcommand{\der}{\operatorname{Der}}

\newcommand{\CC}{\mathcal{C}}
\newcommand{\EE}{\mathcal{E}}
\newcommand{\ee}{\EE^\bul} 

\newcommand{\FF}{\mathcal{F}}
\newcommand{\ff}{\FF^\bul} 
\newcommand{\GG}{\mathcal{G}}
\newcommand{\HH}{\mathcal{H}}
\newcommand{\JJ}{\mathcal{J}} 
\newcommand{\KK}{\mathcal{K}} 
\newcommand{\LL}{\mathcal{L}} 
\newcommand{\MM}{\mathcal{M}}

\newcommand{\OO}{\mathcal{O}} 

\newcommand{\nc}{{\mathbb{C}}}  
\newcommand{\nd}{{\mathbb{D}}}  
\newcommand{\ndtil}{{\widetilde{\mathbb{D}}}}
\newcommand{\ndhat}{{\hat{\mathbb{D}}}}
\renewcommand{\ne}{\mathbb{E}}
\newcommand{\nf}{\mathbb{F}}
\newcommand{\nh}{{\mathbb{H}}}  
\newcommand{\nhtil}{{\widetilde{\nh}}}
\newcommand{\nl}{{\mathbb{L}}}  
\newcommand{\nn}{{\mathbb{N}}}  
\newcommand{\np}{{\mathbb{P}}}  
\newcommand{\nq}{{\mathbb{Q}}}  
\newcommand{\nr}{{\mathbb{R}}}  
\newcommand{\nw}{{\mathbb{W}}}  

\newcommand{\nz}{{\mathbb{Z}}}  
\newcommand{\cL}{\mathcal{L}} 
\newcommand{\cT}{\mathcal{T}}
\newcommand{\cm}{\mathcal{M}}

\newcommand{\KG}{\mathfrak{K}} 
\newcommand{\LG}{\mathfrak{L}^\bullet} 

\newcommand{\nutil}{\tilde{\nu}}
\newcommand{\eps}{{\varepsilon}}

\newcommand{\pitil}{\tilde{\pi}}
\newcommand{\pihat}{\hat{\pi}}
\newcommand{\xibar}{{\bar{\xi}}}
\newcommand{\zetaund}{{\underline{\zeta}}}
\newcommand{\wbar}{{\bar W}}
\newcommand{\tbar}{{\bar T}}
\newcommand{\fbar}{{\bar f}}
\newcommand{\xbar}{{\bar X}}
\newcommand{\ybar}{{\bar Y}}
\newcommand{\omicron}{{\scriptstyle{\mathcal{O}}}}

\newcommand{\lb}{[[}
\newcommand{\rb}{]]}

	\newcommand{\ov}{\mathcal{O}_V}
	\newcommand{\ox}{\mathcal{O}_X}
	\newcommand{\co}{\mathcal{O}}
	\newcommand{\kv}{\mathcal{K}_V}
	\newcommand{\fp}[1]{\Hilb^{#1}_V \times_{Pic^{#1}_V} \Hilb^{k-#1}_V}

        \newcommand{\spin}{Spin^c(4)}
	\newcommand{\spiv}{Spin^c(V)}
	\newcommand{\spim}{Spin^c(M)}
	\newcommand{\oo}{{\scriptstyle{\mathcal{O}}}}
	\newcommand{\oohat}{\hat{\oo}}
	\newcommand{\ooo}{{\scriptscriptstyle{\mathcal{O}}}}
	\newcommand{\ooohat}{{\scriptscriptstyle{\hat{\mathcal{O}}}}}
	\newcommand{\oov}{({\scriptstyle{\mathcal{O}}}_1,{\bf H}_0)}
	\newcommand{\oovhat}{({\scriptstyle{\hat{\mathcal{O}}}}_1,\hat{{\bf H}}_0)}
	\newcommand{\ooov}{({\scriptscriptstyle{\mathcal{O}}}_1,{\bf H}_0)}
	\newcommand{\ooovhat}{({\scriptscriptstyle{\hat{\mathcal{O}}}}_1,\hat{{\bf H}}_0)}
	
	\newcommand{\swv}[1]{SW_{V,\ooo}(#1)}
	\newcommand{\swvhat}[1]{SW_{\hat{V},\ooohat}(#1)}
	\newcommand{\sww}[2]{SW^{#1}_{X,\ooov}(#2)}
	\newcommand{\swwv}[2]{SW^{#1}_{V,\ooov}(#2)}
	\newcommand{\swwvhat}[2]{SW^{#1}_{\hat{V},\ooovhat}(#2)}
	\newcommand{\swwm}[2]{SW^{#1}_{M,\ooov}(#2)}

	\renewcommand{\dfrac}[2]{\genfrac{}{}{0pt}{}{#1}{#2}}

	\newcommand{\INC}[2]{\mathfrak{C}_{#1 / #2 }}
	\newcommand{\im}{\operatorname{Im}}
	\newcommand{\coker}{\operatorname{coker}}

\newtheorem{thm}{Theorem}[section]
\newtheorem{lem}[thm]{Lemma}
\newtheorem{prop}[thm]{Proposition}
\newtheorem{cor}[thm]{Corollary}
\newtheorem{conj}[thm]{Conjecture}

\theoremstyle{definition}

\newtheorem{dfn}[thm]{Definition}
\newtheorem{nota}[thm]{Notation}

\newtheorem{rem}[thm]{Remark} 
 
\newtheorem{exa}[thm]{Example} 

\title{Poincar\'e Invariants}
\author[M.~D\"urr]{Markus D\"urr$^*$}
\author[A.~Kabanov]{Alexandre Kabanov$^\dagger$}
\author[Ch.~Okonek]{Christian Okonek$^*$}
\date{\today}
\thanks{$^*$Partially supported by: EAGER -- European Algebraic Geometry Research Training Network, contract No. HPRN-CT-2000-00099 (BBW 99.0030), and by SNF, nr. 2000-055290.98/1}
\thanks{$^\dagger$Partially supported by NSF grant number DMS-9803553.}
\begin{document}
\begin{abstract}
We construct an obstruction theory for relative Hilbert schemes in the sense of \cite{bf} and compute it explicitly for relative Hilbert schemes of  divisors on smooth projective varieties. In the special case of curves on a surface $V$, our obstruction theory determines a virtual fundamental class $\lb \Hilb^m_V \rb \in A_{\frac{m(m-k)}{2}}(\Hilb^m_V)$, which we use to define Poincar\'e invariants
\[
(P^+_V,P^-_V): H^2(V,\nz) \longrightarrow \Lambda^* H^1(V,\nz ) \times \Lambda^* H^1(V,\nz ).
\]
These maps are invariant under deformations, satisfy a blow-up formula, and a wall crossing formula for surfaces with $p_g(V)=0$. In the case $q(V) \geq 1$, we calculate the wall crossing formula explicitely in terms of fundamental classes of certain Brill-Noether loci for curves. We determine the invariants completely for ruled surfaces, and rederive from this classical results of Nagata and Lange. The invariant $(P^+_V,P^-_V)$ of an elliptic fibration is computed in terms of its multiple fibers.

When the fibered product $\Hilb^m_V \times_{\Pic^m_V}\Hilb^{k-m}_V$ is empty, there exists a more geometric obstruction theory, which gives rise to a second virtual fundamental class $\{ \Hilb^m_V\} \in A_{\frac{m(m-k)}{2}+p_g(V)}(\Hilb^m_V)$. We show that $\{ \Hilb^m_V \} = \lb \Hilb^m_V \rb$ when $p_g(V)=0$, and use the second obstruction theory to prove that $\lb \Hilb^m_V \rb=0$ when $p_g(V) > 0$ and $\Hilb^m_V \times_{\Pic^m_V}\Hilb^{k-m}=\emptyset$.

We conjecture that our Poincar\'e invariants coincide with the full Seiberg-Witten invariants of \cite{ot1} computed with respect to the canonical orientation data. The main evidence for this conjecture is based on the existence of an Kobayashi-Hitchin isomorphism which identifies the moduli spaces of monopoles with the corresponding Hilbert schemes. We expect that this isomorphism identifies also the corresponding virtual fundamental classes. This more conceptual conjecture is true in the smooth case. Using the blow-up formula, the wall crossing formula, and a case by case analysis for surfaces of Kodaira dimension less than 2, we are able to reduce our conjecture to the following assertion: $\deg \lb \Hilb^m_V  \rb = (-1)^{\chi (\ov)}$ for minimal surfaces $V$ of general type with $p_g(V)>0$ and $q(V)>0$.
\end{abstract}
\maketitle
\tableofcontents
\section*{Introduction}
This paper originated from two initially distinct projects: To study the analogue of the Poincar\'e formula for curves in higher dimension, and to develop an algebro-geometric version of Seiberg-Witten theory for projective surfaces.

The Poincar\'e formula relates the geometry of the Abel-Jacobi map to purely topological data. Given a smooth projective curve $C$ of genus $g$ over the field $\nc$, the Abel-Jacobi map
\[
\tau: C_d \longrightarrow \Pic^d_C
\]
sends an effective divisor $\mathfrak{d}$ of degree $d$ to the class of its associated line bundle $\co_C(\mathfrak{d})$. Let $\theta\in \Lambda^2 H^1(V,\nz)^\vee$ be the intersection form
\begin{eqnarray*}
\theta: \Lambda^2 H^1(V,\nz) & \longrightarrow &\nz\\
a \wedge b & \longmapsto & \langle a \cup b,[C] \rangle.
\end{eqnarray*}
The Poincar\'e formula expresses the fundamental class of the Brill-Noether locus $W_d=\tau(C_d)$ in the range $0 \leq d \leq g$ in terms of $\theta$:
\[
[W_d]= \frac{\theta^{g-d}}{(g-d)!} \cap [\Pic^d_C].
\]
Here $\theta$ is considered as a class in $H^2(\Pic^d_C,\nz)$ under the natural identification $H^2(\Pic^d_C,\nz)=\Lambda^2 H^1(V,\nz)$.

Let now $\Delta \subset C_d \times C$ be the universal divisor, choose a point $p \in C$ and set $\eta:= c_1(\co(\Delta)|_{C_d \times \{ p \}})$. The Poincar\'e formula can then be rewritten in the following form:
\begin{eqnarray*}
\tau_* \left( \sum_{i=0}^d \eta^{i} \cap [ C_d] \right) &=& \sum_{i=0}^d [W_{d-i}]\\
&=& \sum_{i=0}^d\frac{\theta^{g-d+i}}{(g-d+i)!} \cap [\Pic^d_C]
\end{eqnarray*}
Note that the expression $\tau_* \left( \sum_{i=0}^d \eta^{i} \cap [ C_d] \right)$ is the Segre class of the projective Abelian cone $\tau : C_d \to \Pic^d_C$.

When one tries to find an analogue of this formula for surfaces, the main difficulty is that the Hilbert schemes $\Hilb^m_V$ parametrizing divisors of topological type $m \in H^2(V,\nz)$ are in general not smooth of the expected dimension. They may have oversized and non-reduced components. Moreover, when $V$ varies in a smooth family, the corresponding family of Hilbert schemes is not in general flat. Hence, in order to define the Segre class of the Abelian cone
\[
\rho: \Hilb^m_V \to \Pic^m_V,
\]
one should replace the fundamental class by a {\em virtual fundamental class}
\[
\lb \Hilb^m_V \rb \in A_{\frac{m(m-k)}{2}} (\Hilb^m_V)
\]
in the Chow group of the expected dimension $\frac{m(m-k)}{2}$. Here $k=c_1(\kv)$ denotes the canonical class of $V$.

The existence of such a virtual fundamental class is a consequence of our first main result: For any flat projective morphism $v: V \to S$ there exists a relative obstruction theory in the sense of Behrend-Fantechi
\[
\varphi: L \bar{\pi}_\sharp ( \LG_{\mathbb{W}/\Hilb^P_{V/S} \times_S V} [-1] ) \to \LG_{\Hilb^P_{V/S}/S}
\]
for the relative Hilbert scheme $\Hilb^P_{V/S} \to S$. Here $\pi: \Hilb^P_{V/S} \times_S V \to \Hilb^P_{V/S}$ is the projection and $\nw \subset \Hilb^P_{V/S} \times_S V$ is the universal subscheme. If in addition $v: V \to S$ is smooth of relative dimension $d$, and $\Hilb^P_{V/S}$ parametrizes divisors, then 
\[
L \bar{\pi}_\sharp ( \LG_{\mathbb{W}/\Hilb^P_{V/S} \times_S V} [-1] ) \cong \left( R^\bul \pi_* \co_\nd (\nd) \right)^\vee,
\]
and $\left( R^\bul \pi_* \co_\nd (\nd) \right)^\vee$ is of perfect amplitude contained in $[1-d,0]$. In this formula $\nd \subset \Hilb^P_{V/S} \times_S V$ denotes the universal divisor. The theory of Behrend and Fantechi yields therefore a virtual fundamental class when $d=2$, i.e.~for curves on surfaces. For a surface $V$ and topological type $m \in H^2(V,\nz)$ we denote this cycle class by
\[
\lb \Hilb^m_V\rb \in A_{\frac{m(m-k)}{2}}(\Hilb^m_V).
\]

Using the virtual fundamental class $\lb \Hilb^m_V\rb$, we define the Poincar\'e invariant of a surface $V$ as follows: Fix $m \in H^2(V,\nz)$ and denote by $\nd^+$ and $\nd^-$ the universal divisors over $\Hilb^m_V$ and $\Hilb^{k-m}_V$ respectively. Let $p \in V$ be an arbitrary point and put $u^+:= c_1(\co(\nd^+)|_{\Hilb^m_V\times \{ p\}}$, $u^-:= c_1(\co(\nd^-)|_{\Hilb^{k-m}_V\times \{ p\}}$. Denote by $\rho^\pm$ the following morphisms:
\begin{eqnarray*}
\rho^+ : \Hilb^m_V & \longrightarrow & \Pic^m_V\\
D & \longmapsto & [\ov(D)]
\end{eqnarray*}
\begin{eqnarray*}
\rho^- : \Hilb^{k-m}_V & \longrightarrow & \Pic^m_V\\
D' & \longmapsto & [\kv(-D)]
\end{eqnarray*}
The Poincar\'e invariant of $V$ is the map
\[
(P^+,P^-): H^2(V,\nz) \longrightarrow \Lambda^* H^1(V,\nz) \times \Lambda^* H^1(V,\nz)
\]
defined by
\[
P^+_V(m)= \rho^+_* \left( \sum_{i\geq 0} (u^+)^i \cap \lb \Hilb^m_V \rb \right)
\]
and
\[
P^-_V(m)= (-1)^{\chi(\ov) + \frac{m(m-k)}{2}}\rho^-_* \left( \sum_{i\geq 0} (u^-)^i \cap \lb \Hilb^{k-m}_V \rb \right)
\]
if $m \in NS(V)$, and by $P^\pm_V(m)=0$ otherwise.

The map
\[
P^-_V: H^2(V,\nz) \longrightarrow \Lambda^* H^1(V,\nz)
\]
is determined by $P^+_V$ in the following way:
\[
[P^-_V(m)]^{2i} = (-1)^{\chi(\ov)+i} [P^+_V(k-m)]^{2i},
\]
where $[P^\pm_V(m)]^{2i}$ denotes the homogeneous component of $P^\pm_V(m)$ of degree $2i$. The reason for this redundant way of defining the Poincar\'e invariant will become clear later.

The Poincar\'e invariant possesses the following four properties:

\begin{itemize}
	\item[I)] It is invariant under smooth deformations of the surface $V$.
	\item[II)] There exists a blow-up formula relating the Poincar\'e invariant of a surface $V$ to the Poincar\'e invariant of the blow-up $\sigma: \hat{V} \to V$ of $V$ in a point.
	\item[III)] The invariant satisfies a wall-crossing formula: For surfaces $V$ with vanishing geometric genus the difference $P^+_V-P^-_V$ is a topological invariant, given by the formula
\[
P^+_V(m)-P^-_V(m)= \sum_{j=0}^{\min \{ q(V), \frac{m(m-k)}{2}\}} \frac{\theta_{2m-k}^{q(V)-j}}{(q(V)-j)!} \cap [\Pic^m_V].
\]
Here $\theta_{2m-k}$ denotes the class in $H^2(\Pic^m_V,\nz)$ corresponding to the map
\begin{eqnarray*}
\Lambda^2 H^2(V,\nz) & \longrightarrow & \nz\\
a\wedge b & \longmapsto & \frac{1}{2} \langle a \cup b \cup (2m-k),[V]\rangle.
\end{eqnarray*}
	\item[IV)] Surfaces with positive geometric genus are of simple type: For surface $V$ with $p_g(V)>0$ we have $(P^+_V(m),P^-_V(m))=(0,0)$ except for finitely many classes $m$ with $m(m-k)$.
\end{itemize}

The Poincar\'e invariant is explicitely computable for many important classes of surfaces, e.g.~for ruled surfaces, or for elliptic fibrations.

If $p: V\to C$ is a ruled surface over a curve of genus $g$, $f\in H^2(V,\nz)$ the class of a fiber, and $m \in H^2(V,\nz)$ with $m(m-k) \geq 0$, then
\[
(P^+_V(m),P^-_V(m))=  \left\{
\begin{array}{cl}
( \sum\limits_{d=0}^{\min \{ g, \frac{m(m-k)}{2}\} } (m \cdot f+1)^{g-d}[W_d],0)& \text{when $m \cdot f \geq -1$}\\
(0,-\sum\limits_{d=0}^{\min \{ g, \frac{m(m-k)}{2}\} } (m \cdot f+1)^{g-d}[W_d])& \text{when $m \cdot f \leq -1$.}
\end{array}\right.
\]
This explicit formula yields classical results of Nagata \cite{na} and Lange \cite{la} about the existence and the number of sections of $p:V \to C$ with minimal self-intersection number.

Let $\pi: V \to C$ be an elliptic fibration over a curve of genus $g$ with general fibre $F$ and multiple fibers $m_1F_1, \ldots , m_r F_r$. Fix a class $m \in H^2(V,\nz)$ with $m^2=n \cdot [F]=0$. Then
\[
P^+_V(m)=\sum_{\dfrac{d[F]+\sum a_i[F_i] =PD(m)}{0\leq a_i <m_i}} (-1)^d 
\begin{pmatrix}
2g-2 + \chi (\ov )\\
d
\end{pmatrix},
\]
\[
P^-_V(m)=\sum_{\dfrac{d[F]+\sum a_i[F_i] =PD(k-m)}{0\leq a_i <m_i}} (-1)^{\chi (\ov)+d} 
\begin{pmatrix}
2g-2 + \chi (\ov )\\
d
\end{pmatrix}.
\]
In particular, we find that
\[
P^+_V=P^-_V
\]
for elliptic surfaces with $p_g(V)>0$. 

The second project started with a question of A.~Parshin. After the talk by one of the authors on the full Seiberg-Witten invariants as defined in \cite{ot1}, he posed the question if there was a purely algebro-geometric analogue of the full Seiberg-Witten invariants for projective surfaces. In order to explain our answer to this question, let us briefly recall the structure of the full Seiberg-Witten invariant; for the construction and details we refer to \cite{ot1}. 

Let $(M,g)$ be a closed oriented Riemannian 4-manifold with first Betti number $b_1$. We denote by $b_+$ the dimension of a maximal subspace of $H^2(M,\nr)$ on which the intersection form is positive definite. The set of isomorphism classes of $\spin$-structures on $(M,g)$ has the structure of a $H^2(M,\mathbb{Z})$-torsor. This torsor does, up to a canonical isomorphism, not depend on the choice of the metric $g$ and will be denoted by $\spim$.

We have the Chern class mapping
\begin{eqnarray*}
c_1:\spim & \longrightarrow & H^2(M,\mathbb{Z}) \\
\mathfrak{c} & \longmapsto & c_1(\mathfrak{c}),
\end{eqnarray*}
whose image consists of all characteristic elements.

When  $b_+>1$, then the Seiberg-Witten invariant is a map
\[
SW_{M,\ooo} : \spim \longrightarrow \Lambda^*H^1(M,\mathbb{Z}),
\]
where $\oo$ is an orientation parameter.

When $b_+=1$, then the invariant depends on a chamber structure and is a map
\[
(SW^+_{M,\ooov},SW^-_{M,\ooov}) : \spim \longrightarrow \Lambda^*H^1(M,\mathbb{Z}) \times\Lambda^*H^1(M,\mathbb{Z}),
\]
where $\oov$ are again orientation data.

Note that the Seiberg-Witten invariant possesses four properties which are completely analogous to the properties of the Poincar\'e invariant:

\begin{itemize}
	\item[I')] It is an invariant of the oriented diffeomorphism type.
	\item[II')] There exists a formula relating the Seiberg-Witten invariant of a 4-manifold $M$ to the invariant of the connected sum $M \# \overline{\np^2}$ with $\overline{\np^2}$ \cite{os}.
	\item[III')] For 4-manifolds with $b_+=1$ the difference 
\[
\sww{+}{\mathfrak{c}}-\sww{-}{\mathfrak{c}}
\]
can be expresses in terms of purely topological data \cite{ot1}.
	\item[IV')] Taubes showed that symplectic 4-manifolds with $b_+>1$ are of simple type, i.e. ~the Seiberg-Witten invariant vanishes except for finitely many classes $\mathfrak{c}$ of virtual dimension $0$ \cite{ta2}.
\end{itemize}

We conjecture that the Seiberg-Witten- and the Poincar\'e invariants coincide for smooth projective surfaces; more precisely: Let $V$ be such a surface. Any {\em Hermitian} metric $g$ on $V$ defines a {\em canonical} $\spin$-structure on $(V,g)$. Its class $\mathfrak{c}_{can} \in \spiv$ does not depend on the choice of the metric. The Chern class of $\mathfrak{c}_{can}$ is $c_1(\mathfrak{c}_{can})=-c_1(\kv )=-k$.

Since $\spiv$ is a $H^2(V,\nz)$-torsor, the distinguished element $\mathfrak{c}_{can}$ defines a bijection:
\begin{eqnarray*}
H^2(V,\nz) & \longrightarrow & \spiv\\
m & \longmapsto & \mathfrak{c}_m
\end{eqnarray*}
The Chern class of the twisted structure $\mathfrak{c}_m$ is $2m-k$. Recall that any surface defines canonical orientation data $\oo$ and $\oov$ respectively.
\begin{conj}\label{conj:noname}
Let $V$ be a smooth projective surface, and denote by $\oo$ or $\oov$ the canonical orientation data. If $p_g(V)=0$, then
\[
P^\pm_V(m)=\swwv{\pm}{\mathfrak{c}_m}\ \ \forall m\in H^2(V,\mathbb{Z}).
\]
If $p_g(V)>0$, then
\[
P^+_V(m)=P^-_V(m)=\swv{\mathfrak{c}_m}\ \ \forall m\in H^2(V,\mathbb{Z}).
\]
\end{conj}
We consider the assertion of this conjecture as the two-dimensional analogue of the Poincar\'e  formula. In contrast to the one-dimensional case, it relates algebraic information about a smooth projective surface to {\em differential-topological} data of the underlying oriented smooth 4-manifold.

Note that if this conjecture holds, then we must have
\[
P^+_V=P^-_V
\]
for all surfaces with $V$ with $p_g(V)>0$. We have seen above that this is true for elliptic surfaces, but we have no a priori proof in the general case.

The conjecture has an important conceptual refinement:

The Kobayashi-Hitchin correspondence identifies monopoles on K\"ahler surfaces with effective divisors. To be precise: Consider a K\"ahler surface $(V,g)$, a class $m \in H^2(V,\nz)$, and a real closed $(1,1)$-form $\beta$. Let $\tau$ be a $\spin$-structure on $(V,g)$ representing the class $\mathfrak{c}_m$, and denote by $\mathcal{W^\tau_\beta}$ the moduli space of solutions to the $\beta$-twisted Seiberg-Witten equations.
\begin{itemize}
	\item[i)] If $(2m-k-[\beta]) \cdot [\omega_g] <0$, then there exists an isomorphism of real analytic spaces
\[
\kappa^+_m: \mathcal{W^\tau_\beta}\stackrel{\cong}{\longrightarrow} \Hilb^m_V.
\]
\item[ii)] If $(2m-k-[\beta]) \cdot [\omega_g] >0$, then there exists an isomorphism of real analytic spaces
\[
\kappa^-_m: \mathcal{W^\tau_\beta}\stackrel{\cong}{\longrightarrow} \Hilb^{k-m}_V.
\]
\end{itemize}

By the work of Brussee \cite{br}, the moduli space of solutions to the Seiberg-Witten equations carries a {\em virtual fundamental class} $[\mathcal{W^\tau_\beta}]_{vir}$. Moreover, the full Seiberg-Witten invariants can be computed by evaluating tautological cohomology classes on $[\mathcal{W^\tau_\beta}]_{vir}$ \cite{ot2}.
Our main conjecture is therefore essentially a consequence of the following more conceptual conjecture:
\begin{conj}
Let $(V,g)$ be a surface endowed with a K\"ahler metric $g$. Fix a class $m \in H^2(V,\nz )$ and a real closed 2-form $\beta$ of type $(1,1)$. Let $\tau$ be a $\spin$-structure on $(V,g)$ representing the class $\mathfrak{c}_m$, and denote by $\mathcal{W^\tau_\beta}$ the moduli space of solutions to the $\beta$-twisted Seiberg-Witten equations. Choose the canonical orientation data $\oo$ or $\oov$. Suppose that $(2m-k-[\beta]) \cdot [\omega_g] <0$. Then the Kobayashi-Hitchin isomorphism
\[
\kappa^+_m: \mathcal{W^\tau_\beta}\stackrel{\cong}{\longrightarrow} \Hilb^m_V
\]
identifies $[\mathcal{W^\tau_\beta}]_{vir}$ with the image of $\lb \Hilb^m_V \rb$ in $H_*(\Hilb^m_V,\mathbb{Z})$.
\end{conj}

Forthcoming work of M.~D\"urr and A.~Teleman will prove this second conjecture in the case when the moduli spaces are smooth but possibly oversized \cite{dt}.
In the present paper we use this result for two purposes: We compute the Seiberg-Witten invariants of elliptic surfaces; this fills a gap in the existing literature. Combining it with the blow-up formula and the wall-crossing formula, we reduce our first conjecture to the proof of the following assertion:

Let $V$ be a minimal surface of general type with $p_g(V)>0$ and $q(V)>0$. Then
\[
\deg\, \lb \Hilb^k_V \rb = (-1)^{ \chi (\ov )}.
\]
An alternative way of proving Conj.~\ref{conj:noname} would be to compare the Poincar\'e invariants with the Gromov invariants defined by Taubes \cite{ta}.

\noindent
{\em Acknowledgement.}
The authors like to thank H. Flenner. He helped us to understand his construction of the cotangent complex, and he suggested to use his result on the existence of a left adjoint functor (\ref{thfl}) for the construction of an obstruction theory for Hilbert schemes.
\section{Obstruction theories and virtual fundamental classes}
\subsection{Background material}
An essential ingredient in our study of Hilbert schemes and their invariants is the construction of a virtual fundamental class. Virtual fundamental classes in the context of complex geometry were first introduced by Li-Tian in \cite{lt}. In our paper we apply the formalism developped by Behrend-Fantechi in \cite{bf}. Note however that Behrend-Fantechi work with Deligne-Mumford stacks, use the \'etale topology, and obtain a virtual fundamental class in the Chow group with rational coefficients. We restrict to schemes with the Zariski topology and obtain  a virtual fundamental class in the usual Chow group, the Chow group with integer coefficients.

In the following, all schemes are separated Noetherian schemes of finite type over $\mathbb{C}$. We denote by $D^-(X)$ the category of complexes of $\ox$-modules bounded from above, and by $D^-_c(X)$ the full subcategory of complexes with coherent cohomology.

Let $X$ and $Y$ be schemes, and let $X \longrightarrow Y$ be a morphism. We denote by $\LG_{X/Y}$ the relative cotangent complex of $X$ over $Y$. It is an object in the derived category $D^-_c(\ox)$, defined up to isomorphism. When $Y=\Spec \mathbb{C}$ we denote $\LG_{X/Y}$ by $\LG_X$.

We start by recalling several facts about the cotangent complex which we will need later.
\begin{itemize}
\item $h^i(\LG_{X/Y}) = 0$ for all $i>0$. Therefore one can choose a
  complex representing $\LG_{X/Y}$ with zero terms in positive
  degrees.
\item $h^0(\LG_{X/Y}) = \Omega_{X/Y}$, the relative cotangent sheaf of
  $X$ over $Y$.
\item Let $\xymatrix@1{ X \ar[r]^f & Y \ar[r]^g & Z}$ be two morphisms of
schemes. They induce a distinguished triangle in $D^-_c(\ox)$:
\[
\xymatrix{
f^* \LG_{Y/Z} \ar[r] & \LG_{X/Z} \ar[d] \\
& \LG_{X/Y} \ar[ul]^{[1]}}
\]
\item Let
\[
\xymatrix{
X' \ar[d]^{f'} \ar[r]^{j'} & X \ar[d]^f \\
Y' \ar[r]^j & Y}
\]
be a commutative square. Then there is a natural morphism
\begin{equation}\label{eq:pb}
{j'}^* \LG_{X/Y} \longrightarrow \LG_{X'/Y'}
\end{equation}
obtained by composing the morphisms
\[ {j'}^* \LG_{X/Y} \longrightarrow \LG_{X'/Y} \ \ \mbox{and} \ \
\LG_{X'/Y} \longrightarrow \LG_{X'/Y'}. \]
If the commutative square is Cartesian, then
\[ h^{i} ({j'}^* \LG_{X/Y}) \longrightarrow h^{i} (\LG_{X'/Y'}) \]
is an isomorphism for $i=0$ and surjective for
$i=-1$ \cite[II.1.1.2.9]{ill}.

If in addition $\mathcal{T}or_i^{\mathcal{O}_Y} (\mathcal{O}_X, \mathcal{O}_{Y'})=0$ for
all $i>0$, then the morphism \eqref{eq:pb} is an isomorphism \cite[Cor.II.2.3.10]{ill}. Note that this condition is satisfied if one of the morphisms $f$
or $j$ is flat.
	\item If $g: X \to Y$ is a regular embedding, then $\LG_{X/Y}  \cong \mathcal{N}^\vee_{X/Y}[1]$ \cite[Cor.III.3.2.7]{ill}.
\end{itemize}
Let $Y$ be a fixed base scheme, $T$ a scheme over $Y$, and $\JJ$ a coherent $\co_T$-module. A closed immersion of $T \to \tbar$ of schemes over $Y$ is a {\em square-zero extension} with ideal $\JJ$ if $\JJ^2_{T/\tbar}=0$ and  $\JJ_{T/\tbar}$ considered as $\co_T$-module is $\JJ$. 
Let now $T \to \tbar$ be a square-zero extension with ideal $\JJ$. Combining the morphism
\[
\LG_{T/Y} \to \LG_{T/\tbar}
\]
with the cut-off morphism
\[
\LG_{T/\tbar} \to \JJ[1]
\]
yields an element $[T \to \tbar] \in \ext^1(\LG_{T/Y},\JJ)$. For a morphism $f: T \to X$ of schemes over $Y$ the associated morphism
\[
f^*\LG_{X/Y} \to \LG_{T/Y}
\]
induces a map $\ext^1(\LG_{T/Y},\JJ) \to \ext^1(f^*\LG_{X/Y},\JJ)$. Let $\omicron [T \to \tbar] \in \ext^1(f^*\LG_{X/Y},\JJ)$ be the image of the class $[T \to \tbar]$. Recall the following facts from deformation theory:
\begin{itemize}
	\item the morphism $f: T \to X$ extends to a morphism $\fbar: \tbar \to  X$ if and only if $\omicron[T \to \tbar]=0$
	\item if $\omicron [T \to \tbar]=0$ then the set of extensions is a torsor under $\ext^0(f^*\LG_{X/Y},\JJ)$.
\end{itemize}
Recall that a complex of sheaves is \emph{of perfect amplitude contained in
$[a,b]$}, where $a,b \in \mathbb{Z}$, if, locally, it is isomorphic in the
derived category to a complex $\FF^a \rightarrow \ldots \rightarrow \FF^b$ of locally free sheaves of finite rank.

If the complex $\ee$ is of perfect amplitude contained in $[a,0]$ for some
$a$, then the assignment, which assigns to every geometric point $j: \{x\}
\hookrightarrow X$ the alternating sum
$\sum_i (-1)^{i} dim\ h^{i}(j_x^* \ee)$, is locally constant. If this number is
globally constant, we will speak of the rank of $\ee$ and denote it by $\rk \ee$.
\begin{dfn} Let $\ee$ be an object in the derived category $D^-_c(\ox)$. A
global resolution of $\ee$ of perfect amplitude contained in $[a,b]$ is an
isomorphism $\ff \stackrel{\cong}{\longrightarrow} \ee$, where
\[
\ff = [ \FF^{a} \rightarrow \FF^{a+1} \rightarrow \cdots \rightarrow
\FF^{b-1} \rightarrow \FF^b ]
\]
is a complex of locally free sheaves. If $[a,b]=[-1,0]$, we say that $\ff$ is a
perfect global resolution of $\ee$.
\end{dfn}

\begin{dfn}
Let $X$ be a scheme. A {\em vector bundle stack} over $X$ is an Artin stack over $X$, which is locally isomorphic to the stack $F_1/F_0$ defined by a morphism $\alpha: F_0 \to F_1$ of vector bundles on $X$.
\end{dfn}
For later use, we note that any complex $\ee$ of perfect amplitude contained in  $[-1,0]$ defines a vector bundle stack, which we will denote by $E$. If $\FF^{-1} \to \FF^0$ is a perfect resolution of $\ee$, then $E$ is isomorphic to the quotient $F_1/F_0$, where $F_i: = \Spec Sym \FF^{-i}$.

By the work of Kresch we have:
\begin{thm}
Let $X$ a scheme, and let $pr: E \to X$ be a vector bundle stack of rank $r$ on $X$. The pull-back morphism
\[
pr^*: A_*(X) \longrightarrow A_{*+r}(E)
\]
is an isomorphism of Abelian groups.
\end{thm}
\begin{proof}
\cite[Thm.2.1.12]{kr}.
\end{proof}
\begin{nota}
Given a scheme $X$ and a vector bundle stack $pr: E \to X$ of rank $r$ on $X$, we denote by
\[
0_E^!: A_{*+r}(E) \longrightarrow A_*(X)
\]
the induced morphism.
\end{nota}

Let again $X$ be a scheme over $Y$. These data define the {\em relative intrinsic normal cone} $\INC{X}{Y}$ (\cite{bf}); it is an Artin stack over $X$ of relative dimension $0$ over $Y$.

\begin{dfn}
Let $X \longrightarrow Y$ be a morphism of schemes, and let $\ee$ be an object in the derived category $D^-(\ox)$.
Suppose, that $h^{i}(\ee) =0$ for $i>0$ and that $h^{i}(\ee)$ is coherent for $i=-1,0$. A morphism $\varphi: \ee \longrightarrow \LG_{X/Y}$ is called a
\emph{relative obstruction theory} for $X$ over $Y$, if $h^0(\varphi) : h^0(\ee)
\longrightarrow h^0(\LG_{X/Y})$ is an isomorphism and $h^{-1}(\varphi) :
h^{-1}(\ee)
\longrightarrow h^{-1}(\LG_{X/Y})$ is an epimorphism. An obstruction theory for $X$ is a relative obstruction theory for $X$ over $Y=\Spec \mathbb{C}$.
\end{dfn}
Let $\varphi: \EE^\bul \to \LG_{X/Y}$ be a relative obstruction theory. Then the induced map
\[
f^*(\varphi) : \ext^1(f^*\LG_{X/Y},\JJ) \to \ext^1(f^*\EE^\bul,\JJ)
\]
is injective and
\[
f^*(\varphi) : \ext^0(f^*\LG_{X/Y},\JJ) \to \ext^0(f^*\EE^\bul,\JJ)
\]
is a bijection. This implies that a morphism $f: T \to X$ extends to a morphism$\fbar: \tbar \to  X$ if and only if $f^*(\varphi)(\omicron[T \to \tbar]) \in \ext^1(f^*\EE^\bul,\JJ)$ vanishes. If this is the case, then the set of extensions is a torsor under $\ext^0(f^*\EE^\bul,\JJ)$. Conversely, we have the following criterion:
\begin{thm}\label{thm:bfcrit}
Let $X$ be a scheme over $Y$. Suppose $\EE^\bul$ is an object in the derived category $D^-(\ox)$ with vanishing cohomology in positive degrees and coherent côhomology $h^{i}(\EE^\bul)$ for $i=-1,0$.

A morphism $\varphi: \EE^\bul \to \LG_{X/Y}$ is a relative obstruction theory if and only if for all morphisms $f: T \to X$, for all coherent $\co_T$-modules $\JJ$, and for all square-zero extensions $T \to \tbar$ over $Y$ with ideal $\JJ$ the following conditions are satisfied:
\begin{itemize}
	\item[i)] The morphism $f: T \to X$ extends to a morphism $\fbar: \tbar \to X$  over $S$ if and only if $f^*(\varphi)(\omicron [T \to \tbar])\in \ext^1(f^* \EE^\bul, \JJ)$ vanishes;
	\item[ii)] The map $f^*(\varphi): \ext^0(f^*\LG_{X/Y},\JJ) \to \ext^0(f^*\EE^\bul, \JJ)$ is a bijection.
\end{itemize}
\end{thm}
\begin{proof}
\cite[Thm.4.5]{bf}.
\end{proof}
\begin{rem} Let $\varphi: \ee \to \LG_{X/Y}$ be a relative obstruction theory. Then for any morphism $f: T \to X$, the functor $\operatorname{Ext}^1(f^*\ee, -)$ is an obstruction theory in the sense of Buchweitz-Flenner \cite[Def. 6.10, 6.14]{buf}. If in addition $\ee$ is of perfect amplitude contained in $[-1,0]$ and $Y=\Spec \nc$, then $h^1({\ee}^\vee)$ is an obstruction theory for $X$ in the sense of Li-Tian \cite{lt}.
\end{rem}
Let $\varphi: \ee \to \LG_{X/Y}$ be a relative obstruction theory, and suppose that $\ee$ is of perfect amplitude contained $[-1,0]$. Then the obstruction theory defines a closed embedding
\[
\INC{X}{Y} \hookrightarrow E.
\]
\begin{dfn}
Let $X \to Y$ be a morphism of schemes, and fix a relative obstruction theory $\varphi : \ee \longrightarrow \LG_{X/Y}$. Suppose that $\ee$ is of perfect amplitude contained in [-1,0] and that $Y$ is of pure dimension $l$. The {\em virtual fundamental class} of $X$ with respect to the obstruction theory $\varphi$ is
\[
[X,\varphi]:= 0_E^! [\INC{X}{Y}] \in A_{l+\rk \ee}(X).
\]
\end{dfn}

Let now $X$ be a scheme over a base scheme $Y$ of pure dimension $l$, and let
$\varphi : \ee \longrightarrow \LG_{X/Y}$ be a relative obstruction theory for $X$ over $Y$. Suppose that $\ee$ admits a perfect global resolution
\[
\ff = [ \FF^{-1} \rightarrow \FF^0 ] \stackrel{\cong}{\longrightarrow}
\ee.
\]
Set $F_i : = \Spec Sym \FF^{-i}$, and denote by
$\varphi_{\ff}$ the induced morphism $\ff \longrightarrow \LG_{X/Y}$.

\[
\xymatrix{
C(\varphi_{\ff}) \ar[d] \ar[r] & F_1 \ar[d] \\
\INC{X}{Y} \ar[r] & F_1/F_0 }
\]
Then $C(\varphi_{\ff})$ is a closed subcone of $F_1$ of pure dimension $l + \rk
F_0$.
The virtual fundamental class $[X,\varphi]$ is the intersection of $C(\varphi_{\ff})$ with the zero section of
$F_1$:
\[
\xymatrix{
X \ar@{=}[d] \ar[r] & C(\varphi_{\ff}) \ar[d]\\
X \ar[r]_{0_{F_1}} & F_1}
\]
\[
[X,\varphi] = 0_{F_1}^! [C(\varphi_{\ff})] \in A_{l+\rk \ee}(X).
\]
\subsection{First properties of virtual fundamental classes}
\begin{prop}[Locally free obstructions]
Let $\varphi : \ee \longrightarrow \LG_{X/Y}$ be a relative obstruction
theory for $X$ over an equidimensional scheme $Y$, and suppose that $\ee$
is of perfect amplitude contained $[-1,0]$.

i) If $h^1 ({\ee}^{\vee})=0$, then $X$ is smooth over
$Y$ and $[X,\varphi]=[X]$, the usual fundamental class of $X$.

ii) If $X$ is smooth over $Y$, then $h^1 ({\ee}^{\vee})$ is locally free
and $[X,\varphi]= c_r(h^1 ({\ee}^{\vee})) \cap [X]$, where $r=\mbox{rk}\
h^1({\ee}^{\vee})$.
\end{prop}
\begin{proof}
\cite[Prop.~7.3]{bf}
\end{proof}
\begin{prop}[Base change]\label{prop:bc}
Let
\[
\xymatrix{X' \ar[r]^{j'} \ar[d]^{f'} & X \ar[d]^f \\
              Y' \ar[r]^j & Y}
\]
be a Cartesian square with equidimensional base schemes $Y$ and $Y'$. If
$\varphi: \ee \longrightarrow \LG_{X/Y}$ is a relative obstruction theory
for $X$ over $Y$, then the induced morphism $\varphi' : {j'}^* \ee
\longrightarrow \LG_{X'/Y'}$ is a relative obstruction theory for $X'$
over $Y'$. If $\varphi$ admits a perfect global resolution, then so does
$\varphi'$.

\noindent
If in addition $j$ is flat, or $j$ is a regular local immersion, then there is an equality of the corresponding virtual fundamental classes
\[
j^! [ X , \varphi] = [ X' , \varphi'],
\]
where $j^!$ denotes the refined Gysin map $A_*(X) \longrightarrow A_*(X')$.
\end{prop}
\begin{proof}
\cite[Prop.~7.2]{bf}
\end{proof}
\begin{prop}\label{prop:besie}
Let $\varphi: \EE^\bullet \to \LG_X$ be a perfect obstruction theory for a scheme $X$, and suppose that $X$ can be embedded into a smooth variety. Then
\[
[X,\varphi] = \left( c({\EE^\bullet}^\vee)^{-1} \cap c_*(X) \right)_{\rk \EE^\bullet},
\]
where $c_*(X)$ is Fulton's canonical class.
\end{prop}
\begin{proof}
\cite[Thm.4.6]{si}
\end{proof}
\begin{rem}
It follows that the virtual fundamental class depends only on the complex $\ee$ and not on the morphism $\varphi: \ee \to \LG_X$ when $X$ can be embedded into a smooth variety. If $X$ is proper, then there is also a direct argument for this observation: Let $\varphi_0: \ee \to \LG_X$ and $\varphi_1: \ee \to \LG_X$ be two obstruction theories, and set $\varphi_t: (1-t) \varphi_0+t\varphi_1$. Then for almost all $t \in \nc$, $\varphi_t$ is an obstruction theory. Hence, if $\FF^{-1} \to \FF^0$ is a global resolution of $\ee$, we obtain a family of cones $C_t$  in the vector bundle $F_1$ dual to $\FF^{-1}$. This family is defined for all $t$ in a Zariski open subset of $\np^1$, which contains $0$ and $1$. By taking the closure in $\np^1 \times F_1$, we obtain a family of cones over $\np^1$. Hence the classes defined by $C_0$ and $C_1$ in the Chow group of $F_1$ agree.
\end{rem}
\subsubsection{The basic example}
Let
\[
\xymatrix{
X \ar[r]^{g'} \ar[d]^{f'} & V \ar[d]^f \\
Y \ar[r]^g & W}
\]
be a Cartesian diagram of schemes. Compose the morphism
\[
{f'}^* \LG_{Y/W} \longrightarrow \LG_{X/V}
\]
with 
\[
\LG_{X/V} \longrightarrow {g'}^* \LG_{V} [1]
\]
and let $\mathcal{E}^{\bullet}$ denote the mapping cone of
\[
{f'}^* \LG_{Y/W}[-1] \longrightarrow {g'}^* \LG{V}.
\]
Let $\mathcal{A}^{\bullet}$ be the mapping cone of
\[
{f'}^* \LG_{Y/W} \longrightarrow \LG_{X/V}.
\]
Then we have the following diagram
\[
\xymatrix{
{f'}^* \LG_{Y/W}  \ar[r] \ar@{=}[d] &
\LG_{X/V} \ar[d] \ar[r] & \mathcal{A}^{\bullet} \ar[r] \ar@{.>}[d] &
{f'}^* \LG_{Y/W}[1] \ar@{=}[d] \\
{f'}^* \LG_{Y/W} \ar[r] & {g'}^*
\LG_{V} [1] \ar[r] \ar[d] & \mathcal{E}^{\bullet}[1]
\ar@{.>}[d]^{\varphi [1]} \ar[r] & {f'}^* \LG_{Y/W}[1] \\
& \LG_{X}[1] \ar[d] \ar@{=}[r] &
\LG_{X}[1] \ar[d] & \\
& \LG_{X/V} [1] \ar[r] & \mathcal{A}^{\bullet}[1] & }
\]
where the dotted arrows exist according to the octahedral axiom \cite[p. 21]{Ha2}.
\begin{prop}
Let
\[
\xymatrix{
X \ar[r]^{g'} \ar[d]^{f'} & V \ar[d]^f \\
Y \ar[r]^g & W}
\]
be a Cartesian diagram of schemes. The induced morphism $\varphi: \ee \longrightarrow \LG_{X}$ is an
obstruction theory for $X$. If $V$ is smooth and $Y
\stackrel{g}{\longrightarrow} W$ is a regular embedding with ideal sheaf
$\mathcal{J}$, then $\mathcal{J}/\mathcal{J}^2$ is a locally free sheaf
on $Y$ and ${f'}^* \mathcal{J}/\mathcal{J}^2 \longrightarrow {g'}^*
\Omega_V$ is a perfect global resolution of $\mathcal{E}^{\bullet}$. If in
addition $W$ is smooth, then we have
\[ [X,\varphi] = g^! [V], \]
where $g^!$ denotes the refined Gysin map $A_*(V) \longrightarrow A_*(X)$.
\end{prop}
\begin{proof}
\cite[p.81]{bf}.
\end{proof}
\begin{cor}
Let $E \to V$ be a vector bundle on a smooth variety, and let $\xi$ be a section. Then the zero locus $Z(\xi)$ comes with a preferred obstruction theory, and the associated virtual fundamental class is the localized Euler class of $E$.
\end{cor}
\begin{proof}
The following diagram is Cartesian:
\[
\xymatrix{
Z(\xi) \ar[r]^j \ar[d]_j & V\ar[d]^\xi\\
V \ar[r]_{0_E} & E}
\]
Therefore our claim is a direct consequence of the above proposition.
\end{proof}
\subsubsection{Associativity}
\begin{lem}\label{lem:ass}
Let
\[
0 \to E' \to E \to E/E' \to 0
\]
be a short exact sequence of vector bundles on a scheme $X$, and let $\xi$ be a
section of $E$. Denote by $\xibar$ the induced section of $E/E'$, and by $\xi'$ the induced section of $E'|_{Z(\xibar)}$. Then the diagram
\[
\xymatrix{A_* (X) \ar[r]^{0^!_{E/E'}} \ar@{=}[d] & 
A_*(Z(\xibar)) \ar[d]^{0^!_{E'}} \\
 A_* (X) \ar[r]^{0^!_E} & A_* (Z(\xi))}
\]
commutes. 
\end{lem}
\begin{proof}
Consider the diagram
\[
\xymatrix{Z(\xi) \ar[r] \ar[d] & Z(\xibar) \ar[r] \ar[d]^{\xi'} &
X \ar[dd]_\xi \ar@/^1pc/[ddd]^{\xibar} \\
Z(\xibar) \ar[r]^{0_{E'}} \ar[d] & E'|_{Z(\xibar)} \ar[d] & \\
X \ar[r]^{0_{X,E'}} & E' \ar[r] \ar[d]^\alpha & E \ar[d] \\
& X \ar[r]^{0_{E/E'}} & E/E',}
\]
where all squares are cartesian. By \cite[Thm.~6.2(c)]{Fu}, we have $0_{E'}^!=0_{X,E'}^! $ and $0_{E/E'}^!=\alpha^!$. Therefore the functoriality of refined Gysin maps \cite[Thm.~6.5]{Fu} implies
\begin{eqnarray*}
0_{E'}^! \circ 0_{E/E'}^! &=&0_{X,E'}^! \circ \alpha^!\\
&=&0_E^!.
\end{eqnarray*}
\end{proof}
\begin{cor}
Let
\[
0 \to E' \to E \to E/E' \to 0
\]
be a short exact sequence of vector bundles on a scheme $X$, and let $\xi$ be a
section of $E$. Denote by $\xibar$ the induced section of $E/E'$, by $\xi'$ the induced section of $E'|_{Z(\xibar)}$, and by $\iota$ the inclusion $Z(\xi) \hookrightarrow Z(\xibar)$. Let $[[Z(\xi)]]$ and $[[Z(\xibar]]$ be the localized Euler classes of the zero loci. Then
\[
[[Z(\xi)]] = 0_{E'}^![[Z(\xibar)]],
\]
and
\[
\iota_* [[Z(\xi)]]= c_{top}(E') \cap [[Z(\xibar)]].
\]
\end{cor}
\begin{proof}
This is an immediate consequence of the above lemma.
\end{proof}
For a more general statement concerning associativity (or functoriality) of virtual fundamental classes, see \cite[Thm.1]{kkp}.
\subsubsection{Excess intersection}
\begin{prop}\label{prop:ex}
Let $\varphi : \ee \to \LG_{X/Y}$ be an obstruction theory for a scheme $X$ over a scheme $Y$ of pure dimension $d$, let $\ff = [\mathcal{F}^{-1} \longrightarrow \mathcal{F}^0]$ be a perfect global resolution of $\ee$, and let $\varphi_{\ff}: \ff \to \LG_{X/Y}$ be the corresponding morphism. Suppose that there exists a subvectorbundle $G_1 \subset F_1:= \Spec Sym \mathcal{F}^{-1}$ such that $C(\varphi_{\ff})$ is contained in $G_1$, and denote by $\{X\} \in A_{d+\rk F_0-\rk G_1} (X)$ the cycle class obtained by intersecting $C(\varphi_{\ff})$ with the zero section of $G_1$ in $G_1$. Then the virtual fundamental class $[X,\varphi]$ of $X$ with respect to the obstruction theory $\varphi$ is given by
\[
[X,\varphi] = c_{top}(F_1/G_1) \cap \{X\}.
\]
\end{prop}
\begin{proof}
Consider the following Cartesian diagram
\[
\xymatrix{
X \ar@{=}[d] \ar[r] & C(\varphi_{\FF^\bullet})\ar[d]\\
X \ar@{=}[d] \ar[r]^{0_{G_1}} & G_1 \ar[d]\\
X \ar[r]^{0_{F_1}} & F_1,}
\]
where $0_{G_1}$ and $0_{F_1}$ are the zero sections of the corresponding vector bundles. By definition, we have
\[
\{ X \} = 0_{G_1}^! [C(\varphi_{\FF^\bullet})]
\]
and
\[
[X,\varphi] = 0_{F_1}^! [C(\phi_{\FF^\bullet})].
\]
Hence our claim follows from \cite[Thm.6.3]{Fu}.
\end{proof}
\section{Hilbert schemes}
\subsection{An obstruction theory for Hilbert schemes}
In this section we want to construct an obstruction theory for Hilbert schemes in the sense of Behrend and Fantechi. An essential ingredient in our construction is the following result of Flenner.
\begin{thm}[Flenner]\label{thfl}
Let $h: M \longrightarrow
N$ be a flat proper morphism of schemes. If $N$ has a dualizing complex,
then there exists a functor
\[
Lh_\sharp : D^-_c (M) \longrightarrow D^-_c (N)
\]
satisfying the following properties:\\
(i) for $\mathcal{F}^{\bullet} \in D^-_c(M)$ and $\mathcal{G}^{\bullet}
\in D^-_c(N)$ there exists a natural isomorphism
\[
Rh_* R\mathcal{H}om_M(\mathcal{F}^{\bullet}, h^*\mathcal{G}^{\bullet} )
\cong R\mathcal{H}om_N(Lh_\sharp \mathcal{F}^{\bullet} ,
\mathcal{G}^{\bullet} );
\]
(ii) if
\[
\xymatrix{M' \ar[r]^{k'} \ar[d]^{h'} & M \ar[d]^h \\
          N' \ar[r]^k & N}
\]
is a Cartesian square, then there exists a natural isomorphism $Lk^*
Lh_\sharp
\cong Lh'_\sharp L{k'}^*$.
\end{thm}
\begin{proof}
\cite[Satz 2.1]{flenner}.
\end{proof}
For the remainder of this section we assume that $v: V \to S$ is a flat projective morphism. We fix a relatively very ample sheaf $\ov(1)$, a polynomial $P$,  and denote by $\Hilb^P_{V/S}$ the corresponding relative Hilbert scheme. Let$\mathbb{W}$ be the universal subscheme of $\Hilb^P_{V/S} \times_S V$:
\[
\xymatrix{ \mathbb{W} \ar[r]^-{i} \ar[dr]^{\bar{\pi}} & 
  \Hilb^P_{V/S} \times_S V \ar[d]^\pi \ar[r]^-{pr}  & V \ar[d]^v \\
           & \Hilb^P_{V/S} \ar[r] & S.} 
\]
We get a natural morphism
\[
\LG_{\mathbb{W}/\Hilb^P_{V/S} \times_S V} [-1]  \longrightarrow i^* \LG_{\Hilb^P_{V/S}\times_S
V/V}.
\]
On the other hand we also have a canonical morphism
\[
\pi^* \LG_{\Hilb^P_{V/S}/S} \longrightarrow \LG_{\Hilb^P_{V/S} \times_S V/V},
\]
which is an isomorphism, since $v$ is flat. So we obtain a morphism
\[
\LG_{\mathbb{W}/\Hilb^P_{V/S} \times_S V} [-1]  \longrightarrow \bar{\pi}^* \LG_{\Hilb^P_{V/S}/S}.
\]
By assumption, the scheme $S$ is of finite type over $\mathbb{C}$, hence admits a dualizing complex \cite[p. 299]{Ha2}. Since the
relative Hilbert scheme $\Hilb^P_{V/S}$ is of finite type over $S$, and $S$ is Noetherian, also $\Hilb^P_{V/S}$ admits a dualizing complex \cite[p. 299]{Ha2}. So we are in the situation of Thm. \ref{thfl} and
may apply the functor $L\bar{\pi}_\sharp$ to the morphism above. We obtain a morphism
\begin{equation}\label{eq:defob}
\varphi : L \bar{\pi}_\sharp ( \LG_{\mathbb{W}/\Hilb^P_{V/S} \times_S V} [-1] ) \longrightarrow L \bar{\pi}_\sharp ( \bar{\pi}^*
\LG_{\Hilb^P_{V/S}/S}) \longrightarrow \LG_{\Hilb^P_{V/S}/S},
\end{equation}
where the second morphism is the canonical morphism associated to a pair
of adjoint functors. Put $\EE^\bullet:=L \bar{\pi}_\sharp ( \LG_{\mathbb{W}/\Hilb^P_{V/S} \times_S V} [-1] )$.
We will prove that the morphism $\varphi: \EE^\bullet \to \LG_{\Hilb^P_{V/S}/S}$ is a relative obstruction theory. To this end we need a few preparations.
First, a result of Illusie:
Consider morphisms
\[
\xymatrix{
X \ar[d]^{f}\\
Y \ar[d]^{q}\\
S,}
\]
coherent sheaves $\mathcal{I}$ and $\mathcal{J}$ on $X$ and $Y$, and a morphism $\nu: \mathcal{J} \to f_* \mathcal{I}$. Define the following map:
\[
\mu : \ext^1(\LG_{Y/S},\mathcal{J}) \to \ext^1(f^*\LG_{Y_/S},f^*\mathcal{J}) \to \ext^1(f^*\LG_{Y/S},\mathcal{I}) \to \ext^2( \LG_{X/Y},\mathcal{I}).
\]
\begin{thm}[Illusie]\label{thm:ill}
Let $Y \to \ybar$ be a square-zero extension over $S$ with ideal $\mathcal{J}$.
There exists a  square zero extension $X \to \xbar$ with ideal $\mathcal{I}$ and a morphism $\fbar: \xbar \to \ybar$ such that
\begin{itemize}
	\item the following diagram commutes
\[
\xymatrix{
X \ar[d]^{f} \ar[r] & \xbar \ar[d]^{\fbar}\\
Y \ar[d]^{q} \ar[r] & \ybar \ar[dl]\\
S}
\]
	\item and gives rise to a morphism of extensions
\[
\xymatrix{
0  \ar[r] & \mathcal{J} \ar[r] \ar[d]_{\nu}& \co_\ybar \ar[r] \ar[d] & \co_{Y} \ar[r] \ar[d] & 0\\
0  \ar[r] & f_* \mathcal{I} \ar[r] & \fbar_* \co_\xbar \ar[r] & \co_{X} &,}
\]
\end{itemize}
if and only if $\mu [Y \to \ybar] \in \ext^2(\LG_{X/Y},\mathcal{I})$ vanishes. If $\mu [Y \to \ybar]=0$, then the set of isomorphism classes of pairs $(\xbar,\fbar)$ is an $\ext^1(\LG_{X/Y},\mathcal{I})$-torsor.
\end{thm}
\begin{proof}
\cite[Thm.III.2.1.7]{ill}.
\end{proof}
Now fix a scheme $T$ over $S$ and a morphism $f: T \to \Hilb^P_{V/S}$, and let $W_T \subset T \times_S V$ be the subscheme corresponding to $f$. We obtain the following commutative diagram:
\begin{equation*}
  \xymatrix@=20pt{ 
& T \times_S V  \ar '[d] [dd]^{\pi_T} \ar[rr]
      & & \Hilb^P_{V/S}\times_S V \ar[dd]^{\pi} \\
W_T \ar[ur]^{i_T} \ar[dr]^{\pibar_T} \ar[rr]^>>>>>>>>{F} 
      & &\nw \ar[ur]^i \ar[dr]^\pibar & \\
& T \ar[rr]^f && \Hilb^P_{V/S} \\
}
\end{equation*}
Consider the composition
\[
\lambda: \ext^1(\LG_{T \times_S V/V}, \pi_T^*\JJ) \to \ext^1( \i_T^*\LG_{T \times_S V/V}, \pibar_T^*\JJ) \to \ext^2( \LG_{W_T/T \times_S V}, \pibar_T^*\JJ).
\]
\begin{lem}\label{lem:obstack}
Let $T \to \tbar$ be a square-zero extension over $S$ with ideal $\JJ$. The morphism $f: T \to \Hilb^P_{V/S}$ extends to $\tbar$ if and only if
\[
\lambda [T \times V \to \tbar \times_S V] \in \ext^2( \LG_{W_T/T \times_S V}, \pibar_T^*\JJ)
\]
vanishes. If $\lambda [T \times V \to \tbar \times_S V]=0$, then the set of extensions  $\fbar: \tbar \to \Hilb^P_{V/S}$ is an $\ext^1( \LG_{W_T/T \times_S V}, \pibar_T^*\JJ)$-torsor.
\end{lem}
\begin{proof}
We apply Illusie's theorem to the following situation:
\[
\xymatrix{
W_T \ar[d]^{i_T}&\\
T \times_S V \ar[r] \ar[d]& \tbar \times_S V \ar[ld]\\
V &}
\]
where
\[
\nu : \pi_T^* \JJ \to (i_T)_*\pibar_T^* \JJ
\]
is the canonical adjoint morphism. If $\wbar_T \to \tbar \times_S V$ defines a morphism $\fbar : \tbar \to \Hilb^P_{V/S}$ which extends $f: T \to \Hilb^P_{V/S}$, then the inlusion $W_T \to \wbar_T$ is a square-zero extension with ideal sheaf $\pibar_T^*\JJ$, since $\wbar_T \to \tbar$ is flat. Conversely, if we have a square-zero extension $W_T \to \wbar_T$ with ideal sheaf $\pibar_T^*\JJ$ and a morphism $\wbar_T \to \tbar \times_S V$ satisfying the two conditions of Illusie's theorem, then $\wbar_T \to \tbar$ is flat \cite[Lemme III.2.1.1.1]{ill}. Therefore,
\[
\xymatrix{
W_T \ar[d] \ar[r]& \wbar_T \ar[d]\\
T \times_S V \ar[r]& \tbar \times_S V }
\]
defines a morphism $\fbar : \tbar \to \Hilb^P_{V/S}$ which extends $f: T \to \Hilb^P_{V/S}$. 
\end{proof}
\begin{thm}\label{thm:ob}
Let $v: V \to S$ be a flat projective morphism, and put $\EE^\bul:=L \bar{\pi}_\sharp ( \LG_{\mathbb{W}/\Hilb^P_{V/S} \times_S V} [-1] )$. The morphism $\varphi: \EE^\bullet \to \LG_{\Hilb^P_{V/S}/S}$ defined in \eqref{eq:defob} is a relative obstruction theory for the relative Hilbert scheme $\Hilb^P_{V/S}$ over $S$.
\end{thm}
\begin{proof}
We will check the conditions $i)$ and $ii)$ of Behrend-Fantechi's criterion (Thm.~\ref{thm:bfcrit}).

$i)$ Note that the canonical morphism
\[
\zeta : F^* \LG_{\nw/\Hilb^P_{V/S}\times_S V} \to \LG_{W_T/T \times_S V}
\]
is an isomorphism since both $V \to S$ and $\nw \to \Hilb^P_{V/S}$ are flat.
Analogously to $\varphi: \EE^\bul \to \LG_{\Hilb^P_{V/S}/S}$ we define a morphism
\[
\varphi_T: L{\pibar_T}_\sharp (\LG_{W_T/T \times_S V}[-1]) \to \LG_{T/S}.
\]
Since the diagram
\[
\xymatrix{
W_T \ar[r]\ar[d]& \nw \ar[d]\\
T \times_S V \ar[r]& \Hilb^P_{V/S} \times_S V
}
\]
is Cartesian, naturality of Flenner's functor yields canonical isomorphisms:
\begin{eqnarray*}
f^*\EE^\bul &=& f^* (L\pibar_\sharp (\LG_{\nw/\Hilb^P_{V/S}\times_S V}[-1]))\\
& \cong & F^*( L {\pibar_T}_\sharp (\LG_{\nw/\Hilb^P_{V/S}\times_S V}[-1]))\\
& \cong & L{\pibar_T}_\sharp (\LG_{W_T/T \times_S V}[-1])
\end{eqnarray*}
Using again the naturality of the functor $L\pibar_\sharp$ we see that the following diagram commutes:
\begin{equation}\label{eq:comob}
\xymatrix{
f^*\EE^\bul \ar[r]^{f^*(\varphi)} \ar[d]_{\cong} & f^* \LG_{\Hilb^P_{V/S}/S} \ar[d] \\
L{\pibar_T}_\sharp (\LG_{W_T/T \times_S V}[-1]) \ar[r]_{\varphi_T} & \LG_{T/S}}
\end{equation}
 Let $T \to \tbar$ be a square-zero extension over $S$ with ideal $\JJ$. We want to show that the element
\[
f^*(\varphi)( \omicron[T \to \tbar ]) \in \ext^1 (f^* \EE^\bul,\JJ)
\]
vanishes, if and only if the morphism $f: T \to \Hilb^P_{V/S}$ extends to a morphism $\fbar : \tbar \to \Hilb^P_{V/S}$. 
The canonical isomorphism
\[
f^*\EE^\bul \stackrel{\cong}{\longrightarrow} L{\pibar_T}_\sharp (\LG_{W_T/T \times_S V}[-1])
\]
induces an isomorphism
\[
\ext^2(\LG_{W_T/T \times_S V}, \pibar^*_T\JJ) \to \ext^1(L{\pibar_T}_\sharp (\LG_{W_T/T \times_S V}[-1]) ,\JJ) \to \ext^1(f^*\EE^\bul,\JJ).
\]
Using Lemma \ref{lem:obstack} it suffices to show that the element $\lambda[T \times_S V \to \tbar \times_S V] \in \ext^2(\LG_{W_T/T \times_S V}, \pibar^*_T\JJ)$ is mapped to $f^*(\varphi)( \omicron[T \to \tbar ]) \in \ext^1 (f^* \EE^\bul,\JJ)$ under this isomorphism. Commutativity of the diagram \eqref{eq:comob} yields that the following diagram commutes:
\[
\xymatrix{
& \ext^1(\LG_{T/S},\JJ) \ar[d]  \ar[dl]\ar[r]& \ext^1(f^*\LG_{\Hilb^P_{V/S}/S},\JJ) \ar[ddd]\\
\ext^1(\pibar_T^*\LG_{T/S},\pibar_T^*\JJ)  & \ext^1(\pi_T^*\LG_{T/S},\pi_T^*\JJ) \ar[l]&\\
\ext^1( i_T^* \LG_{T \times_S V/V},\pibar_T^*\JJ)\ar[u]^{\cong} \ar[d] & \ext^1( \LG_{T \times_S V/V},\pi_T^*\JJ) \ar[l] \ar[u]^{\cong} \ar[dl]^{\lambda}&\\ 
\ext^2( \LG_{W_T/T \times_S V}, \pibar_T^*\JJ) \ar[rr] && \ext^1(f^* \EE,\JJ)}
\]
Since
\[
\xymatrix{T \times_S V \ar[d] \ar[r] & \tbar \times_S V\ar[d]\\
T \ar[r] & \tbar
}
\]
is a Cartesian diagram, the extension classes $[T \to \tbar]$ and $[T \times_S V \to \tbar \times_S V]$ define the same element in $\ext^1(\pi_T^*\LG_{T/S},\pi_T^*\JJ)$. Hence the element $\lambda[T \times_S V \to \tbar \times_S V] \in \ext^2(\LG_{W_T/T \times_S V}, \pibar^*_T\JJ)$ is mapped to\\
$f^*(\varphi)( \omicron[T \to \tbar ]) \in \ext^1 (f^* \EE^\bul,\JJ)$ under the isomorphism
\[
\ext^2(\LG_{W_T/T \times_S V}, \pibar^*_T\JJ) \stackrel{\cong}{\longrightarrow} \ext^1(f^*\EE^\bul,\JJ).
\]
ii) Now we want to show that
\[
f^* (\varphi) : \ext^0(f^*\LG_{\Hilb^P_{V/S}/S},\JJ) \longrightarrow \ext^0(f^*\EE^\bul,\JJ)
\]
is a bijection. By Lemma \ref{lem:obstack} there exists a bijection
\[
\ext^0(f^*\LG_{\Hilb^P_{V/S}/S},\JJ) \longrightarrow \ext^1(\LG_{W_T/T \times_S V},\pibar_T^*\JJ)
\]
which we will denote by $\xi$. Therefore it suffices to show that the following diagram commutes:
\[
\xymatrix{
\ext^0(f^*\LG_{\Hilb^P_{V/S}/S},\JJ) \ar[r]^{\xi} \ar[d]_{f^*(\varphi)}& \ext^1(\LG_{W_T/T \times_S V},\pibar_T^*\JJ) \ar[d]^{\zeta}\\
\ext^0(f^*\EE^\bul, \JJ) \ar@1{=}[r]& \ext^1(F^* \LG_{\nw/\Hilb^P_{V/S}\times_S V},\pibar_T^*\JJ)
}
\]
Denote by $\tbar_0$ the trivial square-zero extension of $T$ with ideal $\JJ$, fix an element
\[
\left( \fbar : \tbar_0 \to \Hilb^P_{V/S} \right) \in \ext^0(f^*\LG_{\Hilb^P_{V/S}/S},\JJ) = \der_{\co_S}(\co_{\Hilb^P_{V/S}},f_*\JJ),
\]
set $\wbar:=\fbar^*\nw$, and let $\bar{F}: \wbar \to \nw$ be the canonical map.

First we want to describe the image of $\fbar: \tbar_0 \to \Hilb^P_{V/S}$ under the composition
\[
\xymatrix{
\ext^0(f^*\LG_{\Hilb^P_{V/S}/S},\JJ) \ar[r]^\xi &  \ext^1(\LG_{W_T/T \times_S V},\pibar_T^*\JJ) \ar[d]^\zeta\\
& \ext^1(F^* \LG_{\nw/\Hilb^P_{V/S}\times_S V},\pibar_T^*\JJ).}
\]
Consider the following (non-commutative) diagram:
\vspace{3mm}

\[
\xymatrix{
\wbar \ar[d] \ar@/^1pc/[rr]& W_T \ar[l]\ar[r]\ar[d] & \nw \ar[r] \ar[d] & V\\
\tbar_0 \ar@<2pt>[r] & T \ar@<2pt>[l]\ar[r] & \Hilb^P_{V/S} &
}
\]
Here the map $\tbar_0 \to T$ is the projection of the trivial square-zero extension. In particular, the composition $\tbar_0 \to T \to \Hilb^p_{V/S}$ is not in general the map $\fbar$. This diagram gives $\wbar$ the structure of a square-zero extension with ideal sheaf $\pibar_T^*\JJ$ over the product $\Hilb^P_{V/S} \times_S V$; we denote its class by
\[
\alpha \in \ext^1(\LG_{W_T/\Hilb^P_{V/S} \times_S V},\pibar_T^*\JJ).
\]
Let
\[
\omicron : \ext^1(\LG_{W_T/\Hilb^P_{V/S} \times_S V},\pibar_T^*\JJ) \to \ext^1(F^* \LG_{\nw/\Hilb^P_{V/S}\times_S V},\pibar_T^*\JJ)
\]
be the canonical map. Note that the map
\[
\zeta: \ext^1(\LG_{W_T/T \times_S V},\pibar_T^*\JJ) \to \ext^1(F^* \LG_{\nw/\Hilb^P_{V/S}\times_S V},\pibar_T^*\JJ)
\]
factorizes as follows:
\[
\xymatrix{
\ext^1(\LG_{W_T/T \times_S V},\pibar_T^*\JJ) \ar[r] \ar[rd]_\zeta & \ext^1(\LG_{W_T/\Hilb^P_{V/S} \times_S V},\pibar_T^*\JJ) \ar[d]^{\omicron}\\
& \ext^1(F^* \LG_{\nw/\Hilb^P_{V/S}\times_S V},\pibar_T^*\JJ)}
\]
Therefore the image of $\fbar: \tbar_0 \to \Hilb^P_{V/S}$ under the map $\zeta \circ \xi$ is
\[
\omicron (\alpha).
\]
Let $\wbar_0:= \tbar_0 \times_T W_T$ be the trivial square-zero extension of $W_T$ with ideal $\pibar_T^*\JJ$. The following diagram
\[
\xymatrix{
\wbar_0 \ar@<2pt>[r] \ar[d]& W_T \ar@<2pt>[l]\ar[r]\ar[d] & \nw \ar[r] \ar[d] & V\\
\tbar_0 \ar@/_1pc/[rr]_{\fbar}& T \ar[l]\ar[r] & \Hilb^P_{V/S} &
}
\]
gives $\wbar_0$ the structure of a square zero-extension over $\Hilb^P_{V/S} \times_S V$. Denote its class by
\[
\beta \in \ext^1(\LG_{W_T/\Hilb^P_{V/S} \times_S V},\pibar_T^*\JJ).
\]
Using \cite[III.1.2.5.4]{ill} we see that
\[
f^*(\varphi)\left( \fbar: \tbar_0 \to \Hilb^P_{V/S} \right) = -\omicron (\beta).
\]
Therefore we have to show that
\[
\omicron(\alpha)+\omicron(\beta)=0.
\]
{\bf Claim:} The element $\alpha +\beta$ is represented by the following commutative diagram:
\[
\xymatrix{
\wbar \ar[d] \ar@/^1pc/[rr]^{\bar F}& W_T \ar[l]\ar[r]\ar[d] & \nw \ar[r] \ar[d] & V\\
\tbar_0 \ar@/_1pc/[rr]_{\fbar} & T \ar[l]\ar[r] & \Hilb^P_{V/S} &
}
\]
{\bf Proof of the claim:} Consider the structural morphisms
\[
0 \longrightarrow \pibar_T^* \JJ \stackrel{j_\alpha}{\longrightarrow} \co_{\wbar} \stackrel{p_\alpha}{\longrightarrow} \co_W \longrightarrow 0
\]
and
\[
0 \longrightarrow \pibar_T^* \JJ \stackrel{j_\beta}{\longrightarrow} \pibar_T^* \JJ \oplus \co_W   \stackrel{p_\beta}{\longrightarrow} \co_W \longrightarrow 0
\]
of the square-zero extensions $W \to \wbar$ and $W \to \wbar_0$. We define the following three maps:
\begin{eqnarray*}
q: \co_{\wbar} \oplus \left( \pibar_T^* \JJ \oplus \co_W \right) & \longrightarrow & \co_W,\\
s_\alpha \oplus s_\beta & \longmapsto & p_\alpha(s_\alpha) -p_\beta (s_\beta),
\end{eqnarray*}
\begin{eqnarray*}
i: \pibar_T^* \JJ & \longrightarrow & \co_{\wbar} \oplus \left( \pibar_T^* \JJ \oplus \co_W \right)\\
t & \longmapsto & j_\alpha (t) \oplus (-j_\beta (t)),
\end{eqnarray*}
and
\begin{eqnarray*}
r: \co_{\wbar} \oplus \left( \pibar_T^* \JJ \oplus \co_W \right) & \longrightarrow & \co_\wbar\\
s_1 \oplus ( t\oplus s_2)& \longmapsto & s_1+ j_\alpha(t).
\end{eqnarray*}
Note that the natural maps
\[
\co_{\Hilb^P_{V/S}} \longrightarrow \co_{\wbar} \oplus \left( \pibar_T^* \JJ \oplus \co_W \right)
\]
and
\[
\co_\nw \longrightarrow \co_{\wbar} \oplus \left( \pibar_T^* \JJ \oplus \co_W \right)
\]
factor through the inclusion $\ker (q) \hookrightarrow \co_{\wbar} \oplus \left( \pibar_T^* \JJ \oplus \co_W \right)$ and hence define morphisms
\[
\co_{\Hilb^P_{V/S}} \longrightarrow \ker (q)/ \im (i)
\]
and
\[
\co_\nw \longrightarrow \ker (q)/ \im (i).
\]
With these maps the extension representing $\alpha + \beta$ is given by the diagram
\[
\xymatrix{
& \co_{\Hilb^P_{V/S}} \ar[d] &\\
\pibar_T^* \JJ \ar[r] & \ker (q) / \im (i) \ar[r] & \co_W.\\
& \co_\nw \ar[u] &}
\]
The composition
\[
\ker (q) \hookrightarrow \co_{\wbar} \oplus \left( \pibar_T^* \JJ \oplus \co_W \right) \stackrel{r}{\longrightarrow} \co_\wbar
\]
induces an isomorphism of square-zero extensions
\[
\ker (q) / \im (i) \stackrel{\cong}{\longrightarrow} \co_\wbar.
\]
Using this isomorphism we obtain maps $\wbar \to \nw$ and $\wbar \to \Hilb^P_{V/S}$ and one checks that they coincide with $\bar{F}: \wbar \to \nw$ and the composition $\wbar \to \tbar \stackrel{\fbar}{\longrightarrow} \Hilb^P_{V/S}$. This proves our claim.

Since $F:W_T \to \nw$ extends to the morphism $\bar{F}: \wbar \to \nw$ over $\Hilb^P_{V/S} \times_S V$, we find
\[
\omicron(\alpha+\beta)=0.
\]
\end{proof}
\begin{lem}
Let $v: V \to S$ be a projective morphism, and let $g: S' \to S$ be a base change. Set $V':=V \times_S S'$. Then there exists a Cartesian diagram:
\[
\xymatrix{
\Hilb^P_{V'/S'} \ar[d] \ar[r] & \Hilb^P_{V/S} \ar[d]\\
S' \ar[r]^g & S}
\]
\end{lem}
\begin{proof}
Composing the forgetful functor
\[
(\operatorname{Schemes}/S') \longrightarrow (\operatorname{Schemes}/S)
\]
with the functor
\[
\Hifu^P_{V/S}: (\operatorname{Schemes}/S) \longrightarrow (\operatorname{Sets})^\circ
\]
yields the functor
\[
\Hifu^P_{V'/S'}: (\operatorname{Schemes}/S') \longrightarrow (\operatorname{Sets})^\circ.
\]
Our claim follows immediately.
\end{proof}
\begin{prop}\label{prop:obthbc}
Let $v: V \to S$ be a flat projective morphism, and let $g: S' \to S$ be a base change. Set $V':=V \times_S S'$. Denote by $\varphi :\mathcal{E}^{\bullet} \longrightarrow \LG_{\Hilb^P_{V/S}/S}$ and $\varphi' :\mathcal{E'}^{\bullet} \longrightarrow \LG_{\Hilb^P_{V'/S'}/S'}$ the relative obstruction theories, and by $g': \Hilb^P_{V'/S'} \to \Hilb^P_{V/S}$ the induced morphism between the Hilbert schemes. Then there exists an isomorphism
${g'}^*
\mathcal{E}^{\bullet}
\stackrel{\cong}{\longrightarrow} \mathcal{E'}^{\bullet}$ such that the
following square commutes:
\[
\xymatrix{
{g'}^* \mathcal{E}^{\bullet} \ar[d]^{\cong} \ar[r]^-{{g'}^* \varphi} &
{g'}^* \LG_{\Hilb^P_{V/S}/S} \ar[d]\\
\mathcal{E'}^{\bullet} \ar[r]^-{\varphi'} & \LG_{\Hilb^P_{V'/S'}/S'} . }
\]
\end{prop}
\begin{proof}
Let $g_V: V'\to V$ be the natural map. Since $v: V \to S$ is flat, the morphism
\[
g_V^* \LG_{V/S} \to \LG_{V'/S'}
\]
is an isomorphism. Hence our claim follows from the functoriality of the functor $L\pibar_\sharp$ \cite[Satz 2.1]{flenner}.
\end{proof}
\subsection{Hilbert schemes of divisors on smooth projective varieties}
Next we want to give a description of $\mathcal{E}^{\bullet}$ in more accessible terms in the case where we are looking at divisors instead of general subschemes. In order to do this, we need the following lemma, which was suggested by Flenner.
\begin{lem} Let $h: M \longrightarrow N$ be a flat proper morphism of
schemes, and assume that $N$ has a dualizing complex. If $h$ is Gorenstein of relative dimension $d$, then there exists for any object $\ff
\in D^-_c (\mathcal{O}_M)$ an isomorphism
\[
Lh_\sharp (\mathcal{F}^{\bullet}) \cong Rh_* (\mathcal{F}^{\bullet}
\otimes \omega_h [d]), 
\]
where $\omega_h$ is the relative dualizing sheaf of $h$. 
\end{lem}
\begin{proof}
Fix a dualizing complex $\KG^\bul_N$ on $N$. By the explicit description of the functor $Lh_\sharp$ given in \cite{flenner}, we have
\begin{eqnarray*}
Lh_\sharp (\FF^\bul) &=& R\HH om_N (Rh_* R\HH om_M( \FF^\bul, h^* \KG_N^\bul), \KG_N^\bul )\\
&\cong& R\HH om_N (Rh_* R\HH om_M( \FF^\bul \otimes \omega_h[d], h^* \KG_N^\bul\otimes \omega_h[d]), \KG_N^\bul ).
\end{eqnarray*}
An application of relative duality \cite[III Thm.11.1]{Ha2} yields
\begin{eqnarray*}
Lh_\sharp (\FF^\bul) &\cong &  R\HH om_N (R\HH om_N (Rh_* 
               (\FF^\bul \otimes \omega_h [n]) , \KG_N^\bul), \KG_N^\bul) \\
     &\cong & Rh_* (\FF^\bul \otimes \omega_h [n]).
\end{eqnarray*}
\end{proof}
\begin{thm}\label{thm:simple}
Let $v:V \longrightarrow S$ be a smooth projective morphism of relative dimension $d$. Fix a polynomial $P$ such that $\Hilb^P_{V/S}$
parametrizes divisors. Denote by $\omega_{V/S}$ the relative dualizing
sheaf of $V$ over $S$ and by $\mathbb{D}$ the universal divisor on
$\Hilb^P_{V/S} \times_S V$. Let $pr: \Hilb^P_{V/S}\times_S V \to V$ and $\pi:\Hilb^P_{V/S}\times_S V \to \Hilb^P_{V/S} $ be the projections, let $i: \mathbb{D} \to \Hilb^P_{V/S}\times_S V$ be the inclusion, and set $\pibar: = \pi \circ i$. Then there are isomorphisms
\begin{eqnarray*}
\mathcal{E}^{\bullet} & \cong & ( R^{\bullet}\pi_* \mathcal{O}_\mathbb{D} (
\mathbb{D}) )^{\vee}, \\
\mathcal{E}^{\bullet} & \cong & R^{\bullet} \pi_* \mathcal{H}om ( \mathcal{O}_\mathbb{D} , \mathcal{O}(-\mathbb{D}) \otimes pr^* \omega_{V/S} )[d],
\end{eqnarray*}
and
\[
\mathcal{E}^{\bullet} \cong R^{\bullet}\bar{\pi}_* (i^* pr^* \omega_{V/S} ) [d-1].
\]
\end{thm}
\begin{proof}
Since $i$ is a regular embedding, we have
\[
\LG_{\nd/\Hilb^P_{V/S} \times_S S}[-1]\cong \co_\nd(-\nd).
\]
Hence, the previous lemma implies
\[
\ee \cong R\pibar_* \left( \omega_{\nd/S}\otimes \co_\nd(-\nd)\right)[d-1],
\]
which, by relative duality, yields
\[
\ee \cong \left(R^\bullet \pibar_* \co_\nd (\nd)\right)^\vee,
\]
or equivalently
\begin{equation}\label{eq:ob}
\ee \cong \left(R^\bullet \pi_* \co_\nd (\nd)\right)^\vee.
\end{equation}
On the other hand, we have
\[
\omega_{\nd/S} \cong i^* pr^* \omega_{V/S} \otimes \co_\nd (\nd ),
\]
which shows
\[
\ee \cong R^\bullet \pibar_* \left( i^* pr^* \omega_{V/S} \right) [d-1].
\]
By applying relative duality with respect to the projection $\pi$ to equation \eqref{eq:ob}, we obtain isomorphisms
\[
\ee \cong R^\bullet \pi_*\HH om \left( \co_\nd (\nd), pr^*\omega_{V/S} \right) [d]
\]
and
\[
\ee \cong R^\bullet \pi_*\HH om \left( \co_\nd , \co (-\nd) \otimes pr^*\omega_{V/S} \right) [d].
\]
\end{proof}
Let $v:V \longrightarrow S$ be a smooth projective morphism of relative dimension $d$, and let $P$ be a polynomial such that $\Hilb^P_{V/S}$
parametrizes divisors. Then, by Thm.~\ref{thm:ob} and Thm.~\ref{thm:simple}, we
obtain an isomorphism
\[
\Omega_{\Hilb^P_{V/S}/S} \cong \mathcal{E}xt_\pi^d (\mathcal{O}_\mathbb{D},
\mathcal{O}(-\mathbb{D}) \otimes pr^* \omega_{V/S} ).
\]
This is a special case of Lehn's description of the cotangent sheaves of Quot-schemes \cite[Thm 3.1.]{lehn}, since
$\Hilb^P_{V/S}$ is a Quot-scheme with universal object
\[
0 \longrightarrow \mathcal{O}(-\mathbb{D}) \longrightarrow \mathcal{O}
\longrightarrow \mathcal{O}_\mathbb{D} \longrightarrow 0.
\]

\begin{prop}
Let $v:V \longrightarrow S$ be a smooth projective morphism of relative dimension $d$. Fix a polynomial $P$ such that the relative Hilbert scheme $\Hilb^P_{V/S}$ parametrizes divisors. Let $k$ be an integer such that for any point $D \in \Hilb^P_{V/S}$ and all $i>k$ we have $H^i(\co_D(D)) =0$. Then for each $p \in S$ there exists an open neighbourhood $U \subset S$ such that $\ee|_U$ has a global resolution of perfect amplitude contained in $[-k,0]$.
\end{prop}
\begin{proof}
First we show that the complex $R^\bul\pi_* \co_\nd (\nd)$ has, locally with respect to the base scheme $S$, a global resolution of perfect amplitude contained in $[0,d-1]$. 

When $d=1$, the higher direct image sheaves $R^i\pi_* \co_\nd(\nd)$ vanish for $i \geq 1$, and $R^\bul\pi_* \co_\nd(\nd) \cong \pi_* \co_\nd(\nd)$, considered as a complex concentrated in degree $0$. Moreover, the sheaf $\pi_* \co_\nd(\nd)$ is locally free \cite[Thm.III.12.11]{Ha1}. Suppose now $d>1$ and fix a point $p \in S$. Let $\ov(1)$ be a relatively ample sheaf. By upper semicontinuity, there exists an $l>>0$ and an open subset $U' \subset S$ containing $p$ such that the following conditions are satisfied:
\begin{itemize}
\item for all $i>0$ and for all $D \in \Hilb^P_{V_{U'}/U'}$ we have $H^i(\co_D(D)(l))=0$;
\item for all $D \in \Hilb^P_{V_{U'}/U'}$ we have $H^0(\co(D)(-l))=0$.
\end{itemize}
Without loss of generality, we may assume that $\ov(l)$ is relatively very ample. Let $j: V \hookrightarrow U' \times \np^n$ be the corresponding embedding. By Bertini's theorem, there exists a hyperplane $H \subset \np^n$ such that $H \cap V_p$ is a smooth connected divisor. Since $v:V \to S$ is smooth and proper, there is an open subset $U'' \subset U' \subset S$ containing the point $p$ such that for all $p'\in U''$ the intersection $H \cap V_{p'}$ is a smooth connected divisor in $V_{p'}$.

Set $\nh:= \Hilb^P_{V_{U''}/U''} \times_{U''}((U'' \times H) \cap V_{U''})$. Since for all $D \in \Hilb^P_{V_{U''}/U''}$ we have $H^0(\co(D)(-l))=0$, the intersection $\nd_{U''} \cap \nh \subset \nh$ is a divisor, flat over $\Hilb^P_{V_{U''}/U''}$. This implies that the following sequence is exact:
\[
0 \to \co_{\nd_{U''}}(\nd_{U''}) \to \co_{\nd_{U''}}(\nd_{U''}+\nh) \to \co_{\nd_{U''} \cap \nh} (\nd_{U''}+\nh) \to 0
\]
The sheaf $\pi_*\co_{\nd_{U''}}(\nd_{U''}+\nh)$ is locally free, while the higher direct image sheaves $R^i\pi_*\co_{\nd_{U''}}(\nd_{U''}+\nh)$ vanish for $i\geq 1$.
Note that the sheaf $\co_{\nd_{U''} \cap \nh} (\nd_{U''}+\nh)$ is the restriction of the invertible sheaf $\co_{\nh} (\nd_{U''}+\nh)$ to the divisor $\nd_{U''} \cap \nh \subset \nh$ which is flat over $\Hilb^P_{V_{U''}/U''}$. By recurrence, we find that there is an open subset $U \subset U''$ containing the point $p$ and a global resolution
\[
(R^\bul\pi_* \co_{\nd_{U''} \cap \nh})|_U \cong \FF^1 \to \ldots \to \FF^{d-1}.
\]
Then
\[
\FF^0:=\pi_*\co_{\nd_{U''}}(\nd_{U''}+\nh)|_U \to \FF^1 \to \ldots \to \FF^{d-1}
\]
is a global resolution of the complex $(R^\bul \pi_* \co_\nd(\nd))|_U$ by locally free sheaves.

Suppose now that for some $k$ we have $H^i(\co_D(D))=0$ for all $i>k$ and for all $D\in \Hilb^P_{V/S}$. Then the direct image sheaves $R^i \pi_* \co_\nd(\nd)$ vanish for all $i>k$. Therefore the kernel of the map $\delta^k: \FF^k \to \FF^{k+1}$ is locally free, and the complex
\[
\FF^0 \to \ldots \to F^{k-1} \to \ker \delta^k
\]
is a global resolution of $(R^\bul \pi_* \co_\nd (\nd))|_U$ by locally free sheaves
\end{proof}
\section{Curves on surfaces}
In this section, all surfaces will be smooth, projective, connected, and defined over $\mathbb{C}$.
\subsection{Virtual fundamental classes for Hilbert schemes of curves on surfaces}
\begin{dfn}
Let $V \to S$ be a smooth family of surfaces, and suppose that $S$ is connected and of pure dimension. Fix a class $\underline{m} \in H^0(S, R^2v_*\underline{\nz})$. Then 
\[
\lb \Hilb^{\underline{m}}_{V/S} \rb \in A_*(\Hilb^{\underline m}_{V/S})
\]
is the virtual fundamental class defined by the obstruction theory
\[
\varphi : \ee \longrightarrow \LG_{\Hilb^{\underline m}_{V/S}/S}.
\]
If $S =\Spec \nc$ and $m\in H^2(V,\nz )$, then we denote by $\lb \Hilb^m_V \rb: = \lb \Hilb^{\underline m}_{V/S}\rb $ the virtual fundamental class of the Hilbert scheme $\Hilb^m_V$.
\end{dfn}
Note that
\[
\lb \Hilb^m_V \rb \in A_{\frac{m(m-k)}{2}} (\Hilb^m_V).
\]
\begin{prop}
Let $V \to S$ be a smooth family of surfaces, and suppose that $S$ is connected and of pure dimension. Fix a class $\underline{m} \in H^0(S, R^2v_*\underline{\nz})$. Let $S'$ be a another connected scheme of pure dimension, and fix a morphism $j:S' \to S$. Set $V' :=V \times_S S'$ and $\underline{m}':= j^*\underline{m}$. If $j$ is flat or a regular local immersion, then
\[
\lb \Hilb^{\underline{m}'}_{V'/S'} \rb = j^! \lb \Hilb^{\underline m}_{V/S} \rb.
\]
\end{prop}
\begin{proof}
This is a direct consequence of Prop.~\ref{prop:bc} and Prop.~\ref{prop:obthbc}.
\end{proof}
The following simple corollary will be of particular interest to us:
\begin{cor}\label{cor:basechange}
Let $V \to S$ be a smooth family of surfaces, and suppose that $S$ is smooth and connected. Fix a class $\underline{m} \in H^0(S, R^2v_*\underline{\nz})$. Let $s \in S$ be a point, and denote by $j_s: \{ s \} \to S$ the inclusion. Then
\[
\lb \Hilb^{\underline{m}(s)}_{V_s} \rb = j_s^! \lb \Hilb^{\underline{m}}_{V/S} \rb .
\]
\end{cor}
\begin{proof}
Since $S$ is smooth, the embedding $j_s$ is regular.
\end{proof}
\subsection{A second obstruction theory on projective surfaces}
Let $V$ be a surface, and let $m \in H^2(V,\mathbb{Z})$. Fix an effective Cartier divisor $H \subset V$, and set $h:=c_1(\ov(H))$. Denote by $\nd$ the universal divisor on $\Hilb^m_V \times V$ and by $\ndtil$ the universal divisor on $\Hilb^{m+h}_V \times V$. Put $\nh:= \Hilb^m_V \times H \subset \Hilb^m_V \times V$ and $\nhtil := \Hilb^{m+h}_V \times H \subset \Hilb^{m+h}_V \times V$. Let $\pi$ be the projection $\Hilb^m_V \times V \to \Hilb^m_V$, and let $\tilde{\pi}$ be the projection $\Hilb^{m+h}_V \times V \to \Hilb^{m+h}_V$.

By adding the fixed divisor $H$ to a divisor $D \in \Hilb^m_V$,  we obtain an inclusion of schemes $j: \Hilb^m_V \hookrightarrow \Hilb^{m+h}_V$. Composing the inclusion $\OO \to \OO(\ndtil)$ with the restriction map $\OO(\ndtil) \to \OO_\nhtil(\ndtil)$ defines a global section $s \in H^0(\pitil_* \OO_\nhtil(\ndtil)$. Let $p: \Hilb^{m+h}_V \times V \to V$ be the projection onto $V$. By relative duality, we have an isomorphism
\[
\pitil_* \OO_\nhtil (\ndtil) \stackrel{\cong}{\longrightarrow} \left( R^{1} \pitil_* p^*\kv (\nhtil -\ndtil) |_\nhtil \right)^\vee.
\] 
\begin{prop}\label{zeroset}
Let $V$ be a surface. The morphism
\begin{eqnarray*}
j: \Hilb^m_V & \longrightarrow & \Hilb^{m+h}_V\\
D & \longmapsto & D+H
\end{eqnarray*}
is a closed embedding, and the following sequence is exact:
\[
R^{1} \pitil_* p^*\kv (\nhtil -\ndtil) |_\nhtil \stackrel{s}{\longrightarrow} \OO_{\Hilb^{m+h}_V} \longrightarrow j_* \OO_{\Hilb^m_V} \longrightarrow 0.
\]
\end{prop}
\begin{proof}
Let $S$ be a scheme, and let $\tilde{D}$ be a divisor on $S \times V$ corresponding to a morphism $\alpha: S \to \Hilb^{m+h}_V$. Put $H_S:= S \times H$ and denote by $\pi_S$ the projection $S \times V \to S$.

Since the sheaf $R^{1} \pitil_* p^*\kv (\nhtil -\ndtil) |_\nhtil$ has the base change property \cite[III, Thm.12.11]{Ha1}, the pull-back of the section
\[
s: R^{1} \pitil_* p^*\kv (\nhtil -\ndtil) |_\nhtil \to \OO
\]
corresponds to the push-forward by $\pi_S$ of the composition of 
\[
\OO \to \OO(\tilde{D})
\]
with the restriction map
\[
\OO(\tilde{D}) \to \OO_{H_S}(\tilde{D}).
\]
So we see that the morphism $\alpha$ factors through the inclusion $Z(s) \hookrightarrow \Hilb^{m+h}_V$ iff the morphism of sheaves $\OO \to \OO (\tilde{D})$ factors through $\OO (\tilde{D}-H_S ) \to \OO (\tilde{D})$. This proves our claims.
\end{proof}
Denote by $k$ the first Chern class of the canonical line bundle $\kv$.
\begin{lem}
Let $V$ be a surface, and let $m \in H^2(V,\mathbb{Z})$. The following conditions are equivalent:
\begin{itemize}
	\item[i)] $H^2(\ov(D))=0\ \ \forall \, D \in \Hilb^m_V$;
	\item[ii)] $R^2\pi_* \OO(\nd)=0$;
	\item[iii)] The fibered product $\Hilb^m_V \times_{\Pic^m_V} \Hilb^{k-m}_V$ is empty.
\end{itemize}
\end{lem}
\begin{proof}
The equivalence of the first two statements follows from \cite[III, Thm.12.11]{Ha1}, while the first and the third statement are equivalent by Serre duality.
\end{proof}
\begin{lem}
Let $V$ be a surface, and fix $m\in H^2(V,\mathbb{Z})$. There exists a smooth effective divisor $H \subset V$ such that $H^i(\LL(H))=0$ for each $[\LL] \in \Pic^m_V$ and for all $i>0$.
\end{lem}
\begin{proof}
Choose an effective ample divisor $E$ on $V$. Upper-semicontinuity of cohomology implies that for $l$ large enough we have
\[
h^i(\mathcal{L}(lE))=0\ \ \forall [\mathcal{L}] \in \Pic^m_V,\ i=1,2.
\]
Moreover, if $lE$ is very ample, Bertini's theorem implies that a general element of the total linear system of $lE$ is a smooth curve on $V$.
\end{proof}
\begin{lem}\label{lem:free+bcp}
Let $V$ be a surface, fix $m \in H^2(V,\mathbb{Z})$, and let $H \subset V$ be a smooth effective divisor such that $H^i(\LL(H))=0$ for each $[\LL] \in \Pic^m_V$ and for all $i>0$. Set $h:= c_1(\ov (H))$, denote by $\ndtil$ the universal divisor over $\Hilb^{m+h}_V$, by $\nhtil$ the divisor $\Hilb^{m+h}_V \times H$, and let $\pitil : \Hilb^{m+h}_V \times V \to \Hilb^{m+h}_V$ be the projection. If the fibered product $\Hilb^m_V \times_{\Pic^m_V} \Hilb^{k-m}_V$ is empty, then the exists a Zariski open neighbourhood $U$ of $\Hilb^m_V \subset \Hilb^{m+h}_V$ such that the restriction $R^1\pitil_* \co_\nhtil(\ndtil)|_U$ vanishes.
\end{lem}
\begin{proof}
Let $D \in \Hilb^m_V$ be a divisor. Then $H^1(\ov(D+H))$ vanishes by our assumption on $H$, while $H^2(\ov(D))$ vanishes since $\Hilb^m_V \times_{\Pic^m_V} \Hilb^{k-m}_V$ is empty. Hence we have $H^1(\co_H(D+H))=0$, and our claim follows by upper semicontinuity.
\end{proof}

The short exact sequence
\[
0 \to \OO \to \OO (\nd) \to \OO_\nd (\nd ) \to 0
\]
induces a morphism $R^\bullet \pi_*\OO_\nd (\nd ) \to R^\bullet \pi_* \OO [1]$. We denote by $\chi$ its composition with the cut-off $R^\bullet\pi_* \OO [1] \to (\sigma_{\geq 2} R^\bullet \pi_*\OO)[1]$, and define $\CC^\bullet$ to be the mapping cone of $\chi$, shifted by $-1$. Then the following triangle is distinguished:
\[
\xymatrix{
\sigma_{\geq 2} R^\bullet \pi_*\OO \ar[r] & \CC^\bullet \ar[d]\\
& R^\bullet\pi_*\OO_\nd (\nd ) \ar[ul]_\chi^{[1]}}
\]
\begin{prop}
Let $V$ be a surface, and let $m \in H^2(V,\mathbb{Z})$. The complex $\CC^\bullet$ has a global resolution of perfect amplitude contained in $[0,2]$. If the fibered product $\Hilb^m_V \times_{\Pic^m_V} \Hilb^{k-m}_V$ is empty, then there exists an obstruction theory $\phi : R \mathcal{H}om (\CC^\bullet ,\OO) \to \LG_{\Hilb^m_V}$ and the complex $R \mathcal{H}om (\CC^\bullet ,\OO)$ admits a perfect global resolution.
\end{prop}
\begin{proof}
The complex $R^\bul \pi_* \co_\nd (\nd )$ has a global resolution of perfect amplitude contained in $[0,1]$, while $\sigma_{\geq 2}R^\bul \pi_* \co$ can be represented by a locally free sheaf in degree 2. This proves the first claim.

Suppose now that the fibered product $\Hilb^m_V \times_{\Pic^m_V} \Hilb^{k-m}_V$ is empty, and choose a smooth effective divisor $H \subset V$ such that $H^i(\LL(H))=0$ for each $[\LL] \in \Pic^m_V$ and for all $i>0$. Consider the following short exact sequence of sheaves on $\Hilb^m_V$:
\[
0  \to \co_\nd (\nd) \to \co_{\nd +\nh}(\nd + \nh) \to \co_\nh (\nd + \nh) \to 0
\]
By our assumption on the divisor $H$, the Hilbert scheme $\Hilb^{m+h}_V$ is smooth, where $h:=c_1(\ov(H))$. In particular, the sheaf $\pi_*\co_{\nd +\nh}(\nd + \nh)$ is locally free and has the base change property. By Lemma \ref{lem:free+bcp} there exists a Zariski open neighbourhood $U$ of $\Hilb^m_V \subset \Hilb^{m+h}_V$ such that
\[
\left(R^1\pitil_*\left ( (p^*\kv (\nhtil-\ndtil))|_\nhtil \right) \right)|_U
\]
is locally free. Using Lemma \ref{zeroset}, we see that $\Hilb^m_V \subset U$ is the zero locus of a section in a vector bundle. Moreover, $\pi_* \co_\nh (\nd + \nh)$ is locally free, has the base change property, and represents the complex $R^\bul \pi_* \co_\nh (\nd + \nh)$. Therefore, there is an obstruction theory
\[
\left( \pi_* \co_{\nd +\nh}(\nd + \nh) \to \pi_*\co_\nh (\nd + \nh) \right)^\vee \longrightarrow \LG_{\Hilb^m_V}.
\]

In order to prove the remaining two assertion, we have to show that
\[
\pi_* \co_{\nd +\nh}(\nd + \nh) \to \pi_*\co_\nh (\nd + \nh)
\]
is a global resolution of the complex $\CC^\bul$.

Let
\[
\nu : R^\bul \pi_* \co_{\nd +\nh}(\nd + \nh) \to \sigma_{\geq 2} R^\bul \pi_* \co[1]
\]
be the composition of the truncation morphism
\[
R^\bul \pi_* \co_{\nd +\nh}(\nd + \nh) \to \sigma_{\geq 1} R^\bul \pi_* \co_{\nd +\nh}(\nd + \nh)\]
 with the isomorphism $\sigma_{\geq 1} R^\bul \pi_* \co_{\nd +\nh}(\nd + \nh) \to \sigma_{\geq 2 } R^\bul \pi_* \co [1]$. We are in the following situation:
\[
\xymatrix{
\CC^\bul \ar@1{.>}[r] \ar[d] & \pi_* \co_{\nd +\nh}(\nd + \nh) \ar[d] \ar[r]& \pi_* \co_\nh (\nd + \nh ) \ar@1{=}[d] \ar@1{.>}[r]& \CC^\bul [1]\\
R^\bul \pi_* \co_\nd (\nd ) \ar[r]\ar[d] & R^\bul \pi_* \co_{\nd +\nh}(\nd + \nh) \ar[d]^\nu \ar[r]& R^\bul \pi_* \co_\nh (\nd + \nh ) \ar[r] & R^\bul \pi_* \co_\nd (\nd )[1]\\
\sigma_{\geq 2} R^\bul \pi_* \co [1] \ar@1{=}[r] \ar[d]& \sigma_{\geq 2} R^\bul \pi_* \co [1]  \ar[d] &&\\
\CC^\bul[1]& \pi_* \co_{\nd +\nh}(\nd + \nh)[1]&&
}
\]
Here, the dotted arrows exist by the octrahedral axiom \cite[p.21]{Ha2}, and our claim follows.
\end{proof}
\begin{dfn}
Let $V$ be a surface, fix a class $m \in H^2(V,\mathbb{Z})$, and suppose that the fibered product $\Hilb^m_V \times_{\Pic^m_V} \Hilb^{k-m}_V$ is empty. Then
\[
\{ \Hilb^m_V \} \in A_*(\Hilb^m_V)
\]
is the virtual fundamental class defined by the obstruction theory
\[
\phi : R \mathcal{H}om (\CC^\bullet ,\OO) \to \LG_{\Hilb^m_V}.
\]
\end{dfn}
Note that
\[
\{ \Hilb^m_V \} \in A_{\frac{m(m-k)}{2}+p_g(V)}(\Hilb^m_V).
\]
\begin{thm}\label{thm:compobth}
Let $V$ be a surface, and let $m \in H^2(V,\mathbb{Z})$ be a class with $m(m-k) \geq 0$. Then the following holds:
\begin{itemize}
	\item[i)] If $p_g(V)=0$, then $[[\Hilb^m_V]]=\{ \Hilb^m_V\}$;
	\item[ii)] If $p_g(V)>0$ and the fibered product $\Hilb^m_V \times_{\Pic^m_V} \Hilb^{k-m}_V$ is empty, then $[[\Hilb^m_V]]=0$.
\end{itemize}
\end{thm}
\begin{proof}
Let $\FF_0 \to \FF_1$ be a global resolution of $R^\bullet \pi_*\OO_\nd (\nd)$. We set:
\begin{eqnarray*}
\GG_0 &:=& \FF_0\\
\GG_1 &:=& \ker (\FF_1 \to R^1\pi_*\OO_\nd (\nd) \to R^2\pi_*\OO)
\end{eqnarray*}
Note that $\GG_1$ is a sub-vector bundle of $\FF_1$ since it is the kernel of a surjective morphism of locally free sheaves. Moreover, the map $\FF_0 \to \FF_1$ factors through $\GG_1$. Therefore $\GG_0 \to \GG_1$ is a global resolution of $\CC^\bullet$, and the morphism $\CC^\bullet \to R^\bullet \pi_*\OO_\nd (\nd)$ is represented by
\[
\xymatrix{
\GG_0 \ar@{=}[d] \ar[r] & \GG_1 \ar[d]\\
\FF_0 \ar[r] & \FF_1 .}
\]
Let $F_0$, $F_1$, $G_0$, and $G_1$ denote the corresponding vector bundles, and let $\varphi'$ be the composition
\[
\varphi': \EE^\bullet \longrightarrow R \mathcal{H}om (\CC^\bullet,\OO) \longrightarrow \LG_{\Hilb^m_V}.
\]
One has the following diagram, where both squares are Cartesian
\[
\xymatrix{
C \ar[d] \ar[r] & G_1 \ar[r] \ar[d] & F_1\ar[d]\\
\mathfrak{C}_{\Hilb^m_V} \ar[r] & G_1/G_0 \ar[r] & F_1/F_0.}
\]
Here $\mathfrak{C}_{\Hilb^m_V}$ is the intrinsic normal cone of the Hilbert scheme, and $C$ is the closed subcone determined by the obstruction theories. We apply Prop.~\ref{prop:besie} and Prop.~\ref{prop:ex} to conclude
\begin{eqnarray*}
[[\Hilb^m_V]] &=& [\Hilb^m_V,\varphi']\\
&=& c_{top}(R^2\pi_*\OO) \cap \{ \Hilb^m_V\}.
\end{eqnarray*}
Since $R^2\pi_*\OO$ is a locally free sheaf of rank $p_g$, our claims follow.
\end{proof}
\subsection{A Porteous' formula}
In this subsection, we prove a Porteous type formula for Hilbert schemes of curves on surfaces. This formula will later play a role in the proof of the wall crossing formula, but is of independent interest. We state our formula in terms of a modified Segre class.

If $E$ is a vector bundle on a scheme, we denote by $P(E)$ the associated projective fiber space in the sense of Fulton, i.e.~$P(E):= \np (E^\vee)$.
\begin{dfn} \label{def:segre}
  Let $V$ be a smooth proper scheme of dimension $d$, and let $[E-F] \in K^0(V)$ be a virtual vector bundle on $V$. The \em{modified Segre class} of $[E-F]$ is
\[
\shat ([E-F]) := \sum_{j=0}^{\min(d,d-1 + \rk([E-F]))} c_{d-j} ([F-E]) \cap [V]. 
\]
\end{dfn}

\begin{rem}
  If $F=0$, then $\shat ([E \oplus 1 -F]) = s(E)$ is the standard Segre
  class \cite[4.1]{Fu}. 
\end{rem}

\begin{prop} \label{prop:segre}
  Let $V$ be a smooth proper scheme of dimension $d$, and let $E$, $F$ be
  vector bundles on $V$. Denote by $\nu: P(E) \to V$ the projection,
  and put $u:= c_1 (\OO_{P(E)}(1))$. Then
\[
\shat ([E-F]) = \nu_* \left( \sum_i u^{i} \cap ( c_{top}( \OO_{P(E)}(1)
  \otimes \nu^*F ) \cap [P(E)] )\right) .
\]
\end{prop}

\begin{proof}
This is a non-classical version of Porteous' formula. For a proof, see \cite[Thm.14.4]{Fu}.
\end{proof}

\begin{cor} \label{cor:segre} 
  Let $\varphi : E \to F$ be a morphism of vector bundles over a
  smooth proper scheme $V$ of dimension $d$. Let $\nu: P(E) \to V $ be
  the projection, and denote by $\tilde{\varphi}$ the induced section in the
  bundle $\OO_{P(E)}(1) \otimes \nu^* F$. Let $\iota: Z(\tilde{\varphi})
  \to P(E)$ be the embedding of the zero scheme of $\tilde{\varphi}$. Then
\[
\nu_* \left( \sum_i u^{i} \cap \iota_* [[ Z(\tilde{\varphi}) ]]
 \right) = \shat ([E-F]). 
\]
\end{cor}
\begin{proof}
Since $i_* [[Z(\tilde{\varphi})]]=c_{top}( \OO_{P(E)}(1) \otimes \nu^*F ) \cap [P(E)]$, this is an immediate consequence of the previous lemma.
\end{proof}
Fix a Poincar\'e line bundle $\nl$ on $\Pic^m_V \times V$, and let $\mu: \Pic^m_V \times V \to \Pic^m_V$ and $pr_V: \Pic^m_V \times V \to V$ be the projections. Consider the projective fibration
\[
\rhotil_+: \np \left( R^2\mu_* (\nl^\vee \otimes pr_V^*\kv) \right) \to \Pic^m_V.
\]
The canonical epimorphism
\[
\rhotil_+^* (R^2\mu_* (\nl^\vee \otimes pr_V^*\kv)) \to \co(1)
\]
defines a section $\Phi$ in the line bundle
\[
(\rhotil_+ \times id_V)^* (\nl) \otimes \pitil^*\co(1)
\]
on $\np \left( R^2\mu_*( \nl^\vee \otimes pr_V^*\kv) \right) \times V$. Here $\pitil$ is the projection
\[
\pitil : \np \left( R^2\mu_* (\nl^\vee \otimes pr_V^*\kv )\right) \times V \to \np \left( R^2\mu_* (\nl^\vee \otimes pr_V^*\kv) \right).
\]
The vanishing locus $D^+$ of $\Phi$ is a divisor, flat over $\np \left( R^2\mu_* (\nl^\vee \otimes pr_V^*\kv) \right)$. Analogously, we obtain a divisor $D^-$ in $\np \left( R^2\mu_* \nl \right) \times V$. 
\begin{lem}\label{lem:finmod}
The pairs
\[
\left( \np \left( R^2\mu_*( \nl^\vee \otimes pr_V^*\kv) \right), D^+\right) \text{ and } \left(\np \left( R^2\mu_* \nl \right),D^- \right)
\]
represent the functors $\Hifu^m_V$ and $\Hifu^{k-m}_V$. When the Poincar\'e line bundle is normalized, i.e.~when $\nl|_{\Pic^m_V \times \{ p \}} \cong \co_{\Pic^m_V}$ for some point $p \in V$, then
\[
\co_{\np \left( R^2\mu_*( \nl^\vee \otimes pr_V^*\kv) \right)}(1) \cong \co(D^+)|_{\np \left( R^2\mu_*( \nl^\vee \otimes pr_V^*\kv) \right) \times \{ p \} }
\]
and
\[
\co_{\np \left( R^2\mu_* \nl \right)}(1) \cong \co(D^-)_{\np \left( R^2\mu_* \nl \right) \times \{ p \}}.
\]
\end{lem}
\begin{proof}
Let $S$ be an arbitrary scheme, and fix a morphism $\varphi : S \to \Hilb^m_V$. Denote the corresponding divisor on $S \times V$ by $D_\varphi$, and set $\psi:=\rho_+ \circ \varphi$, where $\rho_+: \Hilb^m_V \to \Pic^m_V$ is the map which sends a divisor $D$ to the class of its associated line bundle $[\ov(D)]$. By the universal property of the Picard scheme $\Pic^m_V$, there exists a line bundle $T$ on $S$ and an isomorphism
\[
\co(D_\varphi) \stackrel{\cong}{\longrightarrow} (\psi \times Id_V)^* \nl \otimes pr_S^*T.
\]
Form the composition
\[
s:  \co \to \co(D_\varphi) \to (\psi \times Id_V)^* \nl \otimes pr_S^*T.
\]
By relative duality, the section $s$ corresponds to a morphism
\[
\epsilon: \psi^* (R^2\mu_*(\nl^\vee \otimes pr_V^* \kv)) \to T
\]
Moreover, since $D_\varphi$ is flat over $S$, the morphism $\epsilon$ is surjective, and hence defines a map $S \to \np \left( R^2\mu_* (\nl^\vee \otimes pr_V^*\kv) \right)$. This shows that the pair\\
$(\np \left( R^2\mu_* (\nl^\vee \otimes pr_V^*\kv) \right),D^+)$ represents the functor $\Hifu^m_V$. Analogous arguments show that the pair $(\np \left( R^2\mu_* \nl \right),D^-)$ represents the functor $\Hifu^{k-m}_V$.

To prove the second claim, we observe that by construction of the divisor $D^+$ we have
\[
\co(D^+)|_{\np \left( R^2\mu_*( \nl^\vee \otimes pr_V^*\kv) \right) \times \{ p \}} \cong \co_{\np \left( R^2\mu_*( \nl^\vee \otimes pr_V^*\kv) \right)} (1) \otimes \rhotil_+^* (\nl|_{\Pic^m_V \times \{ p \}}),
\]
and analogously
\[
\co(D^-)|_{ \np (R^2 \mu_* \nl) \times \{  p \}} \cong \co|_{ \np (R^2 \mu_* \nl)}(1) \otimes \rhotil_-^* (pr_V^* \kv \otimes \nl^\vee)|_{ \Pic^m_V \times \{ p \}} .
\]
Here $\rhotil_-$ is the projective fibration
\[
\rhotil_-:\np (R^2 \mu_* \nl) \longrightarrow \Pic^m_V.
\]
\end{proof}
\begin{prop}\label{prop:zeroloc}
Let $V$ be a surface, and fix a class $m \in H^2(V,\mathbb{Z})$. Choose a Poincar\'e line bundle $\nl$ on $\Pic^m_V \times V$, and denote by $\mu$ the projection $\Pic^m_V \times V \to \Pic^m_V$. Suppose we have a global resolution
\[
\cm_1 \stackrel{\varphi}{\longrightarrow} \cm_2 \stackrel{\psi}{\longrightarrow} \cm_3
\]
of the complex $R^\bul \mu_* \nl$ by locally free sheaves. Denote by $\nu$ the projection $P(\cm_1) \to \Pic^m_V$, and let $\lambda$ be the section in $\co_{P(\cm_1)}(1)\otimes \nu^*( \ker \psi )$ induced by $\varphi$. Then there is a canonical isomorphism
\[
\Hilb^m_V \stackrel{\cong}{\longrightarrow} Z(\lambda).
\]
\end{prop}
\begin{proof}
By relative duality, the complex
\[
\cm^\vee_3 \stackrel{\psi^\vee}{\longrightarrow} \cm^\vee_2 \stackrel{\varphi^\vee}{\longrightarrow} \cm^\vee_1
\]
is a global resolution of $R^\bul \mu_* ( \nl^\vee \otimes pr_V^*\kv )$.
In particular, we have
\[
\coker \varphi^\vee \cong R^2 \mu_* ( \nl^\vee \otimes pr_V^*\kv ).
\]
On $P(\cm_1)= \np (\cm^\vee_1)$, we form the composition
\[
\chi : \nu^* \cm^\vee_2 \to \nu^* \cm^\vee_1 \to \co_{P(\cm_1)}(1).
\]
Lemma \ref{lem:finmod} implies, that there is a canonical isomorphism
\[
\Hilb^m_V \stackrel{\cong}{\longrightarrow} Z(\chi).
\]
The morphism $\chi$ factorizes through $\nu^* \cm^\vee_2 \to \coker \nu^*\psi^\vee$.
The dual sheaf of $\coker \nu^*\psi^\vee$ is $\ker \nu^* \psi$. Since $\nu: P(\cm_1) \to \Pic^m_V$ is smooth, we have $\ker \nu^*\psi = \nu^* \ker \psi$. This proves our claim.
\end{proof}
\begin{lem}\label{lem:subvec}
Let $V$ be a surface, and fix a class $m \in H^2(V,\mathbb{Z})$. Choose a Poincar\'e line bundle $\nl$ on $\Pic^m_V \times V$, and denote by $\mu$ the projection $\Pic^m_V \times V \to \Pic^m_V$. Let
\[
\cm_1 \stackrel{\varphi}{\longrightarrow} \cm_2 \stackrel{\psi}{\longrightarrow} \cm_3
\]
be a global resolution of the complex $R^\bul \mu_* \nl$. If
\begin{itemize}
	\item[i)] $\Hilb^{k-m}_V=\emptyset$ or
	\item[ii)] the fibered product $\Hilb^m_V \times_{\Pic^m_V} \Hilb^{k-m}_V$ is empty, $m(m-k)\geq 0$, and $\Hilb^m_V\neq \emptyset$,
\end{itemize}
then $\ker \psi \subset \cm_2$ is a subvectorbundle.
\end{lem}
\begin{proof}
If $\Hilb^{k-m}_V=\emptyset$, then by Lemma \ref{lem:finmod} the sheaf $R^2\mu_* \nl$ vanishes, and hence $\psi: \cm_2 \longrightarrow \cm_3$ is an epimorphism. In particular, the kernel of this morphism is a subvectorbundle.

Suppose now that both Hilbert schemes $\Hilb^m_V$ and $\Hilb^{k-m}_V$ are nonempty. Let $U_1$ be the complement of the Brill-Noether locus of the map $\Hilb^{k-m}_V \to Pic^m_V$, and let $U_2$ be the complement of the Brill-Noether locus of the map $\Hilb^m_V \to Pic^m_V$. Then $\psi|_{U_1}$ is an epimorphism, and hence $\ker \psi |_{U_1} \subset \cm_2 |_{U_1}$ is a subvectorbundle. Analogously, $\varphi^\vee|_{U_2}$ is an epimorphism, and hence $\im \varphi |_{U_2} \subset \cm_2|_{U_2}$ is a subvectorbundle.

We claim:
\[
\im \varphi |_{U_2}= \ker \psi |_{U_2}.
\]
Since the fibered product $\Hilb^m_V \times_{\Pic^m_V} \Hilb^{k-m}_V$ is empty, we have $U_1 \cup U_2 = \Pic^m_V$. Hence $\varphi$ is generically injective and $\psi$ is generically surjective. This yields
\[
\rk R^\bullet \mu_* \nl \leq 0.
\]
Conversely, since both Hilbert schemes are nonempty, the surface $V$ is neither rational nor ruled. This implies
\begin{eqnarray*}
\rk R^\bullet \mu_* \nl &=& \chi (\ov)+\frac{m(m-k)}{2}\\
&\geq & \chi (\ov)\\
&\geq &0.
\end{eqnarray*}
We now show that the induced morphism
\[
\bar{\psi} : (\cm_2/\im \varphi)|_{U_2} \to \cm_3|_{U_2}
\]
is a monomorphism, which implies our claim and ends the proof. We already know that $(\cm_2/\im \varphi)|_{U_2}$ is locally free. Moreover, since $\rk R^\bullet \mu_* \nl =0$, $\bar{\psi}|_{U_1 \cap U_2}$ is an isomorphism, which implies that $\bar{\psi}$ is generically injective and hence a monomorphism. This proves our claim, which in turn yields that $\ker (\psi) \subset \cm_2$ is a subvectorbundle.
\end{proof}
\begin{prop}\label{prop:newwall}
Let $V$ be a surface, and let $m \in H^2(V,\mathbb{Z})$. Fix a point $p \in V$ and a normalized Poincar\'e line bundle $\nl$ on $\Pic^m_V \times V$. Denote by $\rho_+$ the morphism $\Hilb^m_V \to \Pic^m_V$. If
\begin{itemize}
	\item[i)] $\Hilb^{k-m}_V=\emptyset$ or
	\item[ii)] the fibered product $\Hilb^m_V \times_{\Pic^m_V} \Hilb^{k-m}_V$ is empty and $m(m-k)\geq 0$,
\end{itemize}
then
\[
(\rho_+)_* \left( \sum_i (c_1( \co (\nd)|_{\Hilb^m_V\times \{ p \}}))^i\cap \{ \Hilb^m_V \} \right) =\shat ( \sigma_{\leq 1} R^\bullet \mu_* \nl).
\]
\end{prop}
\begin{proof}
Fix a smooth effective divisor $H \subset V$, such that $H^i(\cL (H))=0$ for all $[\cL]\in Pic^m_V$ and all $i>0$. Set $H_P:=\Pic^m_V \times H$, $\cm_1:= \mu_*\nl(H_P)$, and fix a global resolution $\cm_2 \stackrel{\psi}{\longrightarrow} \cm_3$ of the complex $R^\bul \mu_* \nl(H_P)|_{H_P}$. Then
\[
\cm_1 \stackrel{\varphi}{\longrightarrow} \cm_2 \stackrel{\psi}{\longrightarrow} \cm_3
\]
is a global resolution of the complex $R^\bul \mu_* \nl$. By Lemma \ref{lem:subvec}, the sheaf $\ker \psi$ is locally free, hence
\[
\cm_1 \to \ker  \psi
\]
is a global resolution of the complex $\sigma_{\leq 1} R^\bul \mu_* \nl$. By Prop.~\ref{prop:zeroloc}, $\Hilb^m_V$ is canonically isomorphic to the zero locus $Z(\lambda)$, where $\lambda$ is the section in $\co_{P(\cm_1)}(1) \otimes \nu^* \ker \psi$ induced by $\psi$. By construction of the virtual fundamental class $\{ \Hilb^m_V\}$, this cycle class is the localized Euler class $\lb Z(\lambda) \rb$. Moreover, by Lemma \ref{lem:finmod}, we have
\[
c_1(\co_{P(\cm_1)}(1))|_{Z(\lambda)}= c_1( \co(\nd )|_{\Hilb^m_V \times \{ p \}}).
\]
Therefore, our claim follows from Cor.~\ref{cor:segre}.
\end{proof}
\section{Poincar\'e invariants of projective surfaces}
In this section, a surface is again a smooth connected projective complex surface.

\subsection{Definition of the Poincar\'e invariant}
Let $V$ be a surface, $p\in V$ an arbitrary point. Fix a class $m \in H^2(V,\mathbb{Z})$, denote by $\nd^+$ the universal divisor over the Hilbert scheme $\Hilb^m_V$, and set
\[
u^+:= c_1\left( \OO (\nd^+)|_{ \Hilb^m_V \times \{ p \}}\right) \in H^2(\Hilb^m_V,\mathbb{Z}).
\]
Since $V$ is connected, the class $u^+$ does not depend on the chosen point $p$. Likewise, we denote by $\nd^-$ the universal divisor over the Hilbert scheme $\Hilb^{k-m}_V$, where $k$ is the first Chern class of the canonical line bundle $\kv$. We put
\[
u^-:= c_1\left( \OO (\nd^-)|_{ \Hilb^{k-m}_V \times \{ p \}}\right) \in H^2(\Hilb^{k-m}_V,\mathbb{Z}).
\]
Denote by $\rho^\pm$ the following morphisms:
\begin{eqnarray*}
\rho^+: \Hilb^m_V & \longrightarrow & \Pic^m_V\\
D & \longmapsto & [\ov (D)]
\end{eqnarray*}
\begin{eqnarray*}
\rho^-: \Hilb^{k-m}_V & \longrightarrow & \Pic^m_V\\
D' & \longmapsto & [\kv (-D')]
\end{eqnarray*}
Recall that $\lb \Hilb^m_V \rb$ denotes the virtual fundamental class of the Hilbert scheme $\Hilb^m_V$ defined in section 3.1; it is an element in the Chow group $A_*(\Hilb^m_V)$. By abuse of notation, we will denote its image in $H_*(\Hilb^m_V,\nz )$ by the same symbol.
\begin{dfn}
Let $V$ be a surface.  The \emph{Poincar\'e invariant} of $V$ is the map
\begin{eqnarray*}
(P^+_V,P^-_V) : H^2(V, \mathbb{Z}) & \longrightarrow & \Lambda^* H^1(V,\mathbb{Z}) \times
\Lambda^* H^1(V,\mathbb{Z})\\
m & \longmapsto & (P^+_V(m),P^-_V(m)),
\end{eqnarray*}
defined by
\[
P^+_V(m):= \rho^+_* \left( \sum_i (u^+)^{i} \cap \lb \Hilb^m_V\rb \right)
\]
and
\[
P^-_V(m):= (-1)^{ \chi (\ov ) + \frac{m(m-k)}{2}} \rho^-_* \left( \sum_i
(-u^-)^{i} \cap \lb \Hilb^{k-m}_V\rb \right),
\]
if $m \in NS(V)$, and by $P^\pm_V(m):=0$ otherwise.
\end{dfn}
\begin{rem}
The map $P^-_V :H^2(V, \mathbb{Z}) \longrightarrow \Lambda^* H^1(V,\mathbb{Z})$ is determined by the map $P^+_V : H^2(V, \mathbb{Z}) \longrightarrow \Lambda^* H^1(V,\mathbb{Z})$ in the following way: When we denote the component of degree $2i$ by $[P^\pm_V(m)]^{2i}$, then we have
\[
[P^-_V(m)]^{2i}= (-1)^{\chi (\ov)+i} [P^+_V(k-m)]^{2i}.
\]
\end{rem}
The following is a first nontrivial example, which will later play a role.
\begin{exa}
Let $\Gamma \subset \nc$ be a lattice, and let $E = \nc/\Gamma$ be the corresponding elliptic curve. We denote by $[z]\in E$ the equivalence
class of $z \in \nc$. Fix an integer $n>1$, a $n$-torsion point $[\zeta] \in E$, and set $\eps := \exp \frac{2\pi i}{n}$. Let the cyclic group
$\langle \eps \rangle$ act on $\np^1 \times E$ by
\[
\eps \cdot ([t_0,t_1],[z]) := ([t_0, \eps t_1],[z +\zeta]).
\]

The quotient $V:=(\np^1 \times E)/\langle \eps \rangle$ is a ruled surface over the elliptic curve $E/\langle [\zeta] \rangle$. We denote by
$[[t_0,t_1],[z]]$ the equivalence class of a point $([t_0,t_1],[z])
\in \np^1 \times E$ in $V$. The surface $V$ admits an elliptic
fibration $\varphi :V \to \np^1$ over the projective line, which sends
a point $[[t_0,t_1],[z]] \in V$ to $[t_0^n,t_1^n] \in \np^1$. This
fibration has exactly two singular fibers of type $nI_0$ over the
points $0$ and $\infty$, which we will denote by $nF_0$ and
$nF_{\infty}$ respectively. Let $F$ be a regular fiber of $\varphi$, and let $m\in H^2(V,\nz )$ be the Poincar\'e dual of $[F]$.

{\bf Claim:} One has
\[ 
\Hilb^m_V \cong |F| \cup \{ a F_0 + (n-a)F_{\infty} \mid a \in\{ 1, \ldots , n-1 \} \},
\]
and
\[
P^+_V(m)=n+1.
\]
\end{exa}
\begin{proof}
Since $V$ is ruled, its homology has no torsion and we infer $[F_0]=[F_\infty]$. Any effective divisor $D \in \Hilb^m_V$ is contained in the fibers of $\varphi: V \to \np^1$, since $D \cdot F=0$. This proves the first claim.

To prove the second claim, we have to compute the degree of the line bundle
\[
\left( R^1 \pi_* \co_\nd (\nd)\right)|_{|F|}
\]
on $|F| \cong \np^1$.
We find isomorphisms:
\begin{eqnarray*}
\left( R^0\pi_*\co \right)|_{|F|} & \cong & \co_{|F|}\otimes H^0(\co_F)\\
\left( R^0\pi_*\co(\nd)\right)|_{|F|}& \cong & \co_{|F|}(1)\otimes H^0(\co_F)\\
\left( R^1\pi_*\co\right)|_{|F|}& \cong & \co_{|F|}\otimes H^1(\co_F)\\
\left( R^1\pi_*\co(\nd)\right)|_{|F|}& \cong & \co_{|F|}(1)\otimes H^1(\co_F)\\
\end{eqnarray*}
Hence the long exact sequence
\[
\xymatrix{
0 \ar[r]& \left( R^0\pi_*\co \right)|_{|F|} \ar[r] & \left( R^0\pi_*\co(\nd)\right)|_{|F|} \ar[r] & \left( R^0 \pi_* \co_\nd (\nd)\right)|_{|F|}\ar[r] &\\
\ar[r] &\left( R^1\pi_*\co\right)|_{|F|} \ar[r] & \left( R^1\pi_*\co(\nd)\right)|_{|F|} \ar[r] & \left( R^1 \pi_* \co_\nd (\nd)\right)|_{|F|}\ar[r] & 0}
\]
yields
\begin{eqnarray*}
\deg \left( R^1 \pi_* \co_\nd (\nd)\right)|_{|F|} &=& \deg \left( R^0 \pi_* \co_\nd (\nd)\right)|_{|F|}.
\end{eqnarray*}
Since $\left( R^0 \pi_* \co_\nd (\nd)\right)|_{|F|}$ is isomorphic to the tangent bundle of $|F| \cong \np^1$, we obtain
\begin{eqnarray*}
P^+_V(m) &=& n-1 + \deg \left( R^1 \pi_* \co_\nd (\nd)\right)|_{|F|}\\
&=& n+1.
\end{eqnarray*}
\end{proof}
\subsection{Deformations}
In this subsection we study the behaviour of the Poin\-car\'e invariants under deformations. In order to make a precise statement, we need a slightly more general definition of the Poincar\'e invariants.

Fix a Poincar\'e line bundle $\nl$ on $\Pic^m_V \times V$, and let $\mu: \Pic^m_V \times V \to \Pic^m_V$ and $pr_V: \Pic^m_V \times V \to V$ be the projections. Recall that we have a projective fibration
\[
\rhotil: \np \left( R^2\mu_* (\nl^\vee \otimes pr_V^*\kv) \right) \to \Pic^m_V
\]
and a canonical section $\Phi$ in the line bundle
\[
(\rhotil \times id_V)^* (\nl) \otimes \pitil^*\co(1)
\]
on $\np \left( R^2\mu_* (\nl^\vee \otimes pr_V^*\kv )\right) \times V$, whose vanishing locus we denoted by $D^+$. In Lemma \ref{lem:finmod}, we have shown that the pair 
\[
\left( \np \left( R^2\mu_* (\nl^\vee \otimes pr_V^*\kv) \right), D^+\right)
\]
represents the functor $\Hifu^m_V$. Analogously, we obtained a pair
\[
\left(\np \left( R^2\mu_* \nl \right),D^- \right)
\]
representing $\Hifu^{k-m}_V$.
From this description, we get relatively ample line bundles on the Hilbert schemes $\Hilb^m_V$ and $\Hilb^{k-m}_V$, which we denote by $\co_\nl^+(1)$ and $\co_\nl^-(1)$. We set
\[
u_\nl^\pm:= c_1(\co_\nl^\pm(1)).
\]
These classes depend on the choice of a Poincar\'e line bundle, but the formal cohomology rings
\[
\Lambda^* H^1(V,\nz)^\vee [u_\nl^\pm ]
\]
are independent up to a canonical isomorphism. To be more precise, if $\nl'$ is a second Poincar\'e line bundle, then there is a line bundle $\cT$ on $\Pic^m_V$ and an isomorphism
\[
\nl' \stackrel{\cong}{\longrightarrow} \nl \otimes \mu^*\cT.
\]
This yields isomorphism
\[
\co_{\nl'}^+(1) \stackrel{\cong}{\longrightarrow} \co_\nl^+(1) \otimes \rho_+^*\cT^\vee
\]
and
\[
\co_{\nl'}^-(1) \stackrel{\cong}{\longrightarrow} \co_\nl^-(1) \otimes \rho_-^*\cT,
\]
and we obtain
\[
 u_{\nl'}^\pm = u_\nl^\pm \mp \rho_\pm^* c_1(\cT).
\]
Therefore, sending $u_\nl^\pm$ to $u_{\nl'}^\pm \pm c_1(\cT)$ gives rise to canonical isomorphisms
\[
\Lambda^* H^1(V,\nz)^\vee [u_\nl^\pm ] \stackrel{\cong}{\longrightarrow} \Lambda^* H^1(V,\nz)^\vee [ u_{\nl'}^\pm ] .
\]
By evaluating cohomology classes on the cycles $\lb \Hilb^m_V \rb$ and\\
$(-1)^{\chi (\ov )+\frac{m(m-k)}{2}} \lb \Hilb^{k-m}_V \rb$ we obtain maps
\[
\mathcal{P}^\pm_{V,\nl}(m): \Lambda^* H^1(V,\nz)^\vee [ u_\nl^\pm ] \to \nz.
\]
\begin{rem}
If $\nl_p$ is a normalized Poincar\'e line bundle, i.e.~if \\
$\nl_p |_{Pic^m_V \times \{ p \}} \cong \co_{\Pic^m_V}$ for some point $p \in V$, then there are isomorphisms $\co( \nd^+)|_{\Hilb^m_V \times \{ p \} } \cong \co^+_{\nl_p}(1)$ and $\co( \nd^-)|_{\Hilb^{k-m}_V \times \{ p \} } \cong \co^-_{\nl_p}(1)$, hence $u_{\nl_p}^\pm = u^\pm$. In this case we have
\[
\mathcal{P}^\pm_{V,\nl_p}(m)(\alpha (\pm u_{\nl_p}^\pm)^{\frac{m(m-k)-\deg \alpha}{2}})= \langle P^\pm_V(m), \alpha \rangle
\]
for any homogeneous element $\alpha \in \Lambda^*H^1(V,\nz)^\vee$.
\end{rem}
Let $v: V \to S$ be a family of surfaces over an irreducible variety $S$. Recall that a class $\alpha \in A_*(\Pic^{\underline m}_{V/S})$ determines a family of classes $\alpha_s \in A_*(\Pic^{\underline{m}(s)}_{V_s})$. We denote the Poincar\'e dual of the homology class associated to $\alpha_s$ by the same symbol.
\begin{prop}
Let $v: V \to S$ be a smooth, connected family of surfaces. Fix a class $\underline{m} \in H^0(S,R^2v_* \underline{\nz})$ and suppose there exists a Poincar\'e line bundle $\nl$ on $\Pic^{\underline m}_{V/S} \times_S V$. For a point $s \in S$, we denote by $\nl_s$ the induced Poincar\'e line bundle on $\Pic^{\underline{m}(s)}_{V_s} \times V_s$. For every element $\alpha \in A_*(\Pic^{\underline m}_{V/S})$ and every $i\in \nn$, the pair
\[
\left( \mathcal{P}^+_{V_s, \nl_s}(\underline{m}(s))( \alpha_s \cdot (u^+_{\nl_s})^{i}),\mathcal{P}^-_{V_s, \nl_s}(\underline{m}(s))( \alpha_s \cdot (-u^-_{\nl_s})^{i}) \right)
\]
is independent of the point $s \in S$.
\end{prop}
\begin{proof}
Denote by $\mathfrak{u}^+$ the first Chern class of the line bundle $\co_\nl (1)$ on $\Hilb^{\underline m}_{V/S}$, and fix a point $s \in S$. By \cite[Prop.10.1]{Fu} and Cor.~\ref{cor:basechange} we have
\begin{eqnarray*}
(u^+_{\nl_s})^{i} \cap \lb \Hilb^{\underline{m}(s)}_{V_s}\rb &=&(u^+_{\nl_s})^{i} \cap \lb \Hilb^{\underline{m}}_{V/S} \rb_s\\
&=&\left( (\mathfrak{u}^+)^i \cap \lb \Hilb^{\underline{m}}_{V/S} \rb \right)_s.
\end{eqnarray*}
Another application of \cite[Prop.10.1]{Fu} yields
\[
(\rho_s^+)_* \left( (u^+_{\nl_s})^{i} \cap \lb \Hilb^{\underline{m}(s)}_{V_s}\rb \right) = \left( \rho^+_* \left( (\mathfrak{u}^+)^i \cap \lb \Hilb^{\underline{m}}_{V/S} \rb \right) \right)_s .
\]
Now \cite[Cor.10.1]{Fu} implies
\[
(\rho_s^+)_* \left( (u^+_{\nl_s})^{i} \cap \lb \Hilb^{\underline{m}(s)}_{V_s}\rb \right) \cdot \alpha_s = \left( \rho^+_* \left( (\mathfrak{u}^+)^i \cap \lb \Hilb^{\underline{m}}_{V/S} \rb \right)  \cdot \alpha\right)_s.
\]
Hence our claim follows by conservation of numbers \cite[Prop.10.2]{Fu}.
\end{proof}
\begin{rem}
Each of the  following conditions is sufficient for the existence of a Poincar\'e line bundle:
\begin{itemize}
	\item the family $v: V \to S$ admits a section;
	\item the base scheme $S$ is a curve.
\end{itemize}
\end{rem}
For the first condition, see \cite[Prop.2.1]{gr}. The second follows from the lower term sequence of the Leray spectral sequence and the vanishing of $H^2(S,v_* \ov^*)$.
\subsection{A blow-up formula}
Let $\sigma: \hat{V} \to V$ be the blow-up of a point $p \in V$, let $E$ be the exceptional curve, and denote by $e \in H^2(\hat{V}, \nz )$ the Poincar\'e dual of the class $[E]$. We want to compare the Poincar\'e invariants of $\hat{V}$ and $V$.

Recall that the push down of an effective divisor $\hat{D} = \hat{D}_0+lE$ on $\hat{V}$ with $E \not \subset \hat{D}_0$ is the unique divisor $\sigma_! \hat{D}$ on $V$, whose strict transform is $\hat{D}_0$; its total transform is $\sigma^* \sigma_! \hat{D} =\hat{D}+ (\hat{D} \cdot E)E$.
Now fix a class $m \in H^2(V,\nz)$, and set $\hat{m}=\sigma^*m$. By pushing down divisors from $\hat V$ to $V$ we obtain maps
\begin{eqnarray*}
\nu_l : \Hilb^{\hat{m}+l\cdot e}_{\hat V} & \longrightarrow & \Hilb^m_V\\
\hat{D} & \longmapsto & \sigma_! \hat{D}
\end{eqnarray*}
for all integers $l$.

We start by observing that for $l \geq 0$, the map $\nu_l$ is an isomorphism: its inverse sends a divisor $D$ in $V$ to $\sigma^*D + lE$.
\begin{prop}\label{prop:blowup1}
Let $V$ be a surface, fix a point $p\in V$, and denote by $\sigma : \hat{V} \to V$ the blow-up of $V$ in $p$. For every class $m \in H^2(V,\nz)$, the isomorphism
\[
\nu_0 : \Hilb^{\hat m}_{\hat V} \longrightarrow \Hilb^m_V
\]
identifies the virtual fundamental classes:
\[
(\nu_0)_* \lb \Hilb^{\hat m}_{\hat V} \rb = \lb \Hilb^m_V \rb.
\]
\end{prop}
\begin{proof}
Denote by $\nd$ the universal divisor on $\Hilb^m_V \times V$ and by $\ndhat$ the universal divisor on $\Hilb^{\hat m}_{\hat V} \times \hat{V}$. Let $\pi : \Hilb^m_V \times V \to \Hilb^m_V$ and $\pihat : \Hilb^{\hat m}_{\hat V} \times \hat{V} \to \Hilb^{\hat m}_{\hat V}$ be the projections. Pulling back the short exact sequence
\[
0 \to \co \to \co(\nd) \to \co_\nd (\nd ) \to 0
\]
from $\Hilb^m_V \times V$ to $\Hilb^{\hat m}_{\hat V} \times \hat{V}$ yields 
\[
\co \to \co(\ndhat ) \to \co_\ndhat (\ndhat ) \to 0.
\]
Therefore, we obtain an isomorphism
\[
\co_\ndhat (\ndhat) \cong (\nu_0\times \sigma)^* \co_\nd(\nd).
\]
Since $\sigma_* \co_{\hat V} \cong \ov$, we have $(\nu_0\times \sigma)_* \co_{\\Hilb^{\hat m}_{\hat V} \times {\hat V}} \cong \co_{\Hilb^m_V \times V}$. Applying the push-pull formula, we find an isomorphism
\[
\co_\nd (\nd) \stackrel{\cong}{\longrightarrow} R^\bul (\nu_0\times \sigma)_* \co_\ndhat (\ndhat ).
\]
Applying $R^\bul \pi_*$ yields an isomorphism
\[
R^\bul\pi_* \co_\nd (\nd ) \stackrel{\cong}{\longrightarrow} R^\bul (\nu_0)_* R^\bul \pihat_* \co_\ndhat (\ndhat ).
\]
This proves our claim.
\end{proof}
Let $\hat{D}\in \Hilb^{\hat{m}-e}_{\hat V}$ be a divisor on $\hat V$. The equality $\hat{D}\cdot E=1$ implies that the point $p$ lies on $\sigma_! \hat{D}$. Conversely, if $D \in \Hilb^m_V$ passes through $p$, then the total transform $\sigma^*D$ can be written as $\sigma^*D= \hat{D}+E$ with $\hat{D} \in \Hilb^{\hat{m}-e}_{\hat{V}}$. Therefore the map
\[
\nu_{-1}: \Hilb^{\hat{m}-e}_{\hat V} \to \Hilb^m_V
\]
is a closed embedding. Its image consists of all divisors $D\in \Hilb^m_V$ which pass through $p$. In particular, it is the zero locus of a section in a line bundle.

In a next step, we want to generalize this observation to arbitrary negative integers $l$. In order to simplify the notation, we set $\nutil_l:= \nu_0^{-1} \circ \nu_l$ for $l \in \nz$. Note that $\nutil_l: \Hilb^{\hat{m}+l\cdot e}_{\hat V} \to \Hilb^{\hat m}_{\hat V}$ sends a divisor $\hat{D}$ on $\hat{V}$ to $\hat{D}-lE$.
\begin{lem}\label{lem:notflat}
Let $V$ be a surface, and fix a point $p \in V$. Denote by $\sigma: \hat{V} \to V$ the blow-up of $p \in V$, and by $E$ the exceptional curve. Then
\[
R^\bullet \sigma_* (\co_E(lE)) \cong \left\{
\begin{array}{cl}
0 & \text{if $l=1$,}\\
\ov/\mathcal{J}_p \otimes H^0(\co_E(lE)) & \text{if $l \leq 0$,}\\
\ov/\mathcal{J}_p \otimes H^1(\co_E(lE))[-1] & \text{if $l \geq 2$.}
\end{array} \right.
\]
\end{lem}
\begin{proof}
Consider the following comutative diagram:
\[
\xymatrix{
E \ar[r]^{j_E} \ar[d]_{\sigma_E} & \hat{V}\ar[d]^\sigma\\
\{ p \} \ar[r]_{j_p} & V}
\]
Since $j_E$ is a closed embedding, the functor $(j_E)_*$ is right exact, and hence
\[
R^\bullet \sigma_* \co_E(lE) \cong R^\bullet (j_p)_* R^\bullet (\sigma_E)_* \co_E(lE).
\]
This proves our claim.
\end{proof}
\begin{lem}
Let $\sigma: \hat{V} \to V$ be the blow-up of a point $p\in V$, and let $E$ be the exceptional curve. Fix a class $m \in H^2(V,\nz)$, denote by $\ndhat$ the universal divisor in $\Hilb^{\hat m}_{\hat V} \times \hat{V}$, and set $\ne := \Hilb^{\hat m}_{\hat V} \times E$. Denote by $\pihat$ the projection $\Hilb^{\hat m}_{\hat V} \times \hat{V} \to \Hilb^{\hat m}_{\hat V}$. For every negative integer $l$, the sheaf $\pihat_* \co_{-l \ne} (\ndhat )$ is locally free and has the base change property. Moreover, if $\zeta_l$ denotes the canonical section in $\pihat_* \co_{-l \ne} (\ndhat )$, then $\nutil_l$ induces an isomorphism
\[
\Hilb^{\hat{m}+l\cdot e}_{\hat V} \stackrel{\cong}{\longrightarrow} Z(\zeta_l) \subset \Hilb^{\hat m}_{\hat V}.
\]
\end{lem}
\begin{proof}
To prove the first claim, we show that $H^1( \co_{-lE}(\hat{D}))=0$ for every divisor $\hat{D} \in \Hilb^{\hat m}_{\hat V}$. We proceed by induction on $-l$.

For $l=0$ there is nothing to show. For the induction step, consider the following short exact sequence:
\[
0 \to \co_E({\hat D}+lE) \to \co_{-(l-1)E}(\hat{D}) \to \co_{-lE}(\hat{D}) \to 0
\]
The push-pull formula yields an isomorphism
\[
R^\bullet \sigma_* \co_E(\hat{D}+lE) \cong \left( R^\bullet \sigma_* \co_E(-lE) \right) \otimes \co_V(\sigma_!\hat{D}),
\]
hence Lemma \ref{lem:notflat} implies $H^1(\co_E(\hat{D}+lE))=0$.
Since by assumption $H^1(\co_{-lE}(\hat{D}))$ vanishes, the long exact cohomology sequence yields\\
 $H^1(\co_{-(l-1)E}(\hat{D}))=0$.

Our second claim follows from the fact that a divisor $\hat{D} \in \Hilb^{\hat m}_{\hat V}$ can be written as $\hat{D}=\hat{D}'+(-l)E$ with $\hat{D}'\geq 0$ iff the composition
\[
\co \to \co(\hat{D}) \to \co_{-lE}(\hat{D})
\]
vanishes.
\end{proof}
\begin{prop}\label{prop:blowup2}
Let $\sigma: \hat{V} \to V$ be the blow-up of a point $p\in V$. For every class $m \in H^2(V,\nz)$ and every negative integer $l$ we have
\[
\lb \Hilb^{\hat{m} + l \cdot e}_{\hat V} \rb = 0^!_{\hat{\pi}_*\co_{-l\ne} ( \ndhat )} \lb \Hilb^{\hat m}_{\hat V} \rb
\]
and
\[
(\nu_l)_* \lb \Hilb^{\hat{m} + l \cdot e}_{\hat V} \rb = c_1 \left(\co (\nd )|_{\Hilb^m_V \times \{ p \} } \right)^{\left( \dfrac{l}{2} \right)} \cap \lb \Hilb^m_V\rb.
\]
\end{prop}
\begin{proof}
Denote by $\ndhat_l$ the universal divisor on $\Hilb^{\hat{m}+l\cdot e}_{\hat V} \times \hat{V}$, and set $\ne: = \Hilb^{\hat{m}+l\cdot e}_{\hat V} \times E$. The short exact sequence
\[
0 \to \co_{\ndhat_l} (\ndhat_l) \to \co_{\ndhat_l -l \ne }(\ndhat_l -l \ne ) \to \co_{-l\ne} (\nd_l -l\ne ) \to 0
\]
gives rise to the following distinguished triangle on $\Hilb^{\hat{m}+l\cdot e}_{\hat V}$:
\[
\xymatrix{
R^\bul \pihat_* \co_{\ndhat_l} (\ndhat_l) \ar[r] & R^\bul \pihat_*\co_{\ndhat_l -l \ne }(\ndhat_l -l \ne )\ar[d] \\
& R^\bul \pihat_*\co_{-l\ne} (\nd_l - l \ne )\ar[lu]^{[1]}}
\]
This is the necessary compatibility datum for the obstruction theories of the Hilbert schemes $\Hilb^{\hat{m}+l\cdot e}_{\hat V}$ and $\Hilb^{\hat{m}}_{\hat V}$  \cite[Thm.1]{kkp}, and hence proves the first claim. To show the second claim, we have to compute the top Chern class of the vector bundle $(\nu_0)_* \pihat_* \co_{-l \ne} (\ndhat_0)$.

\noindent
{\bf Claim:} For all $l\leq 0$, we have
\[
c ((\nu_0)_* \pihat_* \co_{-l \ne} (\ndhat_0))=\left( 1+c_1(\co(\nd)|_{\Hilb^m_V\times \{ p \}})\right)^{\left( \dfrac{l}{2} \right)}.
\]
We proceed by induction on $-l$. For $l=0$, there is nothing to show. For the induction step, consider the following short exact sequence:
\[
0 \to \co_\ne (\ndhat_0+l\ne ) \to \co_{(-l+1)\ne}(\ndhat_0) \to  \co_{-l\ne}(\ndhat_0) \to 0.
\]
This shows that the vector bundle $(\nu_0)_* \pihat_* \co_{(-l+1)\ne}(\ndhat_0)$ is an extension of $(\nu_0)_* \pihat_*\co_{-l\ne}(\ndhat_0)$ by $(\nu_0)_* \pihat_*\co_\ne (\ndhat_0+l\ne )$. By Lemma \ref{lem:notflat}, we have an isomorphism
\[
(\nu_0)_* \pihat_*\co_\ne (\ndhat_0+l\ne ) \stackrel{\cong}{\longrightarrow} \co (\nd)_{\Hilb^m_V \times \{ p\}} \otimes H^0(\co_E(lE)).
\]
Since $H^0(\co_E(lE))$ is a vector space of dimension $-l+1$, we have
\begin{eqnarray*}
c\left( (\nu_0)_* \pihat_* \co_{(-l+1)\ne}(\ndhat_0) \right) &=& \left( 1+c_1(\co(\nd)|_{\Hilb^m_V\times \{ p \}})\right)^{\left( \dfrac{l}{2} \right) + (-l+1)}\\
&=& \left( 1+c_1(\co(\nd)|_{\Hilb^m_V\times \{ p \}})\right)^{\left( \dfrac{l-1}{2} \right)}
\end{eqnarray*}
\end{proof}
\begin{prop}\label{prop:blowup3}
Let $\sigma: \hat{V} \to V$ be the blow-up of a point $p\in V$, and let $E$ be the exceptional curve. Fix a class $m \in H^2(V,\nz)$, denote by $\ndhat$ the universal divisor in $\Hilb^{\hat m}_{\hat V} \times \hat{V}$, and set $\ne := \Hilb^{\hat m}_{\hat V} \times E$. Let $\pihat$ be the projection $\Hilb^{\hat m}_{\hat V} \times \hat{V} \to \Hilb^{\hat m}_{\hat V}$. For every positive integer $l$, we have
\[
(\nu_l)_* \lb \Hilb^{\hat{m}+l\cdot e}_{\hat V} \rb = c_1 \left(\co (\nd )|_{\Hilb^m_V \times \{ p \} } \right)^{\left( \dfrac{l}{2} \right)} \cap \lb \Hilb^m_V\rb.
\]
\end{prop}
\begin{proof}
Let $\ndhat_0$ be the universal divisor on $\Hilb^{\hat m}_{\hat V} \times \hat{V}$, and set $\ne:= \Hilb^{\hat m}_{\hat V}\times E$. The short exact sequence
\[
0 \to \co_{\ndhat_0}(\ndhat_0) \to \co_{\ndhat_0+l\ne } (\ndhat_0+l\ne ) \to \co_{l\ne } (\ndhat_0+l\ne ) \to 0
\]
gives rise to the following distinguished triangle on $\Hilb^{\hat m}_{\hat V}$:
\[
\xymatrix{R^\bul \pihat_* \co_{\ndhat_0}(\ndhat_0) \ar[r] & R^\bul \pihat_*\co_{\ndhat_0+l\ne } (\ndhat_0+l\ne ) \ar[d]\\
&R^\bul \pihat_* \co_{l\ne } (\ndhat_0+l\ne ) \ar[ul]^{[1]}}
\]
Hence our claim follows by excess intersection Prop.~\ref{prop:ex}, once we know that for each $l \geq 0$ we have
\[
H^0(\co_{lE } (\hat{D}+lE ))=0\ \  \forall \hat{D} \in \Hilb^{\hat m}_{\hat V},
\]
and
\[
c \left( (\nu_0)_* R^1 \pihat_*\co_{l\ne } (\ndhat_0+l\ne ) \right) = \left( 1+c_1(\co(\nd)|_{\Hilb^m_V\times \{ p \}})\right)^{\left( \dfrac{l}{2} \right)}.
\]
These claims can be proved by induction on $l$. Since the arguments are very similar to those used in the case $l<0$, we omit the details.
\end{proof}
For a similar computation, see \cite[Prop.~43]{br}.

For an integer $n$ we define a truncation map
\[
\tau_{\leq n} : \Lambda^* H^1(V,\nz ) \longrightarrow \Lambda^* H^1(V,\nz )
\]
as follows: when $P=\sum_i P_i$ is the decomposition of a form $P$ into its homogeneous components $P_i \in \Lambda^i H^1(V,\nz )$, then
\[
\tau_{\leq n}(P):= \sum_{i=0}^n P_i.
\]
\begin{thm}\label{thm:blowup}
Let $\sigma: \hat{V} \to V$ be the blow-up of a point $p \in V$. Using the natural identification $\sigma^*:H^1(V,\nz ) \stackrel{\cong}{\longrightarrow} H^1(\hat{V},\nz )$, we have
\[
P^\pm_{\hat V}(\hat{m}+l\cdot e)= \tau_{\leq m(m-k) - 2\left( \dfrac{l}{2}\right) } P^\pm_V(m)
\]
for every class $m \in H^2(V,\nz)$ and for every integer $l$.
\end{thm}
\begin{proof}
This is an immediate consequence of Prop.~\ref{prop:blowup1}, Prop.~\ref{prop:blowup2} and Prop.~\ref{prop:blowup3}
\end{proof}
\subsection{A wall crossing formula}
Let $V$ be a surface. Recall that an element $c \in H^2(V,\mathbb{Z})$ is called {\em characteristic} iff $c \equiv k \mod{2}$. For a characteristic element $c \in H^2(V,\mathbb{Z})$, we denote by $\theta_c \in \Lambda^2 H^1(V,\mathbb{Z})^\vee$ the mapping
\begin{eqnarray*}
\theta_c: \Lambda^2 H^1(V,\mathbb{Z}) & \longrightarrow & \mathbb{Z}\\
a \wedge b & \longmapsto & \frac{1}{2} \langle a \cup b \cup c, [V]\rangle .
\end{eqnarray*}
\begin{lem}
Let $V$ be a surface with $p_g(V)=0$ and irregularity $q$. Fix a cohomology class $m\in H^2(V,\mathbb{Z})$, choose a normalized Poincar\'e line bundle $\nl$ on $\Pic^m_V\times V$, and denote by $\mu$ the projection $\Pic^m_V \times V \to \Pic^m_V$. Then
\begin{eqnarray*}
ch( \mu_! \nl) &=& \chi (\ov) + \frac{m(m-k)}{2}-\theta_{2m-k},\\
c( \mu_! \nl) &=& \exp (-\theta_{2m-k}).
\end{eqnarray*}
\end{lem}
\begin{proof}
By the Grothendieck-Riemann-Roch theorem \cite[Thm.15.2]{Fu} we have
\[
td(\Pic^m_V) \cdot ch (R\mu_! \nl) = \mu_! \left\{ td (\Pic^m_V \times V) \cdot ch (\nl) \right\}.
\]
Hence we need to compute those components of the expression
\[
\left\{ td (\Pic^m_V \times V) \cdot ch (\nl) \right\}
\]
which have bidegree $(*,4)$ with respect to the decomposition
\begin{eqnarray*}
H^*(\Pic^m_V \times V,\mathbb{Z}) & \cong & H^*(\Pic^m_V,\mathbb{Z}) \otimes H^*(V,\mathbb{Z})\\
& \cong & \Lambda^* H^1(V,\mathbb{Z})^\vee \otimes H^*(V,\mathbb{Z}).
\end{eqnarray*}
Set $f:=c_1(\nl)$. Then
\begin{eqnarray*}
f^{2,0} &=& 0 \in H^2(\Pic^m_V,\mathbb{Z}),\\
f^{1,1} &=& id \in \Hom (H^1(V,\mathbb{Z}),H^1(V,\mathbb{Z})),\\
f^{0,2} &=& m \in H^2(V,\mathbb{Z}),
\end{eqnarray*}
where the first equality holds since $\nl$ is normalized.

Next we compute $g:=f^2$. We obtain
\begin{eqnarray*}
g^{2,2} &=& -2\cdot (a \wedge b \mapsto a\cup b ) \in \Hom (\Lambda^2 H^1(V,\mathbb{Z}), H^2(V,\mathbb{Z})),\\
g^{1,3} &=& 2\cdot (a \mapsto a \cup m) \in \Hom (H^1(V,\mathbb{Z}),H^3(V,\mathbb{Z})),\\
g^{0,4} &=& m\cup m \in H^4(V,\mathbb{Z}),
\end{eqnarray*}
all other components being zero. Here the first equality needs justification. Choose a basis $v_1,\ldots ,v_{2q}$ of $H^1(V,\mathbb{Z})$, and denote by $w_1,\ldots , w_{2q}$ the dual basis of $H^1(V,\mathbb{Z})^\vee$. Then
\[
f^{1,1}= \sum_i w_i\otimes v_i,
\]
hence
\begin{eqnarray*}
g^{2,2} &=& \left( f^{1,1} \right)^2\\
&=& (\sum_i w_i\otimes v_i)\cup (\sum_i w_i\otimes v_i)\\
&=& - \sum_i\sum_j (w_i\wedge w_j) \otimes (v_i \cup v_j)\\
&=& -2 \sum_{i<j} (w_i\wedge w_j) \otimes (v_i \cup v_j).
\end{eqnarray*}
Now we compute the component of $f^3$ of bidegree $(2,4)$, the only component that does not vanish. We obtain
\begin{eqnarray*}
f^3 &=& 3 (f^{1,1})^2 \cup f^{0,2}\\
&=& -6 \cdot ( a \wedge b \mapsto a \cup b \cup m) \in \Hom ( \Lambda^2 H^1(V,\mathbb{Z}),H^4(V,\mathbb{Z})).
\end{eqnarray*}
Since $p_g(V)=0$, we have $\langle a \cup b \cup c \cup d, [V]\rangle =0$ for all $a,b,c,d \in H^1(V,\mathbb{Z})$ \cite{ll}. This implies that the $(4,4)$ part and hence $f^4$ itself vanishes.

Moreover, since $td(\Pic^m_V)=1$, we have
\begin{eqnarray*}
td(\Pic^m_V \times V) &=& pr_V^* td(V)\\
&=& pr_V^*(1 - \frac{1}{2} k + \chi (\ov ) \cdot PD [pt] ),
\end{eqnarray*}
where $pr_V: \Pic^m_V \times V \to V$ denotes the projection onto $V$.

Thus we obtain
\begin{eqnarray*}
ch (\mu_! \nl ) &=& \left\{ \exp f  \cup pr_V^* \left( 1 -\frac{k}{2} + \chi (\ov) \cdot PD[pt] \right) \right\} /[V]\\
&=& \left\{ (\exp f)^{*,4} - (\exp f)^{*,2}\cup pr_V^* \frac{k}{2}+\chi (\ov) \cdot PD[pt] \right\} /[V]\\
&=& \chi (\ov ) + \frac{m \cdot (m-k)}{2}- \theta_{2m-k}.
\end{eqnarray*}
The formula for the Chern class follows immediately since $H^*(\Pic^m_V,\mathbb{Z})$ has no torsion.
\end{proof}
\begin{lem}
Let $V$ be a surface of negative Kodaira dimension. Then there exists a smooth rational curve on $V$ with nonnegative selfintersection.
\end{lem}
\begin{proof}
By the Enriques classification \cite[p.188]{bpv}, the surface $V$ is either the projective plane or a blow-up of a geometrically ruled surface. In the first case any line will do, whereas in the second case we may take a general fibre of the composition
\[
V \to V_{min} \to C,
\]
where $V \to V_{min}$ is a minimal model, and $V_{min} \to C$ is a ruling.
\end{proof}
\begin{cor}\label{cor:34}
Let $V$ be a surface with $p_g(V)=0$, and suppose $m \in H^2(V,\mathbb{Z})$ satisfies $m(m-k) \geq 0$. Then one of the Hilbert schemes $\Hilb^m_V$, $\Hilb^{k-m}_V$ is empty, or we have $kod(V)\geq 0$, $q(V)=1$ and $m(m-k)=0$.
\end{cor}
\begin{proof}
Assume first that the Kodaira dimension of $V$ is negative, and fix a smooth rational curve $C$ on $V$ with $C^2\geq 0$. The adjunction formula yields
\[
\langle k,[C]\rangle \leq -2,
\]
which implies that $\langle m,[C]\rangle$ or $\langle k-m,[C]\rangle$ is negative. Hence $\Hilb^m_V$ or $\Hilb^{k-m}_V$ is empty.

Suppose now that both Hilbert schemes are nonempty. Then $kod(V) \geq 0$, hence $\chi(\ov)\geq 0$, which implies $q(V)=1$. Let $\sigma : V \to V_{min}$ be the minimal model of $V$, and fix elements $D \in \Hilb^m_V$ and $D'\in \Hilb^{k-m}_V$. We have
\begin{eqnarray*}
D \cdot D' &\geq &\sigma_! D \cdot \sigma_! D'\\
&\geq & 0,
\end{eqnarray*}
where the second inequality is a consequence of the fact that the canonical class of $V_{min}$ is numerically effective. This proves our last claim.
\end{proof}
\begin{thm}[Wall crossing formula]\label{thm:wcfpoinc}
Let $V$ be a surface with $p_g(V)=0$ and irregularity $q$. Fix a cohomology class $m\in H^2(V,\mathbb{Z})$, choose a normalized Poincar\'e line bundle $\nl$ on $\Pic^m_V\times V$, and denote by $\mu$ the projection $\Pic^m_V \times V \to \Pic^m_V$. Then
\begin{eqnarray*}
P^+_V(m) -P^-_V(m) & = & \shat (R \mu_! \nl)\\
&=& \sum_{j=0}^{\min \{ q ,\frac{m(m-k)}{2}\} } \frac{\theta_{2m-k}^{q-j}}{(q-j)!} \cap [\Pic^m_V].
\end{eqnarray*}
\end{thm}
\begin{proof}
Without loss of generality we may assume $m(m-k)\geq 0$. We distinguish the following three cases:
\begin{itemize}
	\item[-] $\Hilb^{k-m}_V =\emptyset$;
	\item[-] $\Hilb^m_V = \emptyset$;
	\item[-] both Hilbert schemes are nonempty.
\end{itemize}
In the first case we have
\begin{eqnarray*}
P^+_V(m)-P^-_V(m) &=& P^+_V(m) \\
&=& \shat (\sigma_{\geq 1} R^\bullet\mu_* \nl)\\
&=& \shat (R\mu_! \nl),
\end{eqnarray*}
where the second equality holds by virtue of Prop.~\ref{prop:newwall}.

Assume now $\Hilb^m_V = \emptyset$. By relative duality there is an isomorphism
\[
R^\bullet \mu_* (pr_V^* \kv \nl^\vee ) \stackrel{\cong}{\longrightarrow} \left( R^\bullet \mu_* \nl\right) ^\vee [-2],
\]
where $pr_V$ denotes the projection $\Hilb^{k-m}_V \times V \to V$. Therefore  we obtain analogously
\[
P^+_V(k-m) = \shat ((R\mu_! \nl)^\vee )
\]
and hence
\[
-P^-_V(m) = \shat (R\mu_! \nl).
\]

Suppose finally that both Hilbert schemes are nonempty. By Cor.~\ref{cor:34}, the surface $V$ is neither rational nor ruled, hence $\chi (\ov ) \geq 0$, and we have $q(V)=1$, $m(m-k)=0$.

Fix a global resolution
\[
\cm_1 \stackrel{\varphi}{\longrightarrow} \cm_2 \stackrel{\psi}{\longrightarrow} \cm_3
\]
of the complex $R^\bullet \mu_* \nl$. By Lemma \ref{lem:subvec}, the sheaf $\ker \psi$ is a subvectorbundle of $\cm_2$ and we obtain
\[
P^+_V(m) = c_1(\cm_1-\ker \psi ) \cap [Pic^m_V].
\]
By relative duality
\[
\cm_3^\vee \stackrel{\psi^\vee}{\longrightarrow} \cm_2^\vee \stackrel{\varphi^\vee}{\longrightarrow} \cm_1^\vee
\]
is a global resolution of the complex $R^\bullet \mu_* (pr_V^* \kv \otimes \nl^\vee )$.

We claim:
\[
\left( \cm_2 /\ker \psi \right)^\vee = \ker \varphi ^\vee.
\]
In order to see this, consider the following commutative diagram with exact rows and columns:
\[
\xymatrix{
& & 0 \ar[d] & 0 \ar[d] & &\\
& 0 \ar[r] & \im \varphi \ar[r] \ar[d] & \ker \psi \ar[d] \ar[r] & \ker \psi / \im \varphi \ar[r] & 0\\
& & \cm_2 \ar@{=}[r] \ar[d] & \cm_2 \ar[d]& & \\
0\ar[r] & \ker \psi / \im \varphi \ar[r]  & \coker \varphi \ar[r] \ar[d] & \cm_2/\ker \psi \ar[r] \ar[d] & 0 &\\
&& 0 & 0 &&}
\]
Dualizing the last row yields
\[
0 \longrightarrow \left( \cm_2/ \ker \psi \right)^\vee \longrightarrow \ker \varphi^\vee \longrightarrow \left( \ker \psi / \im \varphi \right)^\vee.
\]
Since $\ker \psi / \im \varphi$ is a skyscraper sheaf, our claim follows.

So Prop.~\ref{prop:newwall} yields
\begin{equation}\label{eq:pdual}
P^+_V(k-m)=c_1(\cm_3^\vee - (\cm_2/\ker \psi )^\vee )\cap [Pic^m_V].
\end{equation}
Since $\cm_2/\ker \psi$ is a vector bundle, equation \ref{eq:pdual} implies
\[
-P^-_V(m) = c_1(\cm_3- \cm_2/\ker \psi )\cap [Pic^m_V].
\]
Hence also in the third case we obtain
\begin{eqnarray*}
P^+_V(m)-P^-_V(m) &=& \left( c_1(\cm_1-\ker \psi )+c_1(\cm_3- \cm_2/\ker \psi )\right) \cap [Pic^m_V]\\
&=&
c_1( \cm_1-\cm_2+\cm_3) \cap [Pic^m_V]\\
&=& \shat(R\mu_! \nl ).
\end{eqnarray*}
\end{proof}
For surfaces $V$ with $p_g(V)>0$, we have no general result comparing $P^+_V(m)$ with $P^-_V(m)$. We expect that the following holds:
\begin{conj}
Let $V$ be a surface with $p_g(V)>0$. Then
\[
P^+_V(m)=P^-_V(m)
\]
for all $m \in H^2(V, \nz )$.
\end{conj}
For a conceptual explanation of this conjecture, we refer to section 6.
\subsection{Basic classes}
\begin{prop}
Let $V$ be a surface with $p_g(V)>0$, $m \in H^2(V,\mathbb{Z})$, and suppose that the fibered product
\[
\fp{m}
\]
is empty. Then
\[
\lb \Hilb^m_V \rb = \lb \Hilb^{k-m}_V \rb =0.
\]
\end{prop}
\begin{proof}
See Thm.~\ref{thm:compobth}
\end{proof}
\begin{dfn}
Let $V$ be a surface. A {\em basic class} of $V$ is an element $m \in H^2(V,\nz )$ with $(P^+_V(m),P^-_V(m) \neq (0,0)$. The surface $V$ is of {\em simple type} if every basic class $m$ satisfies $m(m-k)=0$.
\end{dfn}
\begin{prop}\label{prop:simp}
Every surface $V$ with $p_g(V)>0$ is of simple type and has only finitely many basic classes.
\end{prop}
\begin{proof}
Suppose first that the surface $V$ is minimal, and let $m \in H^2(V,\nz)$ be a basic class. Then the fibered product $\fp{m}$ is nonempty. Fix an element
\[
(D_1,D_2) \in \fp{m}
\]
The sum $K:=D_1+D_2$ is an effective canonical divisor. When $V$ is a $K3$ surface or abelian, then $K=0$ and hence $D_1=D_2=0$. In particular, $V$ has exactly one basic class and is of simple type.
When $V$ is properly elliptic we have
\[
D_1 \cdot D_2 \geq 0
\]
with equality holding iff there exists a rational number $0 \leq \lambda \leq 1$ with
\[
[D_1] = \lambda [K] \in H^2(V,\nq).
\]
This follows from the numerical effectivity of canonical divisors and the Hodge index theorem \cite{mbd}.
Since
\begin{eqnarray*}
m(m-k) &=& -D_1 \cdot D_2\\
&\geq &0
\end{eqnarray*}
we infer that $V$ is of simple type and has only finitely many basic classes.
When $V$ is of general type, canonical divisors are 1-connected \cite[PropVII.6.1]{bpv}. Therefore we have
\[
D_1 \cdot D_2 \geq 0
\]
with equality holding iff $D_1=0$ or $D_2=0$. Hence $V$ has exactly two basic classes, namely $0$ and $k$, and thus is of simple type.

Let now $V$ be a surface of simple type with finitely many basic classes. Let $\sigma: \hat{V} \to V$ be the blow-up of a point, denote by $E$ the exceptional curve and by $e$ the Poincar\'e dual of $[E]$. Fix a basic class $\hat{m} \in H^2(\hat{V},\nz)$. Then $\hat{m}$ can be uniquely written as
\[
\hat{m}=\sigma^*m + l\cdot e
\]
for a class $m \in H^2(V,\nz)$ and an integer $l \in\nz$. Moreover, by Thm.~\ref{thm:blowup}, $m$ is a basic class of $V$. Since by hypothesis $V$ is of simple type, we infer $l=0$ or $1$. Hence also $\hat{V}$ is of simple type and has only finitely many basic classes.
\end{proof}
\begin{prop}
A surface $V$ with $p_g(V)=0$ has infinitely many basic classes and is not of simple type. However, one has
\[
P^+_V(m)=0 \text{ or } P^-_V(m)=0
\]
unless $m(m-k)=0$.
\end{prop}
\begin{proof}
Fix an ample divisor $H \subset V$ and set $h:=c_1(\ov(H))$. Then there exists  $l_0 \in \nz$ such that for all integers $l$ with $l \geq l_0$ we have
\[
(l\cdot h) (l\cdot h -k) \geq 2 q(V).
\]
The wall-crossing formula yields
\[
P^+_V(l\cdot h) - P^-_V (l\cdot h)= [Pic^{l\cdot h}_V] \ + \ \text{terms of lower order}
\]
for all $l \geq l_0$. This proves the first claim.

The second claim is an immediate consequence of Cor.~\ref{cor:34}.
\end{proof}
\section{Examples}
\subsection{Ruled surfaces}
In this subsection we will compute the Poincar\'e invariants of ruled surfaces, and we will show how our methods yield easy proofs of classical results by Nagata \cite{na} and Lange \cite{la}.

To start, we observe that the wall crossing formula implies:
\begin{prop}\label{prop:pivanishing}
Let $V$ be a surface with $p_g(V)=q(V)=0$, and fix a class $m \in H^2(V,\nz )$ with $m(m-k)\geq 0$. Then
\[
(P^+_V(m),P^-_V(m)) = \left\{
\begin{array}{cl}
( 1,0) & \text{if $\Hilb^m_V \neq \emptyset$,}\\
(0,-1) & \text{if $\Hilb^m_V = \emptyset$.}
\end{array} \right.
\]
\end{prop}
\begin{proof}
When $p_g(V)=q(V)=0$, the wall-crossing formula says
\[
P^+_V(m)-P^-_V(m)=1.
\]
The claimed equality follows now from the fact that for every $m \in H^2(V,\nz)$ one of the Hilbert schemes $\Hilb^m_V$ or $\Hilb^{k-m}_V$ is empty.
\end{proof}
\begin{prop}\label{prop:pihirz}
Let $\nf_n \to \np^1$ be the $n$-th Hirzebruch surface, let $F$ be a fiber of the ruling, and choose a class $m\in H^2(V,\nz)$ with $m(m-k)\geq 0$. If $\langle m,[F] \rangle \geq 0$, then
\[
(P^+_V(m),P^-_V(m)) = (1,0).
\]
If $\langle m,[F] \rangle < 0$, then
\[
(P^+_V(m),P^-_V(m)) = (0,-1).
\]
\end{prop}
\begin{proof}
If $\langle m,[F] \rangle \geq 0$, then the adjunction formula yields $\langle k-m,[F] \rangle \leq -2$ and we infer $\Hilb^{k-m}_V = \emptyset$. If $\langle m,[F] \rangle < 0$, then we have $\Hilb^m_V = \emptyset$. Hence our claim is a direct consequence of Prop.~\ref{prop:pivanishing}.
\end{proof}
In order to compute the Poincar\'e invariants of ruled surfaces of irregularity $q\geq 1$, we restate the wall crossing in more accessible terms. This reformulation is of independent interest.

Let $V$ be a surface with $p_g(V)=0$ and $q(V)\geq 1$. Then $V$ admits a map $p: V \to C$ onto a smooth curve such that the induced morphism $p^*: H^1 (C, \nz) \to H^1(V,\nz )$ is an isomorphism: When $kod(V)=-\infty$, any minimal model $V \to V_{min}$ admits a unique geometric ruling $V_{min} \to C$, and we define $p$ to be the composition $V \to V_{min} \to C$. Note that, up to unique isomorphism, the map $p$ does not depend on the choice of a minimal model. When $kod(V) \geq 0$, then we have $q(V)=1$, and we define $p$ to be the Albanese mapping.

Let now $C$ be a smooth curve of genus $g$, and fix a natural number $d$ with $0 \leq d \leq g$. The Brill-Noether locus
\[
W_d:= \{ [\cL] \in Pic^d(C) | h^0(\cL) >0 \}
\]
carries the structure of a subscheme of $Pic^d(V)$ and hence possesses a fundamental class $[W_d] \in H_*(Pic^d(C),\mathbb{Z}) \cong H_*(Pic^0(C),\mathbb{Z})$ \cite{acgh}.
\begin{prop}\label{prop:wallbis}
Let $V$ be a surface with $p_g(V)=0$ and $q(V)\geq 1$, and denote by $F$ a general fiber of the map $p: V \to C$. Then, for every $m\in H^2(V,\nz )$ with $m(m-k)\geq 0$, we have
\[
P^+_V(m)-P^-_V(m)= \sum_{d=0}^{\min \{ q(V), \frac{m(m-k)}{2}\} } \left( \frac{\langle 2m-k, [F] \rangle}{2} \right)^{q(V)-d} [W_d].
\]
\end{prop}
\begin{proof}
By construction, the map $p: V \to C$ induces isomorphisms
\[
p^*: H^1(C,\nz) \stackrel{\cong}{\longrightarrow} H^1(V,\nz)
\]
and
\[
p^*: \Pic^0_C \stackrel{\cong}{\longrightarrow} \Pic^0_V.
\]
In particular, $q(V)=g(C)$, where $g(C)$ is the genus of the curve $C$. We compute:
\begin{eqnarray*}
\theta_{2m-k}(p^*(a) \wedge p^*(b)) &=& \frac{1}{2} \langle p^*(a) \cup p^*(b) \cup (2m-k) ,[V] \rangle\\
&=& \frac{1}{2} \langle p^*(a \cup b) , (2m-k) \cap [V] \rangle\\
&=& \frac{1}{2} \langle a \cup b, p_*( (2m-k) \cap [V]) \rangle\\
&=& \frac{\langle 2m-k,[F] \rangle}{2} \langle a \cup b,[C] \rangle\\
&=& \frac{\langle 2m-k,[F] \rangle}{2} \theta (a \wedge b)
\end{eqnarray*}
Our claim is now a consequence of Thm.~\ref{thm:wcfpoinc} and of the Poincar\'e formula, which asserts that
\[
[W_d] = \frac{\theta^{g(C)-d}}{(g(C)-d)!} [\Pic^d_C]
\]
for $0 \leq d \leq g(C)$.
\end{proof}
\begin{prop}\label{prop:rul}
Let $p: V \to C$ be a ruled surface over a curve of genus $g$, and let $F$ be a fiber of $p$. Fix a class $m \in H^2(V,\nz )$ with $m(m-k)\geq 0$. If $\langle m,[F] \rangle\geq -1$, then
\begin{eqnarray*}
P^+_V(m) &=& \sum_{d=0}^{\min \{g, \frac{m(m-k)}{2}\} } (\langle m,[F] \rangle+1)^{g-d} [W_d],\\
P^-_V(m)&=&0.
\end{eqnarray*}
If $\langle m,[F] \rangle\leq -1$, then
\begin{eqnarray*}
P^+_V(m) &=& 0,\\
P^-_V(m) &=& - \sum_{d=0}^{\min \{ g, \frac{m(m-k)}{2}\} } (\langle m,[F] \rangle+1)^{g-d} [W_d].
\end{eqnarray*}
\end{prop}
\begin{proof}
If $\langle m,[F] \rangle \geq -1$, then $\Hilb^{k-m}_V=\emptyset$, whereas $\langle m,[F] \rangle\leq -1$ implies that $\Hilb^m_V=\emptyset$. Therefore our claim is a consequence of Prop.~\ref{prop:wallbis}.
\end{proof}
Note that the above proposition yields a classical result of Nagata \cite{na}:
\begin{thm}[Nagata]
Let $p: V \to C$ be a geometrically ruled surface over a curve of genus $g$. Then there exists a section $s: C \to V$ with self-intersection number $s^2 \leq g$.
\end{thm}
\begin{proof}
Denote by $F$ a fiber of the ruling, and fix a class $m$ with $\langle m,[F] \rangle =1$. By adding or subtracting the Poincar\'e dual of $[F]$, we can modify $m$ in such a way that
\[
0 \leq \frac{m(m-k)}{2} \leq 1;
\]
then Prop.~\ref{prop:rul} says
\[
P^+_V(m) \neq 0.
\]
In particular, $\Hilb^m_V\neq \emptyset$. Choose a divisor $D \in \Hilb^m_V$, and let $D_0$ be the irreducible component of $D$ with $D_0\cdot F=1$. Then $D_0$ is a smooth curve of genus $g$. By adjunction we have
\begin{eqnarray*}
\frac{m(m+k)}{2} &= & \frac{D_0( D_0+K)}{2}\\
&=& g-1,
\end{eqnarray*}
where $K$ is a canonical divisor. This implies
\[
m^2 \leq g.
\]
Hence $D_0$ defines a section with self-intersection number $D_0^2 \leq m^2 \leq g$.
\end{proof}
For a geometrically ruled surface $p: V \to C$, put
\[
s(V): = \min \{ n \mid \text{$\exists$ a section $s$ with $s^2=n$} \}.
\]
The following proposition is a strengthening of a result of Lange \cite[Cor.5.3]{la}.
\begin{prop}
Let $p: V \to C$ be a geometrically ruled surface over a curve of genus $g$. Suppose either that $s(V)=g-1$ and the number of sections with selfintersection number $g-1$ is finite, or that $s(V)=g$ and the number of sections with selfintersection number $g$ that pass through a fixed point is finite. Then the length of the scheme parametrizing these sections is
\[
2^g.
\]
\end{prop}
\begin{proof}
By our assumption on the invariant $s$, every effective divisor $D \subset V$ of relative degree $1$ over $C$ and intersection number $D^2=g$ or $g-1$ respectively is irreducible and reduced, hence the graph of a section $C \to V$. Therefore, our claim is an immediate consequence of Prop.~\ref{prop:rul}.
\end{proof}
\subsection{Elliptic fibrations and logarithmic transformations}
\begin{lem}\label{lem:fib}
Let $\pi : V \to C$ be an elliptic fibration, denote by $F$ a general fiber and by $m_1F_1,\ldots m_rF_r$ the multiple fibres of $\pi$. For any element $m\in H^2(V,\mathbb{Z})$ with $m^2=\langle m,[F]\rangle =0$, there exists a canonical isomorphism
\[
\coprod_{\dfrac{d[F]+\sum a_i[F_i] =PD(m)}{0\leq a_i <m_i}} C_d \stackrel{\cong}{\longrightarrow} \Hilb^m_V.
\]
\end{lem}
\begin{proof}
Fix an element $D \in \Hilb^m_V$. The equality $D \cdot F=0$ shows that $D$ is contained in the fibres of $\pi$. Since $D^2=0$, Zariski's lemma implies that there exists an effective divisor $\mathfrak{d} \subset C$ and integers $0 \leq a_i <m_i$ for $i=1,\ldots,r$ such that
\[
D = \pi^* \mathfrak{d} + \sum a_i F_i.
\]
Hence the natural map
\[
\coprod_{\dfrac{d[F]+\sum a_i[F_i] =PD(m)}{0\leq a_i <m_i}} C_d \longrightarrow \Hilb^m_V,
\]
which sends an $(r+1)$-tuple $(\mathfrak{d},a_1, \ldots , a_r)$ to the divisor $\pi^* \mathfrak{d} + \sum a_i F_i$, is a bijection. That this map is also an isomorphism of schemes has been proved in \cite[Lemma 1.2.50]{mbd}.
\end{proof}
\begin{prop}\label{prop:pifib}
Let $\pi : V \to C$ be an elliptic fibration, denote by $F$ a general fiber and by $m_1F_1,\ldots m_rF_r$ the multiple fibres of $\pi$. Set $g:=g(C)$, and fix  $m\in H^2(V,\mathbb{Z})$ with $m^2=\langle m,[F]\rangle =0$. Then
\[
P^+_V(m)=\sum_{\dfrac{d[F]+\sum a_i[F_i] =PD(m)}{0\leq a_i <m_i}} (-1)^d 
\begin{pmatrix}
2g-2 + \chi (\ov )\\
d
\end{pmatrix},
\]
\[
P^-_V(m)=\sum_{\dfrac{d[F]+\sum a_i[F_i] =PD(k-m)}{0\leq a_i <m_i}} (-1)^{\chi (\ov)+d} 
\begin{pmatrix}
2g-2 + \chi (\ov )\\
d
\end{pmatrix}.
\]
\end{prop}
\begin{proof}
In \cite[p.473]{fm} it has been shown that
\[
c_{top}(R^1 \pi_* \co_\nd(\nd))|_{C_d} = (-1)^d 
\begin{pmatrix}
2g-2 + \chi (\ov )\\
d
\end{pmatrix}.
\]
Hence our claim is a consequence of the previous lemma.
\end{proof}
\begin{cor}\label{cor:ppluseqpminus}
Let $V$ be an elliptic surface with $p_g(V)>0$. Then
\[
P^+_V(m)=P^-_V(m)
\]
for all $m \in H^2(V, \nz)$.
\end{cor}
\begin{proof}
By Thm.~\ref{thm:blowup}, it suffices to give a proof in the minimal case.

Let $\pi : V \to C$ be a minimal elliptic surface over a curve of genus $g$, and suppose that $p_g(V)>0$. By the canonical bundle formula, there is an effective divisor $\mathfrak{d}$ on $C$ of degree $2g-2+ \chi (\ov )$ such that
\[
K:= \pi^* \mathfrak{d} + \sum_i (m_i-1)F_i
\]
is a canonical divisor on $V$. Fix a class $m \in H^2(V,\nz)$ with $m(m-k)\geq 0$. By Thm.~\ref{thm:compobth} we have $P^+_V(m)=P^-_V(m)=0$ whenever the fibered product $\Hilb^m_V \times_{\Pic^m_V} \Hilb^{k-m}_V$ is empty. Therefore we may suppose that there exists a decomposition $K'=D_1+D_2$ of a canonical divisor $K'$ into two effective divisors, such that $[D_1]$ is Poincar\'e dual to $m$. The inequality
\[
D_1\cdot D_2=-m(m-k)\leq 0
\]
implies $m^2= \langle m, [F]\rangle=0$. This can be seen as follows: since $V$ is minimal and of Kodaira dimension $0$ or $1$, any canonical divisor is numerically 0-connected, hence $D_1 \cdot D_2=0$. On the other hand, $K$ is numerically effective with $K^2=0$, hence $K\cdot D_i=0$ and $D^2_i=0$.  Note that $d[F]+\sum a_i[F_i]$ is Poincar\'e dual to $m$ if and only if $(2g-2+ \chi (\ov )-d) [F]+ \sum (m_i-1-a_i)[F_i]$ is Poincar\'e dual to $k-m$. Therefore we have
\begin{eqnarray*}
P^-_V(m) &=&  \sum_{\dfrac{d[F]+\sum a_i[F_i] =PD(k-m)}{0\leq a_i <m_i}}   (-1)^{2g-2+2 \chi (\ov) -d}
\begin{pmatrix}
2g-2 + \chi (\ov )\\
2g-2 + \chi (\ov )-d
\end{pmatrix}\\
&=& \sum_{\dfrac{d[F]+\sum a_i[F_i] =PD(m)}{0\leq a_i <m_i}} (-1)^d 
\begin{pmatrix}
2g-2 + \chi (\ov )\\
d
\end{pmatrix}\\
&=&P^+_V(m).
\end{eqnarray*}
\end{proof}
Using logarithmic transformations we will construct examples $(V,m)$, where $V$ is a surface with $p_g(V)=0$, and $m\in H^2(V,\nz)$ is a class such that neither $P_V^+(m)$ nor $P^-_V(m)$ vanishes.

Let $\Gamma$ be a lattice in $\nc$, and let $E := \nc/\Gamma$ be the
corresponding elliptic curve. We denote by $[z] \in E$ the point defined by $z \in \nc$. Let $t_1\in \np^1$ be a point. Choose a
positive integer $n_1$ and a complex number $\zeta_1$ such that $[\zeta_1]$ is a $n_1$-torsion point of $E$. Denote by $L_{t_1} (n_1,\zeta_1)(\np^1 \times E)$ the space obtained
by the logarithmic transformation $L_{t_1} (n_1,\zeta_1)$ from $\np^1
\times E$ \cite{eo}. Since a logarithmic transformation is a local analytic
construction one can apply further logarithmic transformations
$L_{t_2}(n_2, \zeta_2)$, \ldots, $L_{t_r}(n_r, \zeta_r)$ at points $t_2, \ldots, t_r \in \np^1$ such that $t_1, \ldots,t_r$ are pairwise distinct. We denote the  resulting space by $L_\tund(\nund,
\zetaund) (\np^1 \times E)$, where $\tund := (t_1, \ldots, t_r)$,
$\nund := (n_1, \ldots, n_r)$, and $\zetaund := (\zeta_1, \ldots,
\zeta_r)$. Note that $\LPN$ is a smooth compact complex surface, but not necessarily algebraic.
\begin{lem} \label{lem:log}
Let $\Gamma=\langle 1, \omega \rangle$ be a lattice in $\nc$, and fix $r$ distinct points $t_1,\ldots,t_r \in \np^1$. Choose integers $n_i$, $u_i$, $v_i$ for $i=1,\ldots, r$ such that\\
$gcd(n_i,u_i,v_i)=1$ for all $i$, and set $\zeta_i:=\frac{u_i+v_i\omega}{n_i}$.
\begin{enumerate}
\item The surface $\LPN$ is projective if and only if $\sum_i
  \zeta_i =0$.

\item Denote by $n_iF_i$ the multiple fibres, and by $F$ a regular fibre. If $\sum_i \zeta_i =0$, then
\[
H_2(\LPN,\nz) \cong \nz \oplus 
\left< [F],[F_1],\ldots , [F_r] \ \left| \
\begin{array}{l} 
n_i[F_i]=[F], \\ 
u_1[F_1]+\ldots +u_r[F_r]=0,\\ 
v_1[F_1]+\ldots +v_r[F_r]=0
\end{array} \right> \right. .
\]
\end{enumerate}
\end{lem}
\begin{proof}
For the first claim see \cite[p.284]{eo}, for the second claim see \cite[Thm.A.2.11]{mbd}.
\end{proof}
\begin{prop}\label{prop:log}
Let $\Gamma=\langle 1, \omega \rangle$ be a lattice in $\nc$, and fix $r$ distinct points $t_1,\ldots,t_r \in \np^1$. Choose integers $n_i$, $u_i$, $v_i$ for $i=1,\ldots, r$ such that\\
$gcd(n_i,u_i,v_i)=1$ for all $i$, and set $\zeta_i:=\frac{u_i+v_i\omega}{n_i}$. Suppose $\sum_i \zeta_i=0$, and denote by $n_iF_i$ the multiple fibres, and by $F$ a regular fibre. Let $P \subset \Pic^0(\LPN)$ be the subgroup generated by the classes of the line bundles $\co (\sum u_iF_i)$ and $\co (\sum v_iF_i)$. Then for all integers $d$, $a_1,\ldots a_r$, we have
\[
\Hilb^m_{\LPN} \cong \coprod_{[\cL ] \in P} | \ov(dF + \sum_i a_iF_i ) \otimes \cL|,
\]
where $m$ is the Poincar\'e dual of the class $d [F] + \sum a_i[F_i]$.
\end{prop}
\begin{proof}
This is a direct consequence of Lemma \ref{lem:fib} and Lemma \ref{lem:log}.
\end{proof}
We will also need the following

\begin{lem} \label{lem:alb}
Let $\Gamma=\langle 1, \omega \rangle$ be a lattice in $\nc$, and fix $r$ distinct points $t_1,\ldots,t_r \in \np^1$. Choose integers $n_i$, $u_i$, $v_i$ for $i=1,\ldots, r$ such that\\
$gcd(n_i,u_i,v_i)=1$ for all $i$, and set $\zeta_i:=\frac{u_i+v_i\omega}{n_i}$. Let $\Gamma' \subset \nc$ be the lattice generated by $1,\omega,\zeta_1,\ldots,\zeta_r$. Let $F$ be a regular fibre of the fibration $\LPN \to \np^1$, and assume that $\sum \zeta_i=0$. Then there exists an isomorphism 
\[
\Alb (\LPN) \stackrel{\cong}{\longrightarrow} \nc /\Gamma'
\]
such that the following diagram commutes:
\[
\xymatrix{
F \ar[r] \ar@1{=}[d] & \Alb (\LPN) \ar[d]^{\cong}\\
\nc/\Gamma \ar[r] & \nc/\Gamma'}
\]
\end{lem}
\begin{proof}
This follows from the explicit description of the fundamental group of $\LPN$, see for instance \cite[p.284]{eo}.
\end{proof}
\begin{rem}
Let $\Gamma' \subset \nc$ be the lattice generated by $1,\omega,\zeta_1,\ldots,\zeta_r$. If $\sum \zeta_i=0$, then $P$ is isomorphic to $\Gamma'/\Gamma$ as an abstract group.
\end{rem}
\begin{proof}
Note first that the group $P$ is the kernel of the restriction map
\[
\Pic^0_{\LPN} \longrightarrow \Pic^0_F.
\]
To see this, consider an element $\cL$ in this kernel. Then $\cL(nF)$ admits global sections for sufficiently large integers $n$. Hence $\cL$ is isomorphic to a line bundle of the form $\co( D-nF)$, where $D$ is an effective divisor contained in the fibers of the projection $\LPN \to \np^1$. Since the homology class $[D-nF]$ vanishes, we get $\co (D-nF) \cong \co(\sum u_iF_i)^{\otimes a} \otimes \co(\sum v_iF_i)^{\otimes b}$ for suitable integers $a,b$.
On the other hand, the group $\Gamma'/\Gamma$ is the kernel of 
\[
F= \Alb(F) \longrightarrow \Alb(\LPN).
\]
The claim follows now from the fact that Albanese and the Picard variety are dual tori.
\end{proof}
\begin{exa}
Choose 4 distinct points $t_1, \ldots , t_4 \in \np^1$, set $\nund=(3,3,3,3)$,
\[
\left(
\begin{array}{c}
\underline{u} \\
\underline{v}
\end{array}
\right)
=
\left(
\begin{array}{cccc}
1 & 1 & 1 & -3 \\
1 & 0 & 0 & -1
\end{array}
\right) ,
\]
and put $V:= \LPN$. Clearly $\Hilb^0_V =\{0\}$, and therefore
$P^+_V(0)=1$. To compute $P^-_V(0)$ we determine $\Hilb^k_V$, where $k$ denotes
$c_1({\KK}_V)$. Prop.~\ref{prop:log} implies that 
\begin{eqnarray*}
\Hilb^k_V & = & |-2F+2F_1+2F_2+2F_3+2F_4|  \cup  |2F_4 | \\
&& \cup \, | -F + F_1 + F_2 + F_3 + 2F_4 | \cup | -F + 2F_2 + 2F_3 +F_4|\\
&& \cup \, | F_1 + F_4 |  \cup  |-F +
2F_2 + F_2 + F_3 + F_4 | \\
&& \cup \, | -F + F_1 + 2F_2 +2F_3|  \cup  |2F_2|  \cup  | F_2 + F_3 |.
\end{eqnarray*}
Thus $\Hilb^k_V$ consists of 4 smooth points, namely $2F_4$, $F_1+F_4$, $2F_2$ and $F_2+F_3$, and we obtain $P^-_V(0)=4$. Of course, we can also use the wall crossing formula to compute the difference $P^+_V(0)-P^-_V(0)=-3$. 

Let $E$ be a fiber of the Albanese mapping $V \to \Alb(V)$. By Lemma \ref{lem:alb} we find $E \cdot F=9$. Since for a canonical divisor $K$ we have
\[
[K] = \frac{2}{3} [F] \in H_2(V,\nq),
\]
Prop.~\ref{prop:wallbis} yields
\[
P^+_V(0)-P^-_V(0) = \frac{1}{2} \cdot \frac{-2}{3} \cdot 9 
=  -3. 
\]
\end{exa}
\section{Comparison with Seiberg-Witten invariants}
\subsection{Three conjectures}
In this section, we will compare our Poincar\'e invariants with the full Seiberg-Witten invariants. The latter are differential-topological invariants, which were defined in \cite{ot1} and refine the invariants introduced by Seiberg and Witten \cite{w}. We briefly recall the structure of the full Seiberg-Witten invariants; for the construction and details, we refer to \cite{ot1}.

Let $(M,g)$ be a closed oriented Riemannian 4-manifold with first Betti number $b_1$. We denote by $b_+$ the dimension of a maximal subspace of $H^2(M,\nr)$ on which the intersection form is positive definite. Recall that the set of isomorphism classes of $\spin$-structures on $(M,g)$ has the structure of a $H^2(M,\mathbb{Z})$-torsor. This torsor does, up to a canonical isomorphism, not depend on the choice of the metric $g$ and will be denoted by $\spim$.
 
We have the Chern class mapping
\begin{eqnarray*}
c_1:\spim & \longrightarrow & H^2(M,\mathbb{Z}) \\
\mathfrak{c} & \longmapsto & c_1(\mathfrak{c}),
\end{eqnarray*}
whose image consists of all characteristic elements.

If $b_+>1$, then the Seiberg-Witten invariants are maps
\[
SW_{M,\ooo} : \spim \longrightarrow \Lambda^*H^1(M,\mathbb{Z}),
\]
where $\oo$ is an orientation parameter.

When $b_+=1$, then the invariants depend on a chamber structure and are maps
\[
SW^\pm_{M,\ooov} : \spim \longrightarrow \Lambda^*H^1(M,\mathbb{Z}) \times\Lambda^*H^1(M,\mathbb{Z}),
\]
where $\oov$ are again orientation data. The difference of the two components is a purely topological invariant. More precisely, we have:
\begin{thm}[Okonek/Teleman]\label{thm:wcfsw}
Let $M$ be a closed connected oriented 4-manifold with $b_+=1$. Fix an orientation $\scriptstyle{\mathcal{O}_1}$ of $H^1(M,\mathbb{R})$, and denote by $l_{\scriptscriptstyle{\mathcal{O}_1}}\in \Lambda^{b_1} H^1(M,\mathbb{Z})$ the generator defining the orientation $\scriptstyle{\mathcal{O}_1}$. For every class $\mathfrak{c}$ of $\spin$-structures of Chern class $c$, the following holds:
\[
\left( \swwm{+}{\mathfrak{c}}-\swwm{-}{\mathfrak{c}} \right)(\lambda) = \frac{1}{\left[ \frac{b_1-r}{2} \right] !}\left\langle \lambda\wedge {\theta_c}^{\left[ \frac{b_1-r}{2} \right]}, l_{\scriptscriptstyle{\mathcal{O}_1}} \right\rangle,
\]
where $r$ is an integer with $0\leq r \leq \min (b_1,\omega_c)$, $r\equiv b_1 \mod{2}$, and $\lambda \in \Lambda^r\left( H_1(M,\mathbb{Z})/Tors \right)$.
\end{thm}
\begin{proof}
\cite[Thm.16]{ot1}.
\end{proof}
\begin{rem}
The original formula in \cite[Thm. 16]{ot1}) contains an incorrect sign, which was detected in \cite{mbd}. The error occurs on page 821, where the authors do not take into account that the cohomology ring of a manifold is graded commutative. The error is of a purely calculatory nature and does not affect the rest of the proof.
\end{rem}

Let now $V$ be a surface. Any {\em Hermitian} metric $g$ on $V$ defines a {\em canonical} $\spin$-structure on $(V,g)$. Its class $\mathfrak{c}_{can} \in \spiv$ does not depend on the choice of the metric. The Chern class of $\mathfrak{c}_{can}$ is $c_1(\mathfrak{c}_{can})=-c_1(\kv )=-k$.

Since $\spiv$ is a $H^2(V,\nz)$-torsor, the distinguished element $\mathfrak{c}_{can}$ defines a bijection:
\begin{eqnarray*}
H^2(V,\nz) & \longrightarrow & \spiv\\
m & \longmapsto & \mathfrak{c}_m
\end{eqnarray*}
The Chern class of the twisted structure $\mathfrak{c}_m$ is $2m-k$. Finally, recall that any surface defines canonical orientation data $\oo$ and $\oov$ respectively.
\begin{conj}\label{conj:main}
Let $V$ be a surface, and denote by $\oo$ or $\oov$ the canonical orientation data. If $p_g(V)=0$, then
\[
P^\pm_V(m)=\swwv{\pm}{\mathfrak{c}_m}\ \ \forall m\in H^2(V,\mathbb{Z}).
\]
If $p_g(V)>0$, then
\[
P^+_V(m)=P^-_V(m)=\swv{\mathfrak{c}_m}\ \ \forall m\in H^2(V,\mathbb{Z}).
\]
\end{conj}

The main evidence for our conjecture comes from the following Kobayashi-Hitchin correspondence:
\begin{thm}[Okonek/Teleman]\label{thm:kh}
Let $(V,g)$ be a surface endowed with a K\"ahler metric $g$. Fix a class $m \in H^2(V,\nz )$ and a real closed 2-form $\beta$ of type $(1,1)$. Let $\tau$ be a $\spin$-structure on $(V,g)$ representing the class $\mathfrak{c}_m$, and denote by $\mathcal{W^\tau_\beta}$ the moduli space of solutions to the $\beta$-twisted Seiberg-Witten equations.
\begin{itemize}
	\item[i)] If $(2m-k-[\beta]) \cdot [\omega_g] <0$, then there exists an isomorphism of real analytic spaces
\[
\kappa^+_m: \mathcal{W^\tau_\beta}\stackrel{\cong}{\longrightarrow} \Hilb^m_V.
\]
\item[ii)] If $(2m-k-[\beta]) \cdot [\omega_g] >0$, then there exists an isomorphism of real analytic spaces
\[
\kappa^-_m: \mathcal{W^\tau_\beta}\stackrel{\cong}{\longrightarrow} \Hilb^{k-m}_V.
\]
\end{itemize}
\end{thm}
\begin{proof}
\cite[Thm.25]{ot1}.
\end{proof}

By the work of Brussee \cite{br}, the moduli space of solutions to the Seiberg-Witten equations carries a {\em virtual fundamental class} $[\mathcal{W^\tau_\beta}]_{vir}$. Moreover, the full Seiberg-Witten invariants can be computed by evaluating tautological cohomology classes on $[\mathcal{W^\tau_\beta}]_{vir}$ \cite{ot2}.
Our main conjecture is essentially a consequence of the following more conceptual conjecture:
\begin{conj}\label{conj:precise}
Let $(V,g)$ be a surface endowed with a K\"ahler metric $g$. Fix a class $m \in H^2(V,\nz )$ and a real closed 2-form $\beta$ of type $(1,1)$. Let $\tau$ be a $\spin$-structure on $(V,g)$ representing the class $\mathfrak{c}_m$, and denote by $\mathcal{W^\tau_\beta}$ the moduli space of solutions to the $\beta$-twisted Seiberg-Witten equations. Choose the canonical orientation data $\oo$ or $\oov$. Suppose that $(2m-k-[\beta]) \cdot [\omega_g] <0$. Then the Kobayashi-Hitchin isomorphism
\[
\kappa^+_m: \mathcal{W^\tau_\beta}\stackrel{\cong}{\longrightarrow} \Hilb^m_V
\]
identifies $[\mathcal{W^\tau_\beta}]_{vir}$ with the image of $\lb \Hilb^m_V \rb$ in $H_*(\Hilb^m_V,\mathbb{Z})$.
\end{conj}
The tautological cohomology classes on $\mathcal{W}^\tau_\beta$ are given by a canonical map
\[
r: \left( \Lambda^* H^1(V,\nz)^\vee \right)[u] \longrightarrow H^*(\mathcal{W}^\tau_\beta,\nz),
\]
where $u$ is a class of degree 2.
\begin{lem}\label{lem:com}
Suppose that $(2m-k-[\beta]) \cdot [\omega_g] <0$. Then the following diagram commutes:
\[
\xymatrix{
H^*(\Pic^m_V, \nz )[u] \ar@1{=}[r]
\ar[d]_{(\rho^+)^*} 
& \left( \Lambda^* H^1(V,\nz)^\vee \right)[u] \ar[d]^r\\
H^*(\Hilb^m_V,\nz) \ar[r]_{(\kappa^+_m)^*} & H^*(\mathcal{W}^\tau_\beta,\nz )}
\]
where $(\rho^+)^*(u):=u^+$.
\end{lem}
\begin{proof}
\cite{dt}.
\end{proof}
Combining this lemma with Conj.~\ref{conj:precise} yields immediately
\[
P^+_V(m) = \swwv{+}{\mathfrak{c}_m}
\]
and
\[
P^+_V(m) = \swv{\mathfrak{c}_m}
\]
respectively.

The second case, i.e.~when $(2m-k-[\beta])\cdot [\omega_g]>0$, can be reduced to the first by the following trick \cite{ta}: By complex conjugation, every $\spin$-structure $\tau$ gives rise to a dual structure $\tau^*$ with the following properties:
\begin{itemize}
	\item[-] If $\tau$ represents the class $\mathfrak{c}_m$, then $\tau^*$ represents the class $\mathfrak{c}_{k-m}$.
	\item[-] For every real closed 2-form $\beta$ there is a canonical isomorphism
\[
\zeta: \mathcal{W}^\tau_\beta \stackrel{\cong}{\longrightarrow} \mathcal{W}^{\tau^*}_{-\beta}.
\]
\end{itemize}
Moreover, $\zeta$ maps $[\mathcal{W}^\tau_\beta]_{vir}$ to $(-1)^{\chi (\ov) +\frac{m(m-k)}{2}}[\mathcal{W}^{\tau^*}_{-\beta}]_{vir}$, and the following diagram commutes
\[
\xymatrix{
\left( \Lambda^* H^1(V,\nz)^\vee \right)[u] \ar[d]_\gamma \ar[r] & H^*(\mathcal{W}^{\tau^*}_{-\beta},\nz) \ar[d]^{\zeta^*}\\
\left( \Lambda^* H^1(V,\nz)^\vee \right)[u] \ar[r]&H^*(\mathcal{W}^\tau_\beta,\nz)}
\]
where $\gamma$ maps an element $t\in H^1(V,\nz)^\vee$ to $-t$ and $u$ to $-u$.

Using this trick we can show that Conj.~\ref{conj:precise} implies that the Kobayashi-Hitchin isomorphism $\kappa^-_{m}$ identifies $[\mathcal{W}^{\tau}_{\beta}]_{vir}$ with the image of \newline
$(-1)^{\chi(\ov ) +\frac{m(m-k)}{2}} \lb \Hilb^{k-m}_V \rb$ in $H_*(\Hilb^{k-m}_V, \nz)$ as follows: We have

\[
(-1)^{\chi(\ov ) +\frac{m(m-k)}{2}} \zeta_*[\mathcal{W}^{\tau^*}_{-\beta}]_{vir}=[\mathcal{W}^\tau_\beta ]_{vir},
\]
and, by assumption, the Kobayashi-Hitchin isomorphism $\kappa^+_{k-m}$ identifies $[\mathcal{W}^{\tau^*}_{-\beta}]_{vir}$ with the image of $\lb \Hilb^{k-m}_V \rb$ in $H_*(\Hilb^{k-m}_V, \nz)$. Since $\kappa^-_m=\kappa^+_{k-m}\circ \zeta$, our claim follows. 

Using again the identity $\kappa^-_m=\kappa^+_{k-m}\circ \zeta$, Lemma \ref{lem:com} shows that also the following diagram commutes:

\[
\xymatrix{
H^*(\Pic^m_V, \nz )[u] \ar@1{=}[r]
\ar[d]_{(\rho^-)^*} 
& \left( \Lambda^* H^1(V,\nz)^\vee \right)[u] \ar[d]^r\\
H^*(\Hilb^{k-m}_V,\nz) \ar[r]_{(\kappa^-_m)^*} & H^*(\mathcal{W}^\tau_\beta,\nz )}
\]
where $(\rho^-)^*(u):=-u^-$.
This commutative diagram and the identification of $[\mathcal{W}^{\tau}_{\beta}]_{vir}$ with the image of $(-1)^{\chi (\ov) +\frac{m(m-k)}{2}} \lb \Hilb^{k-m}_V \rb$ in $H_*(\Hilb^{k-m}_V, \nz)$ under the Kobayashi-Hitchin isomorphism $\kappa^-_m$ yield at once
\[
P^-_V(m) = \swwv{-}{\mathfrak{c}_m}
\]
and
\[
P^-_V(m) = \swv{\mathfrak{c}_m}
\]
respectively.

There is one case, in which Conj.~\ref{conj:precise} is known to hold:
\begin{thm}[D\"urr/Teleman]\label{thm:andreimbd}
Let $(V,g)$ be a surface endowed with a K\"ahler metric $g$. Fix a class $m \in H^2(V,\nz )$ and a real closed 2-form $\beta$ of type $(1,1)$. Let $\tau$ be a $\spin$-structure on $(V,g)$ representing the class $\mathfrak{c}_m$, and denote by $\mathcal{W^\tau_\beta}$ the moduli space of solutions to the $\beta$-twisted Seiberg-Witten equations. Choose the canonical orientation data $\oo$ or $\oov$.

i) If $(2m-k-[\beta]) \cdot [\omega_g] <0$ and the moduli space $\mathcal{W^\tau_\beta}$ is smooth, then the Kobayashi-Hitchin isomorphism
\[
\kappa^+_m: \mathcal{W^\tau_\beta}\stackrel{\cong}{\longrightarrow} \Hilb^m_V
\]
identifies $[\mathcal{W^\tau_\beta}]_{vir}$ with the image of $\lb \Hilb^m_V \rb$ in $H_*(\Hilb^m_V,\mathbb{Z})$.

ii) If $(2m-k-[\beta]) \cdot [\omega_g] >0$ and the moduli space $\mathcal{W^\tau_\beta}$ is smooth, then the Kobayashi-Hitchin isomorphism
\[
\kappa^-_m: \mathcal{W^\tau_\beta}\stackrel{\cong}{\longrightarrow} \Hilb^{k-m}_V
\]
identifies $[\mathcal{W^\tau_\beta}]_{vir}$ with the image of $(-1)^{\chi (\ov ) +\frac{m(m-k)}{2}}\lb \Hilb^{k-m}_V \rb$ in\\
 $H_*(\Hilb^{k-m}_V,\mathbb{Z})$.
\end{thm}
\begin{proof}
\cite{dt}.
\end{proof}
\begin{cor}
Let $V$ be a surface with $q(V)=0$, and denote by $\oo$ or $\oov$ the canonical orientation data. If $p_g(V)=0$, then
\[
P^\pm_V(m)=\swwv{\pm}{\mathfrak{c}_m}\ \ \forall m\in H^2(V,\mathbb{Z}).
\]
If $p_g(V)>0$, then
\[
P^+_V(m)=P^-_V(m)=\swv{\mathfrak{c}_m}\ \ \forall m\in H^2(V,\mathbb{Z}).
\]
\end{cor}
\begin{proof}
The relevant moduli spaces are isomorphic to projective spaces, hence smooth.
\end{proof}
We denote by $\alpha$ the map
\begin{eqnarray*}
\alpha: \fp{m} & \longrightarrow & |\kv |\\
(D_1,D_2) & \longmapsto & D_1+D_2.
\end{eqnarray*}
\begin{thm}[Witten]\label{thm:witten}
Let $V$ be a surface with $p_g(V)>0$, $K$ an effective canonical divisor, and  fix $m\in H^2(V,\mathbb{Z})$. Then
\[
\swv{\mathfrak{c}_m}=\sum_{(D_1,D_2)\in \alpha^{-1}(K)} (-1)^{h^0(\co_{D_1}(D_1))} l(D_1,D_2),
\]
where $l(D_1,D_2)$ is the length of the local ring of the fibre $\alpha^{-1}(K)$ at the point $(D_1,D_2)$.
\end{thm}
\begin{proof}
A complete proof can be found in \cite{mbd}.
\end{proof}
This theorem is a refined version of Witten's trick \cite[p.787]{w}, which allows to compute the Seiberg-Witten invariants even when the relevant Hilbert schemes are oversized. The question arises, if there exists an algebro-geometric analogue. The strongest possible assertion one can hope for is the following
\begin{conj}\label{conj:trick}
Let $V$ be a surface with $p_g(V)>0$, choose an effective canonical divisor $K$, and fix $m \in H^2(V,\mathbb{Z})$. Set
\[
C(m,K):= \sum_{(D_1,D_2)\in \alpha^{-1}(K)} (-1)^{h^0(\co_{D_1}(D_1))} l(D_1,D_2)[D_1,D_2],
\]
where
\[
[D_1,D_2] \in A_0(\fp{m})
\]
denotes the class of the point $(D_1,D_2)$. Let $p_1$ and $p_2$ be the projections from $\fp{m}$ to $\Hilb^m_V$ and $\Hilb^{k-m}_V$ respectively. Then
\[
\lb \Hilb^m_V \rb ={p_1}_*C(m,K)
\]
and
\[
 \lb \Hilb^{k-m}_V \rb =(-1)^{\chi (\ov) +\frac{m(m-k)}{2}} {p_2}_*C(m,K).
\]
\end{conj}
For this conjecture to make sense, it is clearly necessary that the images of the cycle class $C(m,K)$ in the Chow groups of the Hilbert schemes $\Hilb^m_V$ and $\Hilb^{k-m}_V$ do not depend on the choice of the canonical divisor $K$. Indeed, it is possible to show the following stronger result:
\begin{prop}\label{prop:indep}
Let $V$ be a surface with $p_g(V)>0$, choose an effective canonical divisor $K$, and fix $m \in H^2(V,\mathbb{Z})$. The class
\[
C(m,K) \in A_0(\fp{m})
\]
does not depend on the choice of $K \in |\kv|$.
\end{prop}
{\em Sketch of proof:} The argument has two parts: First one shows that the map
\[
\alpha : \fp{m} \to |\kv|
\]
is flat of relative dimension $0$ when the fibered product $\fp{m}$ is non-empty and $m(m-k)=0$. This implies that
\[
\sum_{(D_1,D_2)\in \alpha^{-1}(K)} l(D_1,D_2)[D_1,D_2]
\]
is independent of $K$. The second point is to show that the sign $(-1)^{h^0(\co_{D_1}(D_1))}$ is constant on every connected component of $\fp{m}$. This can be done by a case by case analysis according to the Kodaira dimension.
\mbox \nolinebreak \hfill $\Box$ \medbreak \par

Now we prove that Conj.~\ref{conj:trick} is true in the smooth case:
\begin{prop}
Let $V$ be a surface with $p_g(V)>0$, choose an effective canonical divisor $K$, and fix $m\in H^2(V,\mathbb{Z})$. If $\Hilb^m_V$ is smooth, then
\[
\lb \Hilb^m_V \rb ={p_1}_*C(m,K).
\]
If $\Hilb^{k-m}_V$ is smooth, then
\[
\lb \Hilb^{k-m}_V \rb =(-1)^{\chi(\ov ) +\frac{m(m-k)}{2}}{p_2}_*C(m,K).
\] 
\end{prop}
\begin{proof}
Fix a form $\eta \in H^0(\kv)\setminus \{ 0 \}$ and set $K:=(\eta)$. Recall that $\pi: \Hilb^m_V \times V \to \Hilb^m_V$ and $pr: \Hilb^m_V \times V \to V$ are the projections. Using the restriction morphism
\[
\co_{\Hilb^m_V}\otimes H^0(\kv) \cong \pi_*(pr^*\kv) \longrightarrow \pi_* (pr^*\kv \otimes \co_\nd),
\]
the form $\eta$ defines a section of the coherent sheaf $\pi_* (pr^*\kv \otimes \co_\nd)$, which vanishes exactly at the divisors $D \in \Hilb^m_V$ with $D \leq K$. By relative duality there exists an isomorphism 
\[
\left( R^1\pi_* \co_\nd (\nd)\right)^\vee \stackrel{\cong}{\longrightarrow} \pi_* (pr^*\kv \otimes \co_\nd).
\]
When $\Hilb^m_V$ is smooth, then $R^1\pi_* \co_\nd (\nd)$ is locally free and the virtual fundamental class is given by the formula
\[
\lb \Hilb^m_V \rb = c_{top} (R^1\pi_* \co_\nd (\nd)) \cap [\Hilb^m_V].
\]
The claim follows now since for any locally free sheaf $\EE$ we have
\[
c_i(\EE^\vee) =(-1)^ic_i(\EE).
\]

To prove the second claim, we have to show that for any decomposition $K=D_1+D_2$ with $D_1 \cdot D_2=0$ we have
\begin{equation}\label{eq:vorzeichen}
h^1(\co_{D_1}(D_1))+h^1(\co_{D_2}(D_2)) \equiv \chi (\ov) \mod 2.
\end{equation}
First we reduce to the minimal case: Let $ \sigma : V \to V_{min}$ be the minimal model, fix a canonical divisor $K \geq 0$ on $V$ and let $K= D_1+D_2$ be a decomposition with $D_1 \cdot D_2=0$. Then $\sigma_!K$ is a canonical divisor on $V_{min}$ with decomposition $\sigma_!K = \sigma_!D_1+\sigma_!D_2$. Moreover, we have $\sigma_!D_1\cdot \sigma_!D_2=0$ and $h^1(\co_{\sigma_!D_i}(\sigma_!D_i))=h^1(\co_{D_i}(D_i))$ for $i=1,2$. Since $\chi(\ov)=\chi(\co_{V_{min}})$, we may assume that $V$ is minimal.

Suppose first that $V$ is a $K3$-surface or an abelian variety. Then $K=0$ and equation \eqref{eq:vorzeichen} holds since
\[
\chi (\ov) \equiv 0 \mod 2.
\]
Assume now that $V$ is properly elliptic, denote by $\varphi : V \to C$ the fibration, and by $m_1F_1, \ldots,m_rF_r$ the multiple fibers. Fix a canonical divisor $K$ and a decomposition $K=D_1+D_2$ with $D_1\cdot D_2=0$. Then there are effective divisors $\mathfrak{d}_1$, $\mathfrak{d}_2$ on $C$ and integers $0 \leq a_i <m_i$ for $i=1,\ldots,r$ with $D_1=\varphi^*(\mathfrak{d}_1)+\sum a_iF_i$, $D_2=\varphi^*(\mathfrak{d}_2)+\sum (m_i-a_i)F_i$. Set $d_i:= \deg \mathfrak{d}_i$. We have
\begin{eqnarray*}
h^1(\co_{D_1}(D_1))+ h^1(\co_{D_2}(D_2)) &=& d_1+d_2\\
&=& \chi(\ov)+2g(C)-2.
\end{eqnarray*}

Finally, let $V$ be a minimal surface of general type. If $K=D_1+D_2$ is a decomposition of a canonical divisor with $D_1\cdot D_2=0$, then either $D_1=0$ or $D_2=0$. Hence we have to show that
\[
h^1(\co_K(K)) \equiv \chi (\ov) \mod 2
\]
for all effective canonical divisors $K$. Given an effective canonical divisor$K$, choose a form $\eta$ with $K=(\eta)$, and denote by
\[
\eta \cdot: H^1(\ov) \to H^1(\kv)
\]
the multiplication map. Let
\[
\langle . , . \rangle : H^1(\ov) \times H^1(\kv) \to \nc
\]
be the Serre duality pairing. Since $\eta$ is a form of type $(2,0)$, the induced pairing
\begin{eqnarray*}
\left( H^1(\ov) /\ker \eta \cdot \right) \times \left( H^1(\ov) /\ker \eta \cdot \right) &\to& \nc\\
([\alpha],[\beta]) & \mapsto & \langle \alpha, \eta \cdot \beta \rangle
\end{eqnarray*}
is well-defined, non-degenerate and skew-symmetric. This yields
\[
\dim  \left( H^1(\ov) /\ker \eta \cdot \right) \equiv 0 \mod 2
\]
and 
\[
h^1(\co_K(K)) \equiv \chi (\ov) \mod 2.
\]
\end{proof}
As a simple application of this result one gets a new proof of Prop.~\ref{prop:pifib} in the case $p_g(V)>0$ (compare \cite[Thm.1.2.51]{mbd}).
Note that whenever (a homological version of) Conjecture \ref{conj:trick} holds, then
\[
P^+_V=P^-_V,
\]
and, using Thm. ~\ref{thm:witten},
\[
\swv{\mathfrak{c}_m}=P^\pm_V(m)
\]
for all $m \in H^2(V,\mathbb{Z})$. We do not know if there is a direct proof of the expected identity $P^+_V=P^-_V$, independent of Conj.~\ref{conj:trick}. Note that by Cor.~\ref{cor:ppluseqpminus}, $P^+_V=P^-_V$ holds for all elliptic surfaces with $p_g(V)>0$.

Note that that the three conjectures have a completely different character: While Conj.~\ref{conj:main} could be proved through a case by case analysis, Conj.~\ref{conj:precise} is an instance of a very general principle relating virtual fundamental classes in gauge theory and complex geometry:
{\em Let
\[
KH: \MM_{gauge} \longrightarrow \MM_{complex}
\]
be a Kobayashi-Hitchin correspondence between a gauge theoretical moduli space and a complex geometric moduli space. 
Suppose $\MM_{gauge}$ is the zero locus of a Fredholm section in a Banach bundle over a Banach manifold, and all data involved in the definition of $\MM_{complex}$ are algebraic. Then $\MM_{complex}$ has a preferred perfect obstruction theory, and the Kobayashi-Hitchin correspondence $KH$ maps $[\MM_{gauge}]_{vir}$ to the image of $[[\MM_{complex} ]]$ in Borel-Moore homology.}

For a proof of this general principle in special cases, see \cite{ot2}, \cite{ot3}.
The third conjecture, on the other hand, is of purely algebro-geometric nature. There is a unifying algebraic concept \cite{d2}, the Witten triples, which allows to relate Hilbert schemes of curves on surfaces with sets of decompositions of effective canonical divisors.

\begin{dfn}
Let $V$ be a surface, and fix $m \in H^2(V,\nz )$. A Witten triple of class $m$ is a triple $(\cL , \varphi, \psi )$ consisting of an invertible sheaf $\cL$ with $c_1(\cL )=m$, a morphism $\varphi: \ov \to \cL$, and a morphism $\psi : \cL \to \kv$. Two Witten triples $(\cL , \varphi, \psi )$ and $(\cL' , \varphi', \psi' )$ are equivalent, if there exists an isomorphism $\chi : \cL \to \cL'$ such that the following diagram commutes:
\[
\xymatrix{
\ov \ar@1{=}[d] \ar[r]^\varphi & \cL \ar[d]^\chi \ar[r]^\psi & \kv \ar@1{=}[d] \\
\ov \ar[r]^{\varphi'} & \cL' \ar[r]^{\psi'}& \kv}
\] 
\end{dfn}
For every ample class $h \in H^2(V,\nz )$ and every real number $t\in \nr$, one has a natural stability concept.
\begin{dfn}
A Witten triple $(\cL , \varphi, \psi )$ is {\em $t$-stable} on $(V,h)$, if one of the following three conditions is fullfilled:
\begin{itemize}
	\item[-] $\varphi \neq 0$ and $\psi \neq 0$;
	\item[-] $(2m-k) \cdot h <t$ and $\varphi \neq 0$;
	\item[-] $(2m-k) \cdot h >t$ and $\psi \neq 0$.
\end{itemize}
\end{dfn}
The main result concerning stable Witten triples is the following:
\begin{prop}
Let $V$ be a surface, and choose a class $m \in H^2(V,\nz )$. Fix an ample class $h\in H^2(V,\nz)$ and a real number $t \in \nr$. Then there exists a fine moduli space parametrizing $t$-stable Witten triples of class $m$ on $(V,h)$.
\end{prop}
\begin{proof}
\cite[Thm1.12]{d2}.
\end{proof}
Denote the moduli space of $t$-stable Witten triples by $M^m_{h,t}$. There is a natural morphism
\[
\mu_t: M^m_{h,t} \longrightarrow H^0(\kv ),
\]
which maps the class of a triple $(\cL , \varphi, \psi )$ to the holomorphic 2-form $\psi \circ \varphi$. Recall that above we introduced the addition map
\[
\alpha: \Hilb^m_V \times_{\Pic^m_V} \Hilb^{k-m}_V \to | \kv |.
\]

\begin{prop}
Let $V$ be a surface, and choose a class $m \in H^2(V,\nz )$. Fix an ample class $h\in H^2(V,\nz)$ and a real number $t \in \nr$. For every holomorphic 2-form $\eta \in H^0(\kv)$, there exists a natural isomorphism
\[
\mu_t^{-1}(\eta ) \stackrel{\cong}{\longrightarrow} \left\{
\begin{array}{cl}
\Hilb^m_V & \text{if $\eta=0$ and $(2m-k) \cdot h<t$,}\\
\Hilb^{k-m}_V & \text{if $\eta=0$ and $(2m-k) \cdot h>t$,}\\
\alpha^{-1} [\eta] & \text{if $\eta \neq 0$.}
\end{array}\right.
\]
\end{prop}
\begin{proof}
\cite[Thm.3.3]{d2}.
\end{proof}
\subsection{Some evidence}	
In this subsection we will collect further evidence for Conjecture \ref{conj:main}. We begin with two general facts:
\begin{itemize}
	\item[-] reduction to the minimal case
	\item[-] wall crossing formulas and consequences
\end{itemize}
Thereafter, we proceed with a case by case analysis, organized according to the different Kodaira dimensions.
\subsubsection{Reduction to the minimal case}
With  help of the blow-up formulas for the Seiberg-Witten invariants \cite[Thm.2.2]{os} and the Poincar\'e invariants (Thm.~\ref{thm:blowup}), we reduce the proof of Conj.~\ref{conj:main} to the minimal case:
\begin{thm}
Let $\sigma: \hat{V} \to V$ be the blow-up of a point $p \in V$. Denote by $\oo$ or $\oov$ the canonical orientation data on $V$, and by $\oohat$ or $\oovhat$ the canonical orientation data on $\hat V$. If $p_g(V)=0$ and
\[
P^\pm_V(m)=\swwv{\pm}{\mathfrak{c}_m}\ \ \forall m\in H^2(V,\mathbb{Z}),
\]
then
\[
P^\pm_{\hat V}(\hat{m})=\swwvhat{\pm}{\mathfrak{c}_{\hat{m}}}\ \ \forall \hat{m}\in H^2(\hat{V},\mathbb{Z}).
\]
If $p_g(V)>0$ and
\[
P^+_V(m)=P^-_V(m)=\swv{\mathfrak{c}_m}\ \ \forall m\in H^2(V,\mathbb{Z}),
\]
then
\[
P^+_{\hat V}(\hat{m})=P^-_{\hat V}(\hat{m})=\swvhat{\mathfrak{c}_{\hat m}}\ \ \forall \hat{m}\in H^2(\hat{V},\mathbb{Z}).
\]
\end{thm}
\begin{proof}
Let $E$ be the exceptional curve and set $e:=c_1(\ov(E))$. Let $m \in H^2(V,\nz)$ be a cohomology class and let $l$ be an integer. Then Theorem 2.2 of Ozsv\'ath-Szab\'o, restated in our teminology, says that
\[
\swvhat{\mathfrak{c}_{\sigma^*(m)+l \cdot e}} = \tau_{\leq m(m-k) - 2\left( \dfrac{l}{2}\right)} \swv{\mathfrak{c}_m}
\]
when $b^+(V)>1$, and 
\[
\swwvhat{\pm}{\mathfrak{c}_{\sigma^*(m)+l \cdot e}} = \tau_{\leq m(m-k) - 2\left( \dfrac{l}{2}\right)} \swwv{\pm}{\mathfrak{c}_m}
\]
when $b_+(V)=1$. Hence our theorem is a consequence of Thm.~\ref{thm:blowup}
\end{proof}
\subsubsection{Wall crossing formulas and consequences}
\begin{prop}\label{prop:prwc}
Let $V$ be a surface with $p_g(V)=0$, and denote by  $\oov$ the canonical orientation data. Then
\[
P^+_V(m)-P^-_V(m) = \swwv{+}{\mathfrak{c}_m}-\swwv{-}{\mathfrak{c}_m}
\]
for all $m \in H^2(V,\nz )$.
\end{prop}
\begin{proof}
This is a consequence of the respective wall crossing formulas Thm.~\ref{thm:wcfpoinc} and Thm.~\ref{thm:wcfsw}.
\end{proof}
\begin{cor}\label{cor:prwc}
Let $V$ be a surface with $p_g(V)=0$, denote by  $\oov$ the canonical orientation data, and fix an element $m \in H^2(V, \nz )$. If $\Hilb^m_V$ or $\Hilb^{k-m}_V$ is empty, then
\[
P^\pm_V(m)= \swwv{\pm}{\mathfrak{c}_m}.
\]
\end{cor}
\begin{proof}
When $\Hilb^m_V =\emptyset$, then the Kobayashi-Hitchin correspondence [Thm.~\ref{thm:kh}] yields
\[
\swwv{+}{\mathfrak{c}_m}=P^+_V(m)=0.
\]
Analogously, when $\Hilb^{k-m}_V =\emptyset$, we find
\[
\swwv{-}{\mathfrak{c}_m}=P^-_V(m)=0.
\]
Therefore, our claim is a consequence of Prop.~\ref{prop:prwc}.
\end{proof}
\subsubsection{Case by case analysis}
\begin{prop}
Let $V$ be a surface with $kod(V)=-\infty$, and denote by $\oov$ the canonical orientation data. Then, for any $m \in H^2(V,\mathbb{Z})$,
\[
P^{\pm}_V(m)=\swwv{\pm}{\mathfrak{c}_m}.
\]
\end{prop}
\begin{proof}
As we have seen earlier, the presence of a smooth rational curve $C$ on $V$ with $C^2\geq 0$ implies that for any $m \in H^2(V,\nz )$ one of the Hilbert schemes $\Hilb^m_V$ or $\Hilb^{k-m}_V$ is empty. Therefore, our claim is a consequence of Cor.~\ref{cor:prwc}.
\end{proof}
\begin{prop}
Let $V$ be a surface with $kod(V)=0$, and denote by $\oo$ or $\oov$ the canonical orientation data. If $V$ is a blow-up of a $K3$ surface or of a torus, then
\[
P^\pm_V(m)= \swv{\mathfrak{c}_m}
\]
for any $m \in H^2(V,\mathbb{Z})$. If $V$ is a blow-up of a bielliptic surface or of an Enriques surface, then
\[
P^\pm_V(m)= \swwv{\pm}{\mathfrak{c}_m}
\]
for any $m \in H^2(V,\mathbb{Z})$.
\end{prop}
\begin{proof}
When $V$ is an Enriques surface, then for any $m \in H^2(V,\nz)$ one of the Hilbert schemes $\Hilb^m_V$ or $\Hilb^{k-m}_V$ is empty. When $V$ is a $K3$-surface, an abelian variety, or bielliptic, then $\Hilb^k_V$ consists of one smooth point, the divisor $0$. Hence also in this case one of the Hilbert schemes $\Hilb^m_V$ or $\Hilb^{k-m}_V$ is empty unless $m=0$. Therefore our claim follows from Cor.~\ref{cor:prwc}.
\end{proof}
Next we consider properly elliptic surfaces.
\begin{prop}\label{prop:prell}
Let $V$ be a properly elliptic surface, and denote by $\oo$ or $\oov$ the canonical orientation data. If $p_g(V)>0$, then
\[
P^\pm_V(m)= \swv{\mathfrak{c}_m}
\]
for any $m \in H^2(V,\mathbb{Z})$. If $p_g(V)=0$, then
\[
P^\pm_V(m)= \swwv{\pm}{\mathfrak{c}_m}
\]
for any $m \in H^2(V,\mathbb{Z})$.
\end{prop}
\begin{proof}
It suffices to show that the relevant Hilbert schemes are smooth. The claimed equality follows then from Thm.~\ref{thm:andreimbd}. Without loss of generality, we may assume that $V$ is minimal.

Fix a class $m \in H^2(V,\nz)$ with $m(m-k)\geq 0$. Suppose first that one of the the Hilbert schemes $\Hilb^m_V$ or $\Hilb^{k-m}_V$ is empty. When $p_g(V)=0$, our claim is a consequence of Cor.~\ref{cor:prwc}. In the case $p_g(V)>0$, Thm.~\ref{thm:compobth} yields $P^+_V(m)=P^-_V(m)=0$, while the Kobayashi-Hitchin correspondence implies $\swv{\mathfrak{c}_m}=0$.

Assume now that both Hilbert schemes are non-empty. But then we have $m^2=\langle m,[F]\rangle=0$ and both Hilbert schemes are smooth by Lemma \ref{lem:fib}.
\end{proof}
\begin{cor}
Let $\pi : V \to C$ be an elliptic fibration over a curve of genus $g$. Let $F$ be a general fiber, and let $m_1F_1,\ldots m_rF_r$ be the multiple fibres of $\pi$. Denote by $\oo$ or $\oov$ the canonical orientation data, and fix a class $m \in H^2(V,\nz )$.
\begin{itemize}
	\item[i)] If $V$ is minimal with $p_g(V) >0$, then $\swv{\mathfrak{c}_m}=0$ unless\\
 $m^2=\langle m, [F]\rangle =0$.
	\item[ii)] If $V$ is minimal with $p_g(V) =0$, then $\swwv{+}{\mathfrak{c}_m}=0$ or\\
 $ \swwv{-}{\mathfrak{c}_m}=0$ unless $m^2=\langle m, [F]\rangle =0$.
	\item[iii)] If $m^2=\langle m, [F]\rangle =0$ and $p_g(V) >0$, then
\[
\swv{\mathfrak{c}_m}= \sum_{\dfrac{d[F]+\sum a_i[F_i] =PD(m)}{0\leq a_i <m_i}} (-1)^d 
\begin{pmatrix}
2g-2 + \chi (\ov )\\
d
\end{pmatrix}.
\]
	\item[iv)] If $m^2=\langle m, [F]\rangle =0$ and $p_g(V) =0$, then
\[
\swwv{+}{\mathfrak{c}_m}=
\sum_{\dfrac{d[F]+\sum a_i[F_i] =PD(m)}{0\leq a_i <m_i}} (-1)^d 
\begin{pmatrix}
2g-2 + \chi (\ov )\\
d
\end{pmatrix},
\]
and
\[
\swwv{-}{\mathfrak{c}_m}=
\sum_{\dfrac{d[F]+\sum a_i[F_i] =PD(k-m)}{0\leq a_i <m_i}} (-1)^{\chi (\ov)+d} 
\begin{pmatrix}
2g-2 + \chi (\ov )\\
d
\end{pmatrix}.
\]
\end{itemize}
\end{cor}
\begin{proof}
This follows from Prop.~\ref{prop:pifib} and Prop.~\ref{prop:prell}.
\end{proof}
\begin{rem}
Note that our formulas for the Seiberg-Witten invariants do not agree with the formulas given by Brussee \cite[Prop.~42]{br} and Friedman-Morgan \cite[Prop.~4.4]{fm}. There are two problems with the formulas in \cite{br} and \cite{fm}. The first, conceptual problem is the missing justification of the computation of the Seiberg-Witten invariants in terms of intersection theory on Hilbert schemes. The second problem is of a calculatory nature: in general, the relevant Hilbert schemes are not connected (see example 5.13), but Brussee and Friedman-Morgan find just one of their connected components.
\end{rem}
Let $V$ be a minimal surface of general type. When $q(V)=0$, then all Hilbert schemes are linear systems and in particular smooth. Hence, by applying Thm.~\ref{thm:andreimbd}, we find that the Seiberg-Witten- and the Poincar\'e invariants coincide.

Suppose now that $q(V)>0$. Note that such a surface has $p_g(V)>0$ since $\chi(\ov) >0$ for any surface of general type. Then $V$ has exactly two basic classes, namely $0$ and $k$, as we have shown in the proof of Prop.~\ref{prop:simp}. Likewise we have
\[
\swv{\mathfrak{c}_m}=0
\]
unless $m=0$ or $m=k$ \cite[p.789]{w}. Furthermore we know that
\begin{eqnarray*}
\swv{\mathfrak{c}_{can}}&=&1,\\
\swv{\mathfrak{c}_k}&=&(-1)^{\chi(\ov)},\\
P^+_V(0)&=&1,
\end{eqnarray*}
and
\[
P^-_V(k)= (-1)^{\chi(\ov)}.
\]
Hence, in order to give a case by case proof of Conj.~\ref{conj:main}, it remains to show that $P^-_V(0)= (-1)^{\chi(\ov)} P^+_V(k)=1$.

The results of section 6.2 show that it suffices to prove the following:\\[5mm]
{\bf Assertion:} Let $V$ be a minimal surface of general type with $p_g(V)>0$ and $q(V)>0$. Then
\[
\deg\, \lb \Hilb^k_V \rb = (-1)^{ \chi (\ov )}.
\]

\end{document}